\documentclass[11pt, reqno]{amsart}
\title[Yangians, quantum loop algebras and difference equations]
{Yangians, quantum loop algebras and abelian difference equations}
\author[S. Gautam]{Sachin Gautam}
\address{Department of Mathematics,
Columbia University,
2990 Broadway,
New York, NY 10027}
\email{sachin@math.columbia.edu}
\author[V. Toledano Laredo]{Valerio Toledano Laredo}
\address{Department of Mathematics,
Northeastern University,
567 Lake Hall,
360 Huntington Avenue,
Boston, MA 02115}
\dedicatory{To Andrei Zelevinsky (1953--2013). Advisor, mentor, colleague and friend.}
\thanks
{VTL was supported in part by NSF grants DMS--0854792, DMS--1206305 and PHY--1066293}
\email{V.ToledanoLaredo@neu.edu}

\usepackage{amsmath,amssymb}

\usepackage[usenames,dvipsnames]{color}
\usepackage{hyperref}
\hypersetup{
    colorlinks=true,       			% false: boxed links; true: colored links
    linkcolor=blue, %Cerulean,          	% color of internal links
    citecolor=blue,        			% color of links to bibliography
    filecolor=blue,      				% color of file links
    urlcolor=blue           			% color of external links
}

\usepackage[all]{xy}
\xyoption{arc}
\newtheorem*{thm}{Theorem}
\newtheorem*{prop}{Proposition}
\newtheorem*{lem}{Lemma}

\newenvironment{pf}{\paragraph{{\sc Proof}}}{\qed\par\medskip}
\theoremstyle{definition}

\newtheorem*{rem}{Remark}

\numberwithin{equation}{section}

\numberwithin{figure}{section}

\newcommand{\spec}{\sigma}

\newcommand{\cint}[2]{\frac{1}{2\pi\iota}\int_C\, #1\, d#2}

\newcommand{\Fh}[1]{\mathsf{\Gamma}_{#1}}%{\mathcal{F}_{#1}}

\newcommand{\id}{\mathbf{1}}

\newcommand{\lp}{\left(}
\newcommand{\rp}{\right)}

\newcommand{\ga}{\mathfrak{a}}

\newcommand{\g}{\mathfrak{g}}
\newcommand{\h}{\mathfrak{h}}
\newcommand{\gl}{\mathfrak{gl}}
\newcommand{\Lsl}{\mathfrak{sl}}

\newcommand{\Sym}{\mathfrak{S}}

\newcommand{\bfA}{\mathbf{A}}

\newcommand{\A}{\mathcal{A}}

\newcommand{\LL}{\mathcal{L}}

\newcommand{\calO}{\mathcal{O}}
\newcommand{\PP}{\mathcal{P}}
\newcommand{\QQ}{\mathcal{Q}}

\newcommand{\RR}{\mathcal{R}}

\newcommand{\V}{\mathcal{V}}
\newcommand{\X}{\mathcal X}

\newcommand{\C}{\mathbb{C}}
\newcommand{\nC}{\mathbb{C}^{\times}}
\newcommand{\N}{\mathbb{N}}
\newcommand{\Q}{\mathbb{Q}}
\newcommand{\R}{\mathbb{R}}
\newcommand{\Z}{\mathbb{Z}}

\newcommand{\qchar}{\chi_q^U}
\newcommand{\ychar}{\chi_q^Y}

\newcommand{\fv}{{\mathbf v}}
\newcommand{\hit}{\operatorname{ht}}

\newcommand {\Exp}[1]{e^{2\pi\iota #1}}
\newcommand {\wh}[1]{\widehat{#1}}
\newcommand {\ol}[1]{\overline{#1}}

\newcommand{\End}{\operatorname{End}}
\newcommand{\Hom}{\operatorname{Hom}}
\newcommand{\ad}{\operatorname{ad}}
\newcommand{\Ad}{\operatorname{Ad}}

\newcommand {\aand}{\qquad\text{and}\qquad}
\renewcommand{\dim}{\operatorname{dim}}

\newcommand {\Rep}{\operatorname{Rep}}
\newcommand {\ffd}{\operatorname{fd}}
\newcommand {\fd}{finite--dimensional }
\newcommand {\lhs}{left--hand side }
\newcommand {\rhs}{right--hand side }
\newcommand {\wrt}{with respect to }
\newcommand {\hw}{highest weight }
\newcommand {\irr}{irreducible }

\newcommand{\Rloop}{\Rep_{\ffd}(\qloop)}

\newcommand{\Rlsub}{\Rep_{\ffd}^{\Omega}(\qloop)}

\newcommand{\Ryang}{\Rep_{\ffd}(\Yhg)}

\newcommand{\Rync}{\Rep_{\ffd}^{\scriptscriptstyle{\operatorname{NC}}}(\Yhg)}
\newcommand{\Rysub}{\Rep_{\ffd}^{\Pi}(\Yhg)}

\newcommand{\ds}{\displaystyle}
\newcommand{\wt}[1]{\widetilde{#1}}

\newcommand{\pp}{2\pi\iota}

\newcommand {\qloop}{U_q(L\g)}

\newcommand {\qloopsl}[1]{U_{q}(L\Lsl_{#1})}
\newcommand {\qloopgl}[1]{U_{q}(L\gl_{#1})}

\newcommand {\Yhg}{Y_\hbar(\g)}

\newcommand{\Yhsl}[1]{Y_{\hbar}(\sl_{#1})}

\newcommand {\isom}{\stackrel{\sim}{\rightarrow}}

\newcommand{\qbin}[3]{\left[\begin{array}{c} #1 \\ #2\end{array}\right]_{#3}}

\renewcommand {\gl}{\mathfrak{gl}}
\renewcommand {\sl}{\mathfrak{sl}}

\newcommand {\sfA}{{\mathsf A}}

\newcommand {\KM}{Kac--Moody }

\newcommand {\eg}{{\it e.g., }}

\newcommand{\bfI}{{\mathbf I}}

\usepackage{esint}
\newcommand{\nZ}{\mathbb{Z}_{\neq 0}}
\newcommand{\nN}{\mathbb{Z}_{>0}}

\newcommand{\dwtY}{\mathcal{P}_+^Y}
\newcommand{\dwtU}{\mathcal{P}_+^U}
\renewcommand{\Rloop}{\calO_{\scriptscriptstyle{\operatorname{int}}}(\qloop)}
\renewcommand{\Rlsub}{\calO_{\scriptscriptstyle{\operatorname{int}}}^{\Omega}(\qloop)}

\renewcommand{\Ryang}{\calO_{\scriptscriptstyle{\operatorname{int}}}(\Yhg)}

\renewcommand{\Rync}{\calO_{\scriptscriptstyle{\operatorname{int}}}^{\scriptscriptstyle{\operatorname{NC}}}(\Yhg)}
\renewcommand{\Rysub}{\calO_{\scriptscriptstyle{\operatorname{int}}}^{\Pi}(\Yhg)}

\newcommand{\Rfloop}{\Rep_{\ffd}(\qloop)}
\newcommand{\Rfyang}{\Rep_{\ffd}(\Yhg)}

\newcommand{\Rfysub}{\Rep_{\ffd}^\Pi(\Yhg)}
\newcommand{\Rflsub}{\Rep_{\ffd}^{\Omega}(\qloop)}

\renewcommand {\Re}{\operatorname{Re}}
\renewcommand {\Im}{\operatorname{Im}}

\newcommand {\calF}{\mathcal F}
\newcommand {\calP}{\mathcal P}
\newcommand {\calZ}{\mathcal Z}

\newcommand {\rightleft}{\pm\Re(u)\gg 0}
\newcommand {\relto}{relative to }

\newcommand {\IP}{\mathbb P}
\newcommand {\Y}{\mathcal Y}
\newcommand {\QL}{$\mathcal{QL}$}
\newcommand {\veps}{\varepsilon}

\newcommand {\tord}{toro\"{\i}dal }

\newcommand {\calC}{\mathcal{C}}
\newcommand {\ointC}{\oint_{\calC}}
\newcommand {\ointCi}[2]{\oint_{\calC_{#1}^{#2}}}
\newcommand {\intC}{\int_{\calC}}
\newcommand {\intCi}[2]{\int_{\calC_{#1}^{#2}}}
\newcommand {\eps}{\epsilon}
\newcommand {\calG}{\mathcal G}

\newcommand {\Uqg}{U_q\g}

\newcommand {\Gh}[1]{\mathbf{\daleth}_{#1}} 

\newcommand {\IH}{{\mathbb H}}
\newcommand {\FH}{\mathfrak H}

\newcommand {\bfc}{{\mathbf c}}

\newcommand {\wYhg}{\wt{Y}_\hbar(\g)}
\newcommand {\wqloop}{\wt{U}_q(L\g)}
\newcommand {\EExp}{\mathcal{E}\negthinspace xp}

\newcommand {\Hecked}{\mathfrak H}
\newcommand {\Heckea}{\mathbb H}

\newcommand {\aalpha}{{\mathbf a}}

\begin{document}

\begin{abstract}
Let $\g$ be a complex, semisimple Lie algebra, and $\Yhg$
and $\qloop$ the Yangian and quantum loop algebra of $\g$.
Assuming that $\hbar$ is not a rational number and that $q=
e^{\pi i\hbar}$, we construct an equivalence between the \fd
representations of $\qloop$ and an explicit subcategory of
those of $\Yhg$ defined by choosing a branch of the logarithm.
This equivalence is governed by the monodromy of the abelian,
additive difference equations defined by the commuting fields
of $\Yhg$. Our results are compatible with $q$--characters,
and apply more generally to a symmetrisable \KM algebra
$\g$, in particular to affine Yangians and quantum \tord algebras.
In this generality, they yield an equivalence between the representations
of $\Yhg$ and $\qloop$ whose restriction to $\g$ and $\Uqg$
respectively are integrable and in category $\calO$.
\end{abstract}
\maketitle

\setcounter{tocdepth}{1}
\tableofcontents

\section{Introduction}\label{sec: intro}
%===============

\subsection{} % what we do in this paper
%--------------

Let $\g$ be a complex semisimple Lie algebra, and $\Yhg$ and
$\qloop$ the Yangian and quantum loop algebra of $\g$. Recall
that the latter are deformations of the enveloping algebras of the
current Lie algebra $\g[s]$ of $\g$ and its loop algebra $\g[z,z^{
-1}]$ respectively. We shall assume throughout that $q$ is not a
root of unity, and that $e^{\pi\iota\hbar}=q$.

In this paper, we solve the long--standing problem of relating the
\fd representations of $\Yhg$ and $\qloop$, thus giving a precise
formulation to the widely held belief that these quantum groups
have the same representation theory (see, \eg \cite[p. 283]{varagnolo-yangian}).
Our solution is entirely explicit, and expressed in terms of Drinfeld's
loop generators. Given a \fd representation $V$ of $\Yhg$, we
construct an action of $\qloop$ on $V$ by contour integral
formulae, and viceversa. In particular, this gives rise to an exact,
faithful functor
\[\Fh{}:\Rfyang\longrightarrow\Rfloop\]
which commutes with the forgetful functor to vector spaces. The
functor $\Fh{}$ may be thought of as an exponential cover, in that
it maps simple objects to simple objects, exponentiates the roots
of their Drinfeld polynomials, and restricts to an equivalence on the
subcategory of $\Rfyang$ obtained by requiring these roots to lie
in a domain $\Pi\subset\C$ mapping bijectively to $\C^\times$ via
the map $u\to\Exp{u}$.

\subsection{}
%--------------

Our results apply more generally to an arbitrary symmetrisable \KM
algebra $\g$, in particular to affine Yangians and quantum \tord
algebras, and are proved in that generality in the main body of the
paper. When $\g$ is infinite--dimensional, the relevant categories
consist of representations of $\Yhg,\qloop$ whose restriction to $\g$
and $\Uqg$ respectively is integrable and in category $\calO$. For
simplicity of exposition, however, we assume that $\g$ is \fd for the
rest of this introduction.

\subsection{}
%--------------

Yangians and quantum loop algebras have been known to be
closely related since their introduction by Drinfeld in the mid 80s
\cite{drinfeld-qybe,drinfeld-quantum-groups}. In particular

\begin{enumerate}
\item[(1)] Drinfeld proved that, just as the loop algebra $\g[z,z^{-1}]$
degenerates to the current algebra $\g[s]$ by setting $z=\exp
(\epsilon s)$ and letting $\epsilon\to 0$, the quantum loop algebra
degenerates to the Yangian. Specifically, $\Yhg$ is the graded
algebra associated to the filtration of $\qloop$ by powers of the
ideal of $z=1$ (see \cite{drinfeld-quantum-groups}, and \cite
{guay-degeneration} for a complete proof).
\item[(2)] The highest weight classification of simple modules
due to Drinfeld and Chari--Pressley \cite{drinfeld-yangian-qaffine,
chari-pressley-qaffine} yields a surjection from the parameters
of \fd irreducible representations of $\Yhg$ to those of $\qloop$.
This (set--theoretic) map $\EExp$ is given by exponentiating the 
roots of Drinfeld polynomials.
\item[(3)] Stronger results can be obtained when $\g$ is of type {\sf ADE},
by using Nakajima quiver varieties, and the geometric action of
$\qloop$ on their equivariant $K$--theory \cite{nakajima-qaffine}.
The analogous action of the Yangian $\Yhg$ in cohomology due
to Varagnolo \cite{varagnolo-yangian} allows to relate irreducible
representations of $\qloop$ and $\Yhg$ via the Chern character,
and shows in particular that the map $\EExp$ preserves dimensions.
\end{enumerate}

\subsection{}
%--------------

While (1) does not a priori suggest the existence of a link between
the representations of $\qloop$ and $\Yhg$ (there is after all no
relation between the representations of a filtered algebra and of
its associated graded in general), (2) and (3) clearly do, as does
the fact that, in type $\sfA$, evaluation representations of $\Yhg$
can be explicitly deformed to representations of $\qloop$. 
Such a link is also implicit in the conjecture of
Frenkel--Reshetikhin according to which the monodromy of the
difference analogue of the Knizhnik--Zamolodchikov ($q$KZ for short) 
equations defined by the $R$--matrix
of $\Yhg$ is described by the $R$--matrix of $\qloop$
\cite{frenkel-reshetikhin}.

The belief that \fd representations of $\Yhg$ and $\qloop$ should
be related is further corroborated by the fact that, in the analogous
case of the degenerate and affine Hecke algebras $\Hecked,\Heckea$
of a Weyl group, a complete understanding of the relation between
\fd representations of $\Hecked$ and $\Heckea$ was obtained by
Lusztig \cite{lusztig-graded-affine}, and can be recovered by realising
$\Hecked,\Heckea$ as the equivariant cohomology and $K$--theory
of the corresponding Steinberg variety \cite{ginzburg-geometric-methods}.

\subsection{}
%--------------

Despite the above evidence, however, 
no precise, functorial relation between the categories $\Rfyang$
and $\Rfloop$ was known so far, even conjecturally\footnote{The
Chern character in the geometric realisations (3) above cannot be
considered a functor, since it is not known how to realise arbitrary
\fd representations of $\Yhg$ of $\qloop$ inside the equivariant
cohomology or $K$--theory of a quiver variety.}. One notable
exception is the works \cite{JKOS,etingofmoura-KL,
etingof-schiffmann-hw} which relate 
$\qloopsl{n}$ and Felder's elliptic quantum group $E_{\tau,
\gamma}(\gl_n)$, and construct
in particular a fully faithful functor $\calF$ from $\Rep_{\ffd}(\qloopsl{n})$
to $\Rep_{\ffd}(E_{\tau,\gamma}(\gl_n))$ \cite{etingofmoura-KL}.
A similar construction should give rise to a functor $F:\Rep_{\ffd}
(\Yhsl{n})\to\Rep_{\ffd}(\qloopgl{n})$. These constructions
depend, however, on the RTT presentation of $\qloopsl{n}$,
which is less convenient for other classical Lie algebras, and
currently unavailable for exceptional ones. Moreover,
even in type $\sfA$, the essential surjectivity of the functors
$\calF$ and $F$ seems difficult to establish, and may in fact
fail for numerical values of $\tau,q$.

\subsection{}
%--------------

Our results are in a sense optimal, short of proving that
$\qloop$ and $\Yhg$ are isomorphic. An elementary
analogy may be seen by looking at the exponential covering
$\C\longrightarrow\C^\times$. The latter is a local isomorphism, and possesses a
global section for any determination of the logarithm.

The analogous local statement for $\Yhg$ and $\qloop$ was
obtained in our previous work. Indeed, we proved in \cite{sachin-valerio-1}
that, when the deformation parameters $\hbar,q$ are formal,
the quantum loop algebra and Yangian are isomorphic as
algebras after appropriate completions, specifically \wrt the
natural loop grading on $\Yhg$ and the powers of the ideal
of $z=1$ in $\qloop$. The isomorphism 
\[\Phi:\wh{\qloop}\longrightarrow\wh{\Yhg}\]
constructed in \cite{sachin-valerio-1} deforms the change of
variables $\g[z,z^{-1}]\to\g[[s]]$ given by $z=\Exp{s}$, and exponentiates
the roots of Drinfeld polynomials.

The framework of \cite{sachin-valerio-1} is, however, not suited
for the study of \fd representations since the isomorphism $\Phi$
is given by formal power series which do not converge on these 
representations, and do not allow for the numerical specialisation
of the deformation parameters.

\subsection{}\label{ssec: intro-main}
%--------------

Our first main result is the following

\begin{thm}
If $\hbar\in\C\setminus\Q$ and $q=e^{\pi i\hbar}$, there is an
exact, faithful functor
\[\Fh{}:\Rfyang\longrightarrow\Rfloop\]
which commutes with the forgetful functor to vector spaces,
maps simple objects to simple objects, and exponentiates
the roots of Drinfeld polynomials.\footnote{Strictly speaking,
the functor $\Fh{}$ is only defined on the dense subcategory
of {\it non--congruent} representations of $\Yhg$ introduced
in \S \ref{sec: functor}. See Section \ref{ss:1st main thm} for
a precise statement of Theorem \ref{ssec: intro-main}. We
shall gloss over this point in the introduction.}
\end{thm}

\subsection{}
%--------------

Our second main result asserts the existence of an inverse to
$\Fh{}$ for each determination of the logarithm. Specifically, 
fix a subset $\Pi\subset\C$ such that
\[(\Pi-\Pi)\cap\Z=\{0\}\aand
\Pi\pm \frac{\hbar}{2}\subset \Pi\]
The map $u\mapsto \Exp{u}$ is then a bijection between $\Pi$
and a subset $\Omega\subset\nC$ stable under multiplication
by $q^{\pm 1}$. To such a choice, we attach categories
\[\Rfysub\aand\Rflsub\]
consisting of representations such that the roots of their Drinfeld
polynomials lie in $\Pi$ and $\Omega$ respectively. 

\begin{thm}\label{thm: intro-main}\hfill
\begin{enumerate}
\item[(1)] The functor $\Fh{}$ restricts to a functor
\[\Fh{\Pi}:\Rfysub\longrightarrow\Rflsub\]
\item[(2)] The functor $\Fh{\Pi}$ admits an inverse functor
\[\Gh{\Pi} : \Rflsub\longrightarrow \Rfysub\]
\end{enumerate}
In particular, if $\Pi$ is a fundamental domain for $u\to u+1$, $\Fh
{}$ restricts to an equivalence between $\Rfysub$ and $\Rfloop$.
\end{thm}

\subsection{}\label{ssec: intro-construction}
%--------------

The novel idea underlying our construction is that each of the
commuting fields of $\Yhg$ defines a {\it difference 
equation} on any \fd representation of $\Yhg$. The functor
$\Fh{}$ is governed by the {\it monodromy} of this equation,
an averaging which yields a canonical means of replacing the
rational matrix coefficients of $\Yhg$ by their trigonometric
counterparts for $\qloop$ (for the reader's convenience, a
self--contained account of Birkhoff's theory of difference equations is given in \S \ref{sec: difference}).
The {\it abelian} nature of these equations is essential in two
respects. On the one hand, it leads to explicit formulae
for the action of $\qloop$ on $\Yhg$--modules
in terms of $\Gamma$--functions.
On the other, it guarantees that the corresponding Riemann--Hilbert
problem is always solvable, and therefore that $\Fh{}$ is essentially surjective.

\subsection{} \label{ss:monodromy trick}
%--------------

To illustrate this construction, it is useful to compare the simplest
matrix coefficients of $\Yhg$ and $\qloop$ when $\g=\sl_2$. Let
$V(a),\V(\aalpha)$ be the representations of $\Yhsl{2},\qloopsl{2}$
obtained by evaluating the vector representations of $\sl_2$ and
$U_q\sl_2$ at $a\in\C$ and $\aalpha\in\C^\times$ respectively. If
\[\xi(u) = 1+\hbar\sum_{r\geq 0}\xi_r u^{-u-1} 
\aand
\Psi(z)^{\pm} = \sum_{k\geq 0} \Psi^{\pm}_{\pm k}z^{\mp k}
\]
are the generating functions of the commuting loop generators
of $\Yhsl{2}$ and $\qloopsl{2}$ respectively, the action of $\xi(u),
\Psi(z)^\pm$ on the highest weight vectors $\omega\in V(a)$,
$\Omega\in\V(\aalpha)$ is given respectively by
\[\xi(u)\,\omega=\frac{u+\hbar-a}{u-a}\,\omega
\aand
\Psi(z)^\pm\,\Omega=q^{-1}\frac{q^2z-\aalpha}{z-\aalpha}\,\Omega\]

These formulae are tantalizingly close, yet only seem to be related
by the mechanical, termwise exponentiation $z=\Exp{u}$, $\aalpha=
\Exp{a}$. 
Our construction in this case amounts to the observation that the
matrix coefficient $\psi(z)=(qz-q^{-1}\aalpha)/(z-\aalpha)$ is the monodromy
of the additive difference equation $f(u+1)=A(u)f(u)$ determined
by the matrix $A(u)=(u+\hbar-a)/(u-a)$. In other words, $\psi(z)$ is given
by the {\it averaging} 
\[\frac{qz-q^{-1}\aalpha}{z-\aalpha}=
\cdots\frac{u+1+\hbar-a}{u+1-a}\cdot\frac{u+\hbar-a}{u-a}
\cdot\frac{u-1+\hbar-a}{u-1-a}\cdots\]

\subsection{}
%--------------

This observation may in fact be used to define an action of
the commuting generators of $\qloop$, not only on the \hw
vector of a simple module, but on an arbitrary \fd representation
$V$ of $\Yhg$ for any semisimple $\g$, and in a way which
is consistent with taking $q$--characters.

Specifically, for any vertex $i$ of the Dynkin diagram $\bfI$
of $\g$, the generating function
\[\xi_i(u) = 1+\hbar\sum_{r\geq 0}\xi_{i,r} u^{-u-1}\in GL(V)[[u^{-1}]]\]
of the corresponding commuting elements of $\Yhg$ can be
shown to be the expansion at $u=\infty$ of an $\End(V)$--valued
rational function. One can therefore consider the additive difference equation
\begin{equation}\label{eq: intro-difference}
\phi_i(u+1) = \xi_i(u)\phi_i(u)
\end{equation}
determined by each $\xi_i(u)$, where the unknown $\phi_i$
is a meromorphic function of $u\in\C$ with values in $\End(V)$.
This equation admits two fundamental solutions $\phi_i^{\pm}(u)$
which are uniquely determined by the following requirements 
\begin{align*}
&\bullet\text{$\phi_i^{\pm}(u)$ is holomorphic and invertible for $\pm\Re(u)\gg0$}\\
&\bullet\text{$\phi_i^{\pm}(u)$ is asymptotic to $(1+O(u^{-1}))(\pm u)^{\hbar\xi_{i,0}}$ as
$u\to\infty$ with $\pm\Re(u)\gg0$}
\end{align*}
The solutions $\phi_i^+(u),\phi_i^-(u)$ are regularisations of
the formal infinite products
\[\xi_i(u)^{-1}\xi_i(u+1)^{-1}\cdots
\aand
\xi_i(u-1)\xi_i(u-2)\cdots\]
which symbolically solve \eqref{eq: intro-difference}.

Following Birkhoff \cite{birkhoff-difference}, one considers the
{\it monodromy matrix} of the equation, that is the $GL(V)$--valued
function defined by
\[S_i(u)=\phi_i^+(u)^{-1}\cdot\phi_i^-(u)\]
$S_i(u)$ is a regularisation of the infinite product
\[\cdots\xi_i(u+1)\xi_i(u)\xi_i(u-1)\cdots\]
and, by construction, is a 1--periodic function of $u$. The
asymptotics of $\phi_i^\pm$ further imply that $S_i$ is a rational
function of the variable $z=e^{2\pi\iota u}$, regular at $z=\infty,0$,
and such that
$S_i(\infty)=e^{\pi\iota\hbar\xi_{i,0}}=S_i(0)^{-1}$.

The action of the generating series
\begin{equation}\label{eq:intro-psi}
\Psi_i(z)^{\pm} = \sum_{k\geq 0} \Psi^{\pm}_{i,\pm k}z^{\mp k}
\end{equation}
of the commuting generators of $\qloop$ on the $\Yhg$--module
$V$ is then given by the Taylor series of the monodromy matrix
$S_i(z)$ at $z=\infty,0$. In the example described in \ref{ss:monodromy trick}, 
the (matrix elements of the) Birkhoff solutions are given by
\[\phi^+(u)=\frac{\Gamma(u+\hbar-a)}{\Gamma(u-a)}
\aand
\phi^-(u)=\frac{\Gamma(1-u+a)}{\Gamma(1-u-\hbar+a)}\]
where $\Gamma$ is Euler's Gamma function. The corresponding
monodromy is
\[S(z)=
\frac{\Gamma(u-a)}{\Gamma(u+\hbar-a)}\cdot
\frac{\Gamma(1-u+a)}{\Gamma(1-u-\hbar+a)}=
\frac{qz-q^{-1}\aalpha}{z-\aalpha}
\]
where we used the identity $\Gamma(u)\Gamma(1-u)=\pi/\sin(\pi u)$
\cite[\S 12.14]{whittaker-watson}.

\subsection{}\label{ss:forward}
%--------------

The difference equation \eqref{eq: intro-difference} can also be used
to define an action of the raising and lowering loop generators $\{\X^
{\pm}_{i,k}\}_{i\in\bfI,k\in \Z}$ of $\qloop$ on the \fd $\Yhg$--module
$V$ as follows. For each $i\in\bfI$ and $k\in\Z$, $\X^{\pm}_{i,k}$ acts
as the operator defined by the contour integral
\[\frac{c_i}{2\pi\iota}\int_{\calC_i^\pm}\Exp{ku}g_i^{\pm}(u)x_i^{\pm}(u)\,du\]
where
\begin{itemize}
\item $x_i^{\pm}(u)=\hbar\sum_{r\geq 0}x_{i,r}^\pm u^{-r-1}$ is
the generating function of the raising/lowering generators of
$\Yhg$, and acts on $V$ as a rational function of $u$.
\item The functions $g_i^+(u)$ and $g_i^-(u)$ are regularisations of
the products
\[\qquad\qquad\quad
\cdots \xi_i(u+2)\xi_i(u+1)
\aand
\xi_i(u-1)\xi_i(u-2)\cdots\]
respectively, and are given in terms of the fundamental solutions
$\phi_i^\pm(u)$ of the difference equation \eqref{eq: intro-difference}
by
\begin{equation}\label{eq: intro-g}
\qquad\quad\quad
g_i^+(u) = \phi_i^+(u+1)^{-1} \aand 
g_i^-(u) = \phi_i^-(u)
\end{equation}
\item The constant $c_i$ is given in terms of the symmetrising
integers $d_i$ of the Cartan matrix of $\g$ by
$c_i=\sqrt{d_i}\,\Gamma(d_i\hbar)$.
\item The contour $\calC_i^\pm$ contains all the poles of $x_i^
{\pm}(u)$, and none of their $\nZ$--translates.
\end{itemize}

The proof that the above operators satisfy the commutation relations
of $\qloop$ follows from contour integral manipulations reminiscent of
those used in Conformal Field Theory. One notable exception, however,
is that the integrands have singularites on the shifted diagonals $u-v=a$,
where $a$ is an integral multiple of $\hbar/2$.

\subsection{}\label{ssec: intro-i-construction}
%--------------

The construction of the inverse functor
\[\Gh{\Pi}:\Rflsub\longrightarrow\Rfysub\]
given by Theorem \ref{thm: intro-main} is obtained in a similar spirit
by solving the inverse monodromy problem corresponding to the
difference equation \eqref{eq: intro-difference}. It is worth stressing that, as for
its differential counterpart, this Riemann--Hilbert problem possesses
a solution for any monodromy matrix only because the equations
\eqref{eq: intro-difference} are abelian. This property is therefore
crucial in establishing the essential surjectivity of the functor $\Fh{}$.

Specifically, let $\V\in\Rfloop$, and consider the $\End(\V)$--valued
generating series $\Psi_i(z)^{\pm}$ given by \eqref{eq:intro-psi},
\[\X_i^{\pm}(z)^+ = \sum_{k\geq 0}\X^{\pm}_{i,k}z^{-k} \aand 
\X_i^{\pm}(z)^- = -\sum_{k\geq 1}\X^{\pm}_{i,-k}z^k\]
It is well--known that these are the expansions at $z=\infty,0$ of
rational functions $\Psi_i(z)$ and $\X^{\pm}_i(z)$ \cite{beck-kac,
hernandez-affinizations}. We show in \S \ref{sec: s-cat} that $\V$
is in $\Rflsub$ if, and only if the poles of $\Psi_i(z)^{\pm 1}$ and
$\X_i^\pm(z)$ are contained in $\Omega$.

The operators giving the action of $\Yhg$ on $\Gh{\Pi}(\V)=\V$
are then defined as follows.
\begin{enumerate}
\item[(1)] As explained in \S \ref{sec: difference}, there is a unique rational
function $\C\to GL(V)$, which we take to be $\xi_i(u)$, such that
\begin{itemize} 
\item 
$\xi_i(u)=1+\hbar d_i\alpha_i^\vee u^{-1}+\cdots$
\item
the connection matrix of the difference equation \eqref
{eq: intro-difference} defined by $\xi_i(u)$ is $S_i(z)=\Psi_i(z)$
\item
the poles of $\xi_i(u)^{\pm 1}$ lie in $\Pi$. 
\end{itemize}
\item[(2)] Let $\phi_i^{\pm}(u)$ be the canonical fundamental
solutions of \eqref{eq: intro-difference}, and define the
functions $g_i^{\pm}(u)$ by \eqref{eq: intro-g}. Then, for
any $i\in\bfI$ and $r\in\N$, the generators $x^\pm_{i,r}$
act on $\V$ as the operators
\[\frac{1}{c_i\hbar}\int_{\calC_i^\pm} v^{-r}
g_i^{\pm}(v)^{-1}\X^{\pm}_i\lp\Exp{v}\rp\, dv\]
where the constants $c_i$ are as in \S \ref{ss:forward},
and the contours $\calC_i^\pm$ enclose all the poles
of $\X^{\pm}_i\lp\Exp{v}\rp$ contained in $\Pi$, and none
of their $\nZ$--translates.
\end{enumerate}

\subsection{}
%--------------

The relevance of difference equations to Yangians and quantum loop
algebras is widely known. For example, the $R$--matrix of $\Yhg$ (resp.
$\qloop$) gives rise to the additive (resp. multiplicative) $q$KZ equations
satisfied by products of intertwiners \cite{frenkel-reshetikhin,smirnov}. In
a related vein, Baxter's $Q$--operator for $\Yhg$ (resp. $\qloop$) satisfies
the TQ relations, which are difference equations with step a multiple of
$\hbar$ (resp. a power of $q$) (see, \eg \cite{baxter,bazhanovetal,frenkel-hernandez}).
Moreover, the idea of averaging a rational or trigonometric solution of
the quantum Yang--Baxter equations to obtain a trigonometric or elliptic
solution has also been extensively used \cite{bazhanov85,bazhanov87,
faddeev11,TakhtajanZn}.

It is worth pointing out that the difference equations used in this paper
are not of $q$KZ or TQ type. Indeed, they are abelian and therefore
far simpler. Their relevance seem to be new, although they are very
closely related to the $\Fh{}$--factors which appear in quantum
cohomology \cite{iritani,okounkov-pandharipandeQDE}.

\subsection{} 
%--------------

The use of abelian difference equations to define a functor between
\fd representations of $\Yhg$ and $\qloop$ is the main discovery of
the present paper. This is a potentially far reaching principle, which
we have so far applied to two further problems
\begin{enumerate}
\item[(1)] The construction of a tensor structure on the functor
$\Fh{}$, which yields a Kazhdan--Lusztig equivalence of
meromorphic braided tensor categories between
$\Yhg$ and $\qloop$, when both are
endowed with the deformed Drinfeld coproduct and the
commutative part $R_0$ of the $R$--matrix. This tensor
structure arises from the monodromy of the abelian $q$KZ
equations defined by $R_0$ of $\Yhg$ \cite{sachin-valerio-III}.
\item[(2)] The constuction of a faithful functor $\mathsf{\Theta}$ between \fd
representations of $\qloop$ and of Felder's elliptic quantum
group \cite{sachin-valerio-V}. The functor $\mathsf{\Theta}$ is governed
by the monodromy of the difference equations defined by
the commuting fields of $\qloop$. This is a very promising
avenue, since elliptic quantum groups and their representations
are not very well understood outside of type ${\mathsf A}$
\cite{felder-icm}.
\end{enumerate}

We also plan to apply similar ideas to relate the affine Yangian
$Y_\hbar(\wh{\gl_1})$ appearing in the recent work of Maulik--Okounkov
\cite{maulik-okounkov} to the quantum toro\"{\i}dal algebra
$U_q(\wh{\wh{\gl_1}})$ \cite{feiginetal-toroidal1,miki-toroidal,
schiffmannvasserot-hall} (see also \cite{tsymbaliuk-gl1}, where
an isomorphism between appropriate completions of $U_q(\wh{\wh{\gl_1}})$
and $Y_\hbar(\wh{\gl_1})$ is obtained in the formal setting
by using the method of \cite{sachin-valerio-1}).

\subsection{}
%--------------

The picture described in this paper is consistent with Cherednik's philosophy
of obtaining equivalences of appropriate categories of \fd representations
of the affine Hecke algebra $\IH$ and degenerate affine Hecke algebra
$\FH$ of a finite Coxeter group through the monodromy of a flat
connection with values in $\FH$ \cite[\S 1.1--1.2]{cherednik-hecke}.
We note that, whereas Cherednik's affine KZ connection is differential
and not abelian, our equations are difference and abelian. It seems
an interesting problem to relate these two approaches.

\subsection{Outline of the paper}\label{ssec: intro-outline}
%---------------------------------------

In Section \ref{sec: yqla}, we review the definitions of the Yangian
$\Yhg$ and quantum loop algebra $\qloop$ associated to a
symmetrisable Kac--Moody algebra $\g$.

In Section \ref{sec: s-cat}, we consider the categories $\Ryang$
and $\Rloop$ of representations whose restriction to $\g$ and
$\Uqg$ respectively is integrable and in category $\calO$. We
review the classification of their simple objects, and prove that
the generating series of the loop generators of
$\Yhg$ and $\qloop$ act on objects of $\Ryang,\Rloop$ as rational functions.
We then define the subcategories
$\Rysub$ and $\Rlsub$ in terms of
the roots of Drinfeld polynomials. The main result of this section,
Theorem \ref{thm: first-main-theorem}, gives an alternative
characterisation of these categories in terms of the poles
of the fields $\xi_i(u),x_i^\pm(u)$ of $\Yhg$ (resp. $\Psi_i(z),
\X_i^\pm(z)$ of $\qloop$).

Section \ref{sec: difference} gives a self--contained account
of the theory of additive difference equations including canonical
fundamental solutions, monodromy, and a detailed discussion of
the inverse problem for abelian equations.

Section \ref{sec: functor} gives the definition of the functor $\Fh{}:
\Ryang\to\Rloop$ and the proof that the operators defined in \S \ref
{ss:forward} give an action of $\qloop$ on any 
category $\Ryang$ representation of $\Yhg$. We also prove the compatibility of
$\Fh{}$ with the shift automorphisms of $\Yhg$ and $\qloop$.

Section \ref{sec: inverse-functor} contains the contruction of
the inverse functor $\Gh{\Pi}:\Rlsub$ $\to\Rysub$, and the proof
that $\Gh{\Pi}\circ\Fh{\Pi}$ and $\Fh{\Pi}\circ\Gh{\Pi}$ are the
identity functors on $\Rysub$ and $\Rlsub$ respectively.

In Section \ref{sec: q-characters}, we review the definition of
the $q$--characters of $\Yhg$ and $\qloop$, and show that
the functor $\Fh{}$ is compatible with these.

\subsection{Acknowledgments} We are very grateful to Andrei
Okounkov and Igor Krichever for many stimulating discussions,
to David Hernandez for asking about the compatibility between
our construction and $q$--characters, and to Edward Frenkel
for his numerous comments on an earlier version of this paper.
Part of this project was completed while we were visiting ETH Z\"{u}rich,
while the second author was on a sabbatical leave at Columbia
University, and during a month long visit at the Aspen Center
for Physics. We are grateful to these institutions for their
support and wonderful working conditions, and to D. Calaque,
G. Felder, A. Okounkov, D. Freed, G. Moore, A. Neitzke
and H. Ooguri for their invitations.

\section{Yangians and quantum loop algebras}\label{sec: yqla}
%==================================

\subsection{}\label{ssec: kma}
%--------------

Let $\bfA = (a_{ij})_{i,j\in\bfI}$ be a symmetrisable generalised Cartan
matrix \cite{kac}. Thus, $a_{ii}=2$ for any $i\in\bfI$, $a_{ij}\in\Z_{\leq
0}$ for any $i\neq j\in \bfI$, and there exists a diagonal matrix $D$
with positive integer entries $\{d_i\}_{i\in\bfI}$ such that $D\bfA$ is 
symmetric. We assume that $(d_i)$ are relatively prime.

Let $(\h,\{\alpha_i\}_{i\in\bfI},\{\alpha_i^{\vee}\}_{i\in\bfI})$ be the unique
realization of $\bfA$. Thus, $\h$ is a complex vector space of dimension
$2|\bfI|-\text{rank}(\bfA)$, $\{\alpha_i\}_{i\in\bfI}\subset \h^*$ and $\{\alpha
_i^{\vee}\}_{i\in\bfI} \subset \h$ are linearly independent sets and, for any
$i,j\in\bfI$, $\alpha_j(\alpha_i^{\vee}) = a_{ij}$.

Let $\g$ be the Kac--Moody algebra associated to $\bfA$.

\subsection{The Yangian $\Yhg$}\label{ssec: yangian}
%========================

Let $\hbar\in\C$. The Yangian $\Yhg$ is the $\C$--algebra generated
by elements $\{x^{\pm}_{i,r},\xi_{i,r}\}_{i\in\bfI,r\in\N}$ and $h\in \h$,
subject to the following relations

\begin{enumerate}
\item[(Y0)] For any $i\in\bfI$,
\[\xi_{i,0} = d_i\alpha_i^{\vee}\]
\item[(Y1)] For any $i,j\in\bfI$, $r,s\in\N$ and $h,h'\in\h$
\[[\xi_{i,r}, \xi_{j,s}] = 0\qquad\qquad [\xi_{i,r},h]=0 \qquad\qquad [h,h']=0\]
\item[(Y2)] For $h\in\h$, $j\in\bfI$ and $s\in \N$
\[[h, x_{j,s}^{\pm}] = \pm \alpha_j(h) x_{j,s}^{\pm}\]
\item[(Y3)] For $i,j\in\bfI$ and $r,s\in\N$
\[[\xi_{i,r+1}, x^{\pm}_{j,s}] - [\xi_{i,r},x^{\pm}_{j,s+1}] =
\pm\hbar\frac{d_ia_{ij}}{2}(\xi_{i,r}x^{\pm}_{j,s} + x^{\pm}_{j,s}\xi_{i,r})\]
\item[(Y4)] For $i,j\in\bfI$ and $r,s\in \N$
\[
[x^{\pm}_{i,r+1}, x^{\pm}_{j,s}] - [x^{\pm}_{i,r},x^{\pm}_{j,s+1}]=
\pm\hbar\frac{d_ia_{ij}}{2}(x^{\pm}_{i,r}x^{\pm}_{j,s} + x^{\pm}_{j,s}x^{\pm}_{i,r})
\]
\item[(Y5)] For $i,j\in\bfI$ and $r,s\in \N$
\[[x^+_{i,r}, x^-_{j,s}] = \delta_{ij} \xi_{i,r+s}\]
\item[(Y6)] Let $i\not= j\in\bfI$ and set $m = 1-a_{ij}$. For any
$r_1,\cdots, r_m\in \N$ and $s\in \N$
\[\sum_{\pi\in\Sym_m}
\left[x^{\pm}_{i,r_{\pi(1)}},\left[x^{\pm}_{i,r_{\pi(2)}},\left[\cdots,
\left[x^{\pm}_{i,r_{\pi(m)}},x^{\pm}_{j,s}\right]\cdots\right]\right]\right.=0\]
\end{enumerate}

\subsection{}\label{ssec: formal-series-y}
%--------------

Assume now that $\hbar\neq 0$, and define $\xi_i(u),x^\pm_i(u)\in
\Yhg[[u^{-1}]]$ by
\[\xi_i(u)=1 + \hbar\sum_{r\geq 0} \xi_{i,r}u^{-r-1}\aand
x^{\pm}_i(u)=\hbar\sum_{r\geq 0} x_{i,r}^{\pm} u^{-r-1}\]

\begin{prop}\label{pr:Y fields}
The relations (Y1),(Y2)--(Y3),(Y4),(Y5),(Y6) are equivalent to the following
identities in $\Yhg[u,v;u^{-1},v^{-1}]]$
\begin{enumerate}
\item[($\Y$1)] For any $i,j\in\bfI$ and $h,h'\in\h$,
\[[\xi_i(u), \xi_j(v)]=0\qquad\qquad [\xi_i(u),h]=0\qquad\qquad [h,h']=0\]\\[-4ex]
\item[($\Y$2)] For any $i\in\bfI$, and $h\in\h$,
\[[h,x^\pm_i(u)]=\pm\alpha_i(h)x^\pm_i(u)\]\\[-4ex]
\item[($\Y$3)] For any $i,j\in \bfI$, and $a = \hbar d_ia_{ij}/2$
\[(u-v\mp a)\xi_i(u)x_j^{\pm}(v)=
(u-v\pm a)x_j^{\pm}(v)\xi_i(u)\mp 2a x_j^{\pm}(u\mp a)\xi_i(u)\]
\\[-4ex]
\item[($\Y$4)] For any $i,j\in \bfI$, and $a = \hbar d_ia_{ij}/2$
\begin{multline*}
(u-v\mp a) x_i^{\pm}(u)x_j^{\pm}(v)\\
= (u-v\pm a)x_j^{\pm}(v)x_i^{\pm}(u)
+\hbar\lp [x_{i,0}^{\pm},x_j^{\pm}(v)] - [x_i^{\pm}(u),x_{j,0}^{\pm}]\rp
\end{multline*}\\[-4ex]
\item[($\Y$5)] For any $i,j\in \bfI$
\[(u-v)[x_i^+(u),x_j^-(v)]=-\delta_{ij}\hbar\left(\xi_i(u)-\xi_i(v)\right)\]
\item[($\Y$6)] For any $i\neq j\in\bfI$, and $m=1-a_{ij}$, 

\[\sum_{\pi\in\Sym_m}
\left[x^{\pm}_i(u_{\pi(1)}),\left[x^{\pm}_i(u_{\pi(2)}),\left[\cdots,
\left[x^{\pm}_i(u_{\pi(m)}),x^{\pm}_j(v)\right]\cdots\right]\right]\right.=0\]
\end{enumerate}
\end{prop}

\rem When $i=j$, ($\Y$4) may be rewritten as follows. Taking
$u=v$ yields
\begin{equation}\label{eq:yang i=j}
[x_{i,0}^{\pm},x_i^{\pm}(u)] = \mp d_i x_i^{\pm}(u)^2
\end{equation}
and therefore the relation ($\Y4^\prime$)
\begin{multline*}
(u-v\mp \hbar d_i)x_i^{\pm}(u)x_i^{\pm}(v)\\
= (u-v\pm \hbar d_i)
x_i^{\pm}(v)x_i^{\pm}(u) \mp \hbar d_i \lp
x_i^{\pm}(u)^2 + x_i^{\pm}(v)^2\rp
\end{multline*}
Conversely, equating coefficients of $v^0$ in this identity yields
\eqref{eq:yang i=j}, and
therefore ($\Y$4).

\subsection{Proof of Proposition \ref{pr:Y fields}}
%----------------------------------------------------------

The equivalences (Y1)$\equiv$($\Y$1) and (Y6)$\equiv$($\Y$6)
are clear.

To see that (Y4) implies ($\Y$4), multiply (Y4) by $\hbar u^{-r-1}
\hbar v^{-s-1}$. Summing over $r,s\geq 0$, and using $\hbar\sum
_{r\geq 0}x_{i,r+1}^\pm u^{-r-1}=u x_i^\pm(u)-\hbar x^\pm_{i,0}$
yields
\begin{multline*}
(u-v)[x_i^\pm(u),x_j^\pm(v)]-
\hbar\left([x_{i,0}^\pm,x_j^\pm(v)]-[x_i^\pm(u),x_{j,0}^\pm]\right)\\
=\pm a\left(x_i^\pm(u)x_j^\pm(v)+x_j^\pm(v)x_i^\pm(u)\right)
\end{multline*}
as claimed. Conversely, taking coefficients of $u^{-r-1}v^{-s-1}$
in ($\Y$4) for $r,s\geq 0$ yields (Y4).\footnote{equating coefficients
of $u^0v^{-s-1}$, and of $u^{-r-1}v^{0}$ yields the tautological identities
$[x_{i,0}^\pm,x_{j,s}^\pm]=[x_{i,0}^\pm,x_{j,s}^\pm]$ and $[x_{i,r}
^\pm,x_{j,0}^\pm]=[x_{i,r}^\pm,x_{j,0}^\pm]$ respectively.}

(Y2) is clearly equivalent to ($\Y$2). To see that
(Y2)--(Y3) imply ($\Y$3), set $\xi_{i,-1}=\hbar^{-1}$ so that
\[\xi_i(u)=\hbar\sum_{r\geq -1}\xi_{i,r}u^{-r-1}\]
and (Y2) for $h=\xi_{i,0}$ coincides with (Y3) for $r=-1$. Multiplying
both sides of (Y3) by $\hbar u^{-r-1}\hbar v^{-s-1}$, and summing
over $r\geq -1,s\geq 0$ yields
\[(u-v\mp a)\xi_i(u)x^{\pm}_j(v) - (u-v\pm a)x_j^{\pm}(v)\xi_i(u)=
-\hbar[\xi_i(u),x_{j,0}^\pm]\]
Since the \rhs is independent of $v$, so is the left--hand side.
The latter is therefore equal to its value when $v=u\mp a$, that
is $\mp 2a x_j^{\pm}(u\mp a)\xi_i(u)$ and ($\Y$3) follows. Conversely,
taking coefficients of $u^{-r-1}v^{-s-1}$, $r,s\geq 0$ in ($\Y$3)
yields (Y3)\footnote{equating coefficients of $u^1$, $u^0$ yields
the identities $x^\pm_j(v)=x^\pm_j(v)$ and (Y2) for $h=\xi_{i,0}$
respectively, while equating coefficients of $v^0$ yields the identity 
$-\hbar[\xi_i(u),x_{j,0}^\pm]=\mp 2a x_j^{\pm}(u\mp a)\xi_i(u)$.}

To see that (Y5) implies ($\Y$5), multiply (Y5) by $\hbar u^{-r-1}
\hbar v^{-s-1}$. Summing over $r,s\geq 0$, and using
\[(u-v)\sum_{r+s=t}u^{-r-1}v^{-s-1}=v^{-t-1}-u^{-t-1}\]
then yields ($\Y$5). Conversely, the above identity implies that
$\xi_i(u)-\xi_i(v)=\hbar(v-u)\sum_{r,s\geq 0}u^{-r-1}v^{-s-1}\xi_
{i,r+s}$. Equating this to $\hbar^{-1}(v-u)[x^+_i(u),x^-_i(v)]$
yields (Y5) since $\C[u,v;u^{-1},v^{-1}]]$ has no zero divisors.

\subsection{}\label{ssec: Y6-reduction}
%--------------

The following result is due to Levendorskii \cite{levendorskii}.
We reproduce his proof below for completeness.

\begin{lem}\label{lem: Y6-reduction}
The relation (Y6) follows from (Y1)--(Y3) and the special case of (Y6)
when $r_1=\cdots=r_m=0$.
\end{lem}

\begin{pf}
The proof of this lemma is based on the construction of commuting
elements $\{\wt{t}_{i,r}\}_{i\in\bfI, r\in\N}$. These are polynomials
in $\{\xi_{i,r}\}$ such that the following commutation relations
hold
\begin{align*}
[\wt{t}_{i,r}, x^{\pm}_{i,s}] &= \pm 2d_i x^{\pm}_{i,r+s}\\
[\wt{t}_{i,r}, x^{\pm}_{j,s}] &= \sum_{k=0}^{r+s} b^{\pm}_{ij}(r,s,k)  
x^{\pm}_{j,r+s-k}
\end{align*}
where $b^{\pm}_{ij}(r,s,k)\in\C$ are scalars. We postpone their
construction to \ref{ssec: Y6-red-const}, and proceed with the
proof of the lemma.

Denote the \lhs of (Y6) by $S^{\pm}_{ij}(r_1,\cdots, r_m; s)$.
We prove that this element is zero by induction on the number
of non--zero $r_j$. The case $S_{ij}^{\pm}(0,\cdots, 0; s)=0$
holds by assumption. Assume now that
\[S_{ij}^{\pm}(r_1,\cdots,r_k,0,\cdots,0;s)=0\]
for some $0\leq k < m$ and for every $r_1,
\cdots, r_k,s\in\N$. Applying $\ad(\wt{t}_{i,r})$ to this identity yields
\begin{align*}
0 &= \pm 2d_i(m-k) S_{ij}^{\pm}(r_1,\cdots,r_k,r,0,\cdots,0;s)\\
&\phantom{=} \pm 2d_i \sum_{p=1}^k S_{ij}^{\pm}(r_1,\cdots,r_{p-1},
r_p+r,r_{p+1},\cdots, r_k, 0,\cdots,0; s)\\
&\phantom{=}+ \sum_{p=0}^{r+s} b_{ij}^{\pm}(r,s,p) 
S_{ij}^{\pm}(r_1,\cdots, r_k,0,\cdots, 0; r+s-p)
\end{align*}
Since the terms in the last two sums are zero by induction, it
follows that $S_{ij}^{\pm}(r_1,\cdots, r_k,r,0,\cdots,0;s)=0$.
\end{pf}

\subsection{}\label{ssec: Y6-red-const}
%--------------

Define $\{t_{i,r}\}_{i\in\bfI, r\in\N}$ by
\[
\hbar\sum_{r\geq 0} t_{i,r}u^{-r-1} = \log\lp \xi_i(u)\rp
\]
The following commutation relation was obtained in \cite[Lemma 1.9]
{levendorskii} (see also \cite[Remark 2.9]{sachin-valerio-1})
\[
[t_{i,r},x^{\pm}_{j,s}] = \pm d_ia_{ij}\sum_{l=0}^{\lfloor r/2\rfloor}
\left(\begin{array}{c} r\\ 2l\end{array}\right)
\frac{(\hbar d_ia_{ij}/2)^{2l}}{2l+1} x_{j,r+s-2l}^{\pm}\]
Now it is clear that one can define
$\wt{t}_{i,r}$ as 
\[
\wt{t}_{i,r} = t_{i,r} + \text{ linear combination of } t_{i,r-k}
\ (1\leq k\leq r)
\]
such that $[\wt{t}_{i,r}, x_{i,s}^{\pm}] = \pm 2d_i x_{i,r+s}^{\pm}$.
The commutation relation above, and the lower triangularity of 
the transformation matrix between $\{\wt{t}_{i,r}\}$ and 
$\{t_{i,r}\}$ 
implies that we also have
\[
[\wt{t}_{i,r},x^{\pm}_{j,s}] = \sum_{p=0}^{r+s} b^{\pm}_{ij}(r,s,p)
x^{\pm}_{j,r+s-p}
\]
for some scalars $b_{ij}^{\pm}(r,s,p)$.

\subsection{}\label{ssec: Y6-red-cor}
%------------------------------------

It follows from Lemma \ref{lem: Y6-reduction} that, on an integrable
representation, relation (Y6) can be deduced from the other relations.
More precisely, let $\wYhg$ be a unital associative algebra generated
by  $\{h,\xi_{i,r},x_{i,r}^{\pm}\}$ subject to the relations (Y0)--(Y5) of 
Section \ref{ssec: yangian}. Let $V$ be a representation of $\wYhg$
which is integrable in the following sense
\begin{itemize}
\item $V$ is $\h$--diagonalizable with \fd weight spaces.
\item For each $i\in\bfI$, the action of $x_{i,0}^{\pm}$ on $V$ is locally nilpotent.
\end{itemize}

\begin{prop}\label{prop: Y6-reduction}
The action of $\wYhg$ on $V$ factors through one of $\Yhg$.
\end{prop}
\begin{pf}
By Lemma \ref{lem: Y6-reduction}, it suffices to prove that relation
(Y6) with $r_1=\cdots=r_m=0$ holds in $\End(V)$. Let $\sl_2^i\subset
\Yhg$ be the Lie subalgebra generated by $\{\xi_{i,0},x_{i,0}^\pm\}$.
Let $V[\mu]$, $\mu\in\h^*$ be a weight subspace of $V$, and set
\[\ol{V}=\bigoplus_{m\in\Z}V[\mu+m\alpha_i]\aand
\ol{V}^\pm_j=\bigoplus_{m\in\Z}V[\mu+m\alpha_i\pm\alpha_j]\]
$\ol{V},\ol{V}_j^\pm$ are $\sl_2^i$--invariant subspaces of $V$ which,
by the integrability of $V$, are finite--dimensional. Since $[\xi_{i,0},x_{j,s}
^\pm]=\pm d_ia_{ij}x_{j,s}^\pm$ and $\ad(x_{i,0}^\mp)x_{j,s}^\pm=0$,
the restriction of $x_{j,s}^\pm$ to $\ol{V}$ is a lowest (resp. highest)
weight vector in the \fd $\sl_2^i$--module $\Hom(\ol{V},\ol{V}_j^\pm)$,
which implies that $\ad(x_{i,0}^\pm)^{1-a_{ij}}x_{j,s}^\pm$ acts by 0
on $\ol{V}$.
\end{pf}

\subsection{Shift automorphism}\label{ssec: shift-yangian}
%----------------------------------------------------

The group of translations of the complex plane acts on 
$\Yhg$ by
\[\tau_a(y_r) = \sum_{s=0}^r
\left(\begin{array}{c}r\\s\end{array}\right)
a^{r-s}y_s
\aand
\tau_a(h) = h\]
where $a\in\C$, $y$ is one of $\xi_i,x_i^\pm$, $r\in\N$ and $h\in\h$. In terms of
the generating series introduced in \ref{ssec: formal-series-y}, 
\[\tau_a(y(u)) = y(u-a)\]
Given a representation $V$ of $\Yhg$ and $a\in \C$, we set
$V(a)=\tau_a^*(V)$.

\subsection{The quantum loop algebra $\qloop$}\label{ssec: qla}
%=====================================

Let $q\in\C^\times$ be of infinite order. For any $i\in\bfI$, set $q
_i=q^{d_i}$. We use the standard notation for Gaussian integers
\begin{gather*}
[n]_q = \frac{q^n - q^{-n}}{q-q^{-1}}\\[.5ex]
[n]_q! = [n]_q[n-1]_q\cdots [1]_q\qquad
\qbin{n}{k}{q} = \frac{[n]_q!}{[k]_q![n-k]_q!}
\end{gather*}

The quantum loop algebra $\qloop$ is the $\C$--algebra generated
by elements $\{\Psi_{i,\pm r}^\pm\}_{i\in\bfI,r\in\N}$, $\{\X_{i,k}^\pm\}
_{i\in\bfI,k\in\Z}$ and $\{K_h\}_{h\in\h}$, subject to the following relations
\begin{itemize}
\item[(QL0)] For any $i\in\bfI$
\[\Psi_{i,0}^\pm=K_{\pm d_i\alpha^\vee_i}\]
\item[(QL1)] For any $i,j\in\bfI$, $r,s\in\N$ and $h,h'\in\h$,
\begin{gather*}
[\Psi_{i,\pm r}^\pm,\Psi_{j,\pm s}^\pm]=0
\qquad
[\Psi_{i,\pm r}^\pm,\Psi_{j,\mp s}^\mp]=0
\qquad
[\Psi_{i,\pm r}^\pm,K_h]=0\\[.7ex]
K_hK_{h'} = K_{h+h'}
\qquad
K_0=1
\end{gather*}
\item[(QL2)] For any $i\in\bfI$, $k\in\Z$ and $h\in\h$,
\[K_h\X^\pm_{i,k}K_h^{-1}=q^{\pm\alpha_i(h)}\X^\pm_{i,k}\]
\item[(QL3)] For any $i,j\in\bfI$, $\veps\in\{\pm\}$ and $l\in\Z$
\[\Psi^\veps_{i,k+1}\X^\pm_{j,l} - q_i^{\pm a_{ij}}\X^\pm_{j,l}\Psi^\veps_{i,k+1}
=
q_i^{\pm a_{ij}}\Psi^\veps_{i,k}\X^\pm_{j,l+1}-\X^\pm_{j,l+1}\Psi^\veps_{i,k}\]
for any $k\in\Z_{\geq 0}$ if $\veps=+$ and $k\in\Z_{<0}$ if $\veps=-$
\item[(QL4)] For any $i,j\in\bfI$ and $k,l\in \Z$
\[\X^\pm_{i,k+1}\X^\pm_{j,l} - q_i^{\pm a_{ij}}\X^\pm_{j,l}\X^\pm_{i,k+1}=
q_i^{\pm a_{ij}}\X^\pm_{i,k}\X^\pm_{j,l+1}-\X^\pm_{j,l+1}\X^\pm_{i,k}\]
\item[(QL5)] For any $i,j\in\bfI$ and $k,l\in \Z$
\[[\X^+_{i,k},\X^-_{j,l}] = \delta_{ij} \frac{\Psi^+_{i,k+l} - \Psi^-_{i,k+l}}{q_i - q_i^{-1}}\]
where we set $\Psi^{\pm}_{i,\mp k}=0$ for any $k\geq 1$.
\item[(QL6)] For any $i\neq j\in\bfI$, $m=1-a_{ij}$, $k_1,\ldots, k_m\in\Z$
and $l\in \Z$
\[\sum_{\pi\in \Sym_m} \sum_{s=0}^m (-1)^s\qbin{m}{s}{q_i}
\X^\pm_{i,k_{\pi(1)}}\cdots \X^\pm_{i,k_{\pi(s)}} \X^\pm_{j,l}\X^\pm_{i,k_{\pi(s+1)}}\cdots \X^\pm_{i,k_{\pi(m)}} = 0\]
\end{itemize}

\subsection{}\label{ssec: qla-fields}
%--------------
Define $\Psi_i(z)^+,\X^\pm_i(z)^+\in\qloop[[z^{-1}]]$ and 
$\Psi_i(z)^-,\X^\pm_i(z)^-\in \qloop[[z]]$ by
\begin{align*}
\Psi_i(z)^+ 	&= \sum_{r\geq 0}\Psi_{i,r}^+ z^{-r} & \Psi_i(z)^- 	&= \sum_{r\leq 0}\Psi_{i,r}^-z^{-r}\\
\X_i^\pm(z)^+	&= \sum_{r\geq 0}\X_{i,r}^\pm z^{-r}& \X_i^\pm(z)^-&=-\sum_{r< 0}\X_{i,r}^\pm z^{-r}
\end{align*}

\begin{prop}\label{pr:qloop fields}
The relations (QL1),(QL2)--(QL3),(QL4),(QL5),(QL6) imply the following
relations in $\qloop[z,w;z^{-1},w^{-1}]]$
\begin{enumerate}
\item[(\QL1)] For any $i,j\in\bfI$, and $h,h'\in\h$,
\begin{gather*}
[\Psi_i(z)^+,\Psi_j(w)^+]=0
\qquad
[\Psi_i(z)^+,K_h]=0\\[.7ex]
K_hK_{h'} = K_{h+h'}
\qquad
K_0=1
\end{gather*}
\item[(\QL2)] For any $i\in\bfI$, and $h\in\h$,
\[K_h\X_i^\pm(z)^+K_h^{-1}=q^{\pm\alpha_i(h)}\X_i^\pm(z)^+\]
\item[(\QL3)] For any $i,j\in\bfI$
\begin{multline*}
(z-q_i^{\pm a_{ij}}w)\Psi_i(z)^+\X_j^\pm(w)^+\\
=(q_i^{\pm a_{ij}}z-w)\X_j^\pm(w)^+\Psi_i(z)^+-
(q_i^{\pm a_{ij}} - q_i^{\mp a_{ij}})q_i^{\pm a_{ij}}w\X_j^\pm(q_i^{\mp a_{ij}}z)^+\Psi_i(z)^+
\end{multline*}
\item[(\QL4)] For any $i,j\in\bfI$
\begin{multline*}
(z-q_i^{\pm a_{ij}}w)\X_i^\pm(z)^+\X_j^\pm(w)^+-
(q_i^{\pm a_{ij}}z-w)\X_j^\pm(w)^+\X_i^\pm(z)^+\\
=z\lp\X_{i,0}^\pm\X_j^\pm(w)^+-q_i^{\pm a_{ij}}\X_j^\pm(w)^+\X_{i,0}^\pm\rp
+ w\lp\X_{j,0}^\pm\X_i^\pm(z)^+-q_i^{\pm a_{ij}}\X_i^\pm(z)^+\X_{j,0}^\pm\rp
\end{multline*}
\item[(\QL5)] For any $i,j\in\bfI$
\[(z-w)[\X^+_i(z)^+,\X^-_j(w)^+] =
\frac{\delta_{ij}}{q_i-q_i^{-1}}\left(z\Psi_i(w)^+-w\Psi_i(z)^+-(z-w)\Psi_{i,0}^-\right)\]
\item[(\QL6)] For any $i\neq j\in\bfI$, and $m=1-a_{ij}$
\begin{multline*}
\sum_{\pi\in \Sym_m} \sum_{s=0}^m (-1)^s\qbin{m}{s}{q_i}
\X^\pm_i(z_{\pi(1)})^+\cdots \X^\pm_i(z_{\pi(s)})^+ \X^\pm_j(w)^+\\
\cdot \X^\pm_i(z_{\pi(s+1)})^+\cdots \X^\pm_i(z_{\pi(m)})^+ = 0
\end{multline*}
\end{enumerate}
\end{prop}

\rem When $i=j$, (\QL4) may be rewritten as follows. Taking
$z=w$ yields
\begin{equation}\label{eq:loop i=j}
\lp\X_{i,0}^\pm\X_i^\pm(z)-q_i^{\pm 2}\X_i^\pm(z)\X_{i,0}^\pm\rp
=
(1-q_i^{\pm 2})\X_i^\pm(z)^2
\end{equation}
and therefore the relation (\QL4$'$)
\begin{multline*}
(z-q_i^{\pm 2}w)\X_i^\pm(z)\X_i^\pm(w)-
(q_i^{\pm 2}z-w)\X_i^\pm(w)\X_i^\pm(z)\\
=z(1-q_i^{\pm 2})\X_i^\pm(w)^2+
w(1-q_i^{\pm 2})\X_i^\pm(z)^2
\end{multline*}
Conversely, equating coefficients of $z$ in this identity yields
\eqref{eq:loop i=j} and therefore (\QL4).

\subsection{Proof of Proposition \ref{pr:qloop fields}}
%----------------------------------------------------------------

The implications (QL$n$)$\Rightarrow$(\QL$n$), for $n=1,2,6$ are clear.

Multiplying (QL4) by $z^{-k}w^{-l}$, summing over $k,l\geq 0$ and using
the fact that $\sum_{k\geq 0}\X^\pm_{i,k+1}z^{-k}=z(\X^\pm_i(z)-\X^\pm
_{i,0})$ yields(\QL4).

Multiplying (QL3) by $z^{-k}w^{-l}$, summing over $k,l\geq 0$ yields similarly
\begin{multline*}
(z-q_i^{\pm a_{ij}}w)\Psi_i(z)^+\X_j^\pm(w)^+-(q_i^{\pm a_{ij}}z-w)\X_j^\pm(w)^+\Psi_i(z)^+\\
=
z(\Psi^+_{i,0}\X_j^\pm(w)^+-q_i^{\pm a_{ij}}\X_j^\pm(w)^+\Psi^+_{i,0})
-w(q_i^{\pm a_{ij}}\Psi_i(z)^+\X_{j,0}^\pm-\X_{j,0}^\pm\Psi_i(z)^+)
\end{multline*}
Since the first summand on the \rhs is equal to 0 by (QL2), the ratio of the
\lhs by $w$ is independent of $w$. It is therefore equal to its value for $w=
q_i^{\mp a_{ij}}z$, that is $-q_i^{\pm a_{ij}} \left(
q_i^{\pm a_{ij}} - q_i^{\mp a_{ij}}
\right)
\X_j^{\pm}(q_i^{\mp a_{ij}}z)^+
\Psi_i(z)^+$  
and (\QL3) holds.

Multiplying (QL5) by $(z-w)z^{-k}w^{-l}$, summing over $k,l\geq 0$, and using
$(z-w)\sum_{k+l=m}z^{-k}w^{-l}=zw^{-m}-wz^{-m}$ yields (\QL5).

\begin{rem}\label{rm:missing-}
The relations (\QL 1)--(\QL 6) of Proposition \ref{pr:qloop fields}
for $\{\Psi_i(z)^-, \X_i^{\pm}(z)^-\}$ continue to hold in 
$\qloop[[z,w]]$. Similarly these relations remain true between 
$+$ and $-$ fields, in the appropriate algebra. For example, 
(\QL 3) for $\Psi_i(z)^+$ and $\X_i^{\pm}(w)^-$ is true
in $\qloop[z,w^{-1}; z^{-1},w]]$.
\end{rem}

\subsection{}\label{ssec: QL6-reduction}
%--------------

The following result is the analogue of Lemma \ref{lem: Y6-reduction}
for $\qloop$.

\begin{lem}\label{lem: QL6-reduction}
The relation (QL6) follows from (QL1)--(QL3) and the special case of
(QL6) when $k_1=\cdots=k_m=0$.
\end{lem}
\begin{pf} We proceed as in Lemma \ref{lem: Y6-reduction}. Define
elements $\{H_{i,k}\}_{i\in\bfI, k\in\nZ}\in\qloop$ by
\[
\Psi_i(z)^{\pm} = \Psi_{i,0}^{\pm} 
\exp\lp\pm (q_i-q_i^{-1})\sum_{n\geq 1} H_{i,\pm n}z^{\mp n}
\rp
\]
Then we have the following relation
\[
[H_{i,k}, \X_{j,l}^{\pm}] = \pm \frac{[ka_{ij}]_{q_i}}{k} \X_{j,k+l}^{\pm}
\]
The proof then proceeds by an induction argument, as in that of 
Lemma \ref{lem: Y6-reduction}.
\end{pf}

\subsection{}\label{ssec: QL6-red-cor}
%------------------------------------

Let $\wqloop$ be the unital associative algebra generated by elements
$\{K_h,\Psi_{i,\pm r}^{\pm}, \X^{\pm}_{i,l}\}$ subject to the relations (QL0)--(QL5)
of Section \ref{ssec: qla}. Let $\V$ be a representation of $\wqloop$
which is integrable in the following sense
\begin{itemize}
\item $\V$ is $\h$--diagonalizable with \fd weight spaces.
\item For each $i\in\bfI$, the action of $\X_{i,0}^{\pm}$ on $\V$
is locally nilpotent.
\end{itemize}
The proof of the following result is similar to that of Proposition \ref
{prop: Y6-reduction} and therefore omitted.

\begin{prop}\label{prop: QL6-reduction}
The action of $\wqloop$ on $\V$ descends to that of $\qloop$.
\end{prop}

\subsection{Shift automorphism}\label{ssec: shift-qla}
%--------------------------------------

The group $\C^*$ of dilations of the complex plane acts on $\qloop$
by
\[\tau_{\alpha}(Y_k)=\alpha^k Y_k
\aand
\tau_{\alpha}(K_h)=K_h\]
where $\alpha\in\C^*$, $Y$ is one of $\Psi^\pm_i,\X^\pm_i$, $k\in\Z$ and $h
\in\h$. In terms of the generating series of \ref{ssec: qla-fields}, we
have
\[\tau_{\alpha}\lp Y(z)^\pm\rp = Y(\alpha^{-1}z)^\pm\]

Given a representation $\V$ of $\qloop$ and $\alpha\in\nC$, we
denote $\tau_{\alpha}^*(\V)$ by $\V(\alpha)$.

\section{Integrable representations of $\Yhg$ and $\qloop$}
\label{sec: s-cat} 
%==============================================

\subsection{Integrable representations \cite{hernandez-affinizations}}\label{ssec: int-rep}
%-------------------------------------------------------------------------------------

By definition, a representation of $\qloop$ (resp. $\Yhg$) is in $\Rloop$
(resp. $\Ryang$) if, and only if its restriction to $U_q(\g)$ (resp. $\g$)
is an integrable representation in the corresponding category $\calO$
(see \cite{kac, lusztig-book}). Specifically,
\begin{enumerate}
\item A representation $\V$ of $\qloop$ is in $\Rloop$ if 
\begin{enumerate}
\item $\V = \bigoplus_{\mu\in\h^*} \V_{\mu}$ and $\dim(\V_{\mu})<\infty$, where
\[\V_{\mu} = \{v\in \V|\,K_hv = q^{\mu(h)}v \text{ for every } h\in\h\}\]
\item There exist $\lambda_1,\ldots, \lambda_r\in \h^*$ such that $\V_{\mu}\neq 0$
implies that $\mu \leq \lambda_i$ for some $i=1,\ldots, r$. Recall that we say 
$\mu\leq \lambda$ if $\lambda-\mu\in \sum_{i\in\bfI} \N\alpha_i$.
\item For each $\mu\in\h^*$ such that $\V_{\mu}\neq 0$ and $i\in\bfI$, there exists
$N>0$ such that $\V_{\mu-n\alpha_i}=0$ for every $n\geq N$.
\end{enumerate}

\item A representation $V$ of $\Yhg$ is in $\Ryang$ if
\begin{enumerate}
\item $V = \bigoplus_{\mu\in\h^*} V_{\mu}$ and $\dim(V_{\mu})<\infty$, where
\[
V_{\mu} = \{v\in V|\,hv = \mu(h)v \text{ for every } h\in\h\}
\]
\item There exist $\lambda_1,\ldots, \lambda_r\in \h^*$ such that $V_{\mu}\neq 0$
implies that $\mu \leq \lambda_i$ for some $i=1,\ldots, r$.
\item For each $\mu\in\h^*$ such that $V_{\mu}\neq 0$ and $i\in\bfI$, there exists
$N>0$ such that $V_{\mu-n\alpha_i}=0$ for every $n\geq N$.
\end{enumerate}
\end{enumerate}

\begin{rem}
When $\g$ is a semisimple Lie algebra, the categories $\Ryang$ and $\Rloop$
are the categories of \fd representations of $\Yhg$ and $\qloop$. Note that in
this case, the definition of an integrable representations in category $\calO$
for $\qloop$ given above differs from the usual one for $U_q\hat{\g}$ given
in \cite{kac, lusztig-book}.
\end{rem}

\subsection{Drinfeld polynomials: Yangians}\label{ssec: dp}
%-----------------------------------------------------

We recall below the classification of irreducible representations in $\Ryang$.
For $\g=\sl_2$, this result is proved in \cite{chari-pressley-yangian}. For the
general case, a proof can be obtained by adapting the arguments of \cite
{chari-pressley-qaffine-rep}, see \cite{hernandez-affinizations} and \cite
[Chaper 12]{chari-pressley}.

Let $\bfc = \{c_{i,r}\}_{r\in\N, i\in\bfI}\subset\C$ be a collection of complex
numbers, and let $\lambda\in\h^*$ such that $c_{i,0}=d_i\lambda(\alpha_i
^{\vee})$.  A representation $V$ of $\Yhg$ is said to be a \hw
representation of \hw $(\lambda,\bfc)$ if there exists $\fv\in V$ such that
\begin{enumerate}
\item $V = \Yhg\fv$.
\item $x_{i,r}^+\fv=0$ for every $i\in \bfI$, $r\in\N$.
\item $\xi_{i,r}\fv = c_{i,r}\fv$ and $h\,\fv=\lambda(h)\fv$, for every $i\in\bfI,r\in\N$
and $h\in \h$.
\end{enumerate}
To the pair $(\lambda,\bfc)$, we associate the Verma module 
$M(\lambda,\bfc)$ in an obvious way and denote by $L(\lambda,
\bfc)$ its unique irreducible quotient.

\begin{thm}\label{thm: dp}\hfill
\begin{enumerate}
\item Every irreducible representation in $\Ryang$ is a \hw representation
for a unique highest weight $(\lambda,\bfc)$.
\item The irreducible representation $L(\lambda,\bfc)$ is in $\Ryang$
if, and only if there exist (unique) monic polynomials $\{P_i(u)\in\C[u]\}
_{i\in\bfI}$ such that
\[1+\hbar\sum_{r\geq 0}c_{i,r}u^{-r-1} = \frac{P_i(u+d_i\hbar)}{P_i(u)}\]
\end{enumerate}
\end{thm}
It follows from (ii) that $\lambda(\alpha_i^{\vee})=\deg(P_i)\in\N$ and hence
that $\lambda$ is a dominant integral weight. The polynomials $\{P_i(u)\}$
are called Drinfeld polynomials.
By Theorem \ref{thm: dp}, the set of isomorphism classes of simple objects
in $\Ryang$ is in bijection with the set $\dwtY$ of pairs $(\lambda\in\h^*,\{P
_i(u)\in\C[u]\}_{i\in\bfI})$ such that $\lambda(\alpha_i^{\vee}) = \deg(P_i)$
and $P_i$ are monic. If $(\lambda,\{P_i\})\in\dwtY$, we denote by $L(\lambda,
\{P_i\})$ the corresponding \irr $\Yhg$--module.

\subsection{Drinfeld polynomials: quantum loop algebra}\label{ssec: dp-qla}
%-------------------------------------------------------------------------

Let $\gamma=\{\gamma_{i,\pm m}^{\pm}\}_{i\in\bfI, m\in \N}$ be a collection
of complex numbers and $\lambda\in\h^*$ such that 
$\gamma^{\pm}_{i,0} = q_i^{\pm\lambda(\alpha_i^{\vee})}$. 
A  representation $\V$ of $\qloop$ is said to be an $l$--\hw representation
of $l$--\hw $(\lambda,\gamma)$ if there exists
$\fv\in \V$ such that
\begin{enumerate}
\item $\V = \qloop \fv$.
\item $\X^+_{i,k}\fv = 0$ for every $i\in\bfI$ and $k\in\Z$.
\item $\Psi_{i,\pm m}^{\pm}\fv = \gamma_{i,\pm m}^{\pm}\fv$ and $K_h\fv=
q^{\lambda(h)}\fv$ for every $i\in \bfI,m\in\N$ and $h\in\h$.
\end{enumerate}

As for Yangians, for any $(\lambda,\gamma)$, there is a unique irreducible
representation, $\mathcal{L}(\lambda,\gamma)$, with highest weight $(\lambda,\gamma)$.

\begin{thm}\label{thm: qla-dp}\hfill
\begin{enumerate}
\item Every irreducible representation in $\Rloop$ is a \hw representation
for a unique highest weight $(\lambda,\gamma)$.
\item The \irr representation $\LL(\lambda,\gamma)$ is in $\Rloop$ if, and
only if there exist monic polynomials $\{\PP_i(w)\in \C[w]\}_{i\in\bfI}$, $\PP
_i(0)\neq 0$, such that
\[\sum_{m\geq 0} \gamma_{i,m}^{+} z^{-m}=
q_i^{-\deg(\PP_i)} \frac{\PP_i(q^{2}z)}{\PP_i(z)}=
\sum_{m\leq 0} \gamma_{i,m}^{-} z^m\]
\end{enumerate}
\end{thm}

Again, we have $\lambda(\alpha_i^{\vee}) = \deg(\PP_i)$ and hence $\lambda$
is a dominant integral weight. The polynomials $\PP_i$ are again called Drinfeld
polynomials. The set of isomorphism classes of simple objects in $\Rloop$ is in
bijection with the set $\dwtU$ of pairs $(\lambda\in\h^*,\{\PP_i\})$ such that $\PP_i$
are monic, $\PP_i(0)\neq 0$ and $\lambda(\alpha_i^{\vee})=\deg(\PP_i)$. We
denote by $\LL(\lambda,\{\PP_i\})$ the \irr $\qloop$--module corresponding to
$(\lambda,\{\PP_i\})\in\dwtU$.

\subsection{Composition series}\label{ssec: kac-lemma}
%----------------------------------------------------

We shall need the analogue for $\Yhg$ and $\qloop$ of the existence of
composition series in category $\calO$ for an arbitrary symmetrisable
Kac--Moody algebra $\g$, (see \cite[Lemma 9.6]{kac} and \cite[Prop. 15]
{hernandez-affinizations}).

\begin{lem}\label{lem: kac-lemma}\hfill
\begin{enumerate}
\item Let $\V\in\Rloop$ and $\lambda\in\h^*$. Then, there exists a filtration
of $\qloop$--modules
\[0=\V_0\subset\cdots\subset\V_t=\V\]
such that the following holds for any $\V_j$, $j=1,\ldots,t$
\begin{itemize}
\item either $\V_j/\V_{j-1}=\LL(\lambda_j,\{\PP^{(j)}_i\})$, for some
$\lambda_j\geq \lambda$,
\item or $\lp \V_j/\V_{j-1}\rp_{\mu} = 0$ for every $\mu\geq\lambda$.
\end{itemize}
Moreover, given $\mu\geq \lambda$ and $\lp\mu,\{\PP_i\}\rp\in\dwtU$, the number
of times $\LL\lp\mu,\{\PP_i\}\rp$ appears as $\V_j/\V_{j-1}$ is independent of $\lambda$
and the filtration chosen, and is denoted by  $[\V:\LL\lp\mu,\{\PP_i\}\rp]$.
\item Let $V\in\Ryang$ and $\lambda\in\h^*$. Then, there exists a filtration of
$\Yhg$--modules
\[0=V_0\subset\cdots\subset V_t = V\]
such that the following holds for any $V_j$, $j=1,\ldots,t$
\begin{itemize}
\item either $V_j/V_{j-1}=L(\lambda_j,\{P^{(j)}_i\})$ for some $\lambda_j\geq\lambda$.
\item or $\lp V_j/V_{j-1}\rp_{\mu} = 0$ for every $\mu\geq \lambda$.
\end{itemize}
Moreover, given $\mu\geq \lambda$ and $\lp\mu,\{P_i\}\rp\in\dwtY$, the number of 
times $L\lp\mu,\{P_i\}\rp$ appears as $V_j/V_{j-1}$ is independent of $\lambda$ and
the filtration chosen, and is denoted by  $[V:L\lp\mu,\{P_i\}\rp]$.
\end{enumerate}
\end{lem}

\subsection{Definition of the categories $\Rysub$ and $\Rlsub$}\label{ssec: s-cat}
%--------------------------------------------------------------------------------

Let $\Pi\subset\C$ be a subset such that $\Pi\pm\frac{\hbar}{2}\subset\Pi$.
We define $\Rysub$ to be the full subcategory of $\Ryang$ consisting of the
representations $V$ such that, for every $(\lambda,\{P_i\})\in\dwtY$
for which $[V:L(\lambda,\{P_i\})]\neq 0$, the roots of $P_i$ lie in $\Pi$ for every
$i\in\bfI$.\\

Let similarly $\Omega\subset\nC$ be a subset stable under multiplication by
$q^{\pm 1}$. We define $\Rlsub$ to be the full subcategory of $\Rloop$ consisting
of those $\V$ such that for every $(\lambda,\{\PP_i\})\in\dwtU$ for which $[\V:\LL
(\lambda,\{\PP_i\})]\neq 0$, the roots of $\PP_i$ lie in $\Omega$ for every
$i\in\bfI$.

\begin{prop}\label{prop: obvious-prop}\hfill
\begin{enumerate}
\item $\Rysub$ and $\Rlsub$ are Serre subcategories of $\Ryang$ and
$\Rloop$ respectively. That is, they are closed under taking direct sums,
subobjects, quotient and extensions.
\item When $\g$ is a simple Lie algebra, $\Rysub$ and $\Rlsub$ are closed
under tensor product.
\end{enumerate}
\end{prop}

\begin{pf}
The first part is obvious from the definition. The second part will follow from
the alternate characterizations of these categories given in Theorem \ref
{thm: first-main-theorem} (iv), since by \cite[Lemma 1]{knight} and \cite
{frenkel-reshetikhin-qchar} $\xi_i(u)$ and $\Psi_i^{\pm}(z)$ are group--like
modulo off--diagonal entries.
\end{pf}

\subsection{Rationality of fields}\label{ssec: rationality}
%--------------------------------------

The aim of this paragraph is to prove that the fields giving the action of the
generators of $\Yhg,\qloop$ on category $\calO$ integrable representations
are rational functions. The statement for $\qloop$ is well known, see \cite
[\S 6]{beck-kac} and \cite[Proposition 3.8]{hernandez-drinfeld-coproduct}.
Our proof is similar, but yields in addition the explicit form of these rational
functions.

\begin{prop}\label{prop: rationality}\hfill
\begin{enumerate}
\item Let $V$ be a $\Yhg$--module on which $\h$ acts semisimply
with \fd weight spaces. Then, for every weight $\mu$ of $V$, the
generating series 
\[\xi_i(u)\in\End(V_{\mu})[[u^{-1}]]
\qquad
x_i^{\pm}(u)\in\Hom(V_{\mu},V_{\mu\pm\alpha_i})[[u^{-1}]]\]
defined in \ref{ssec: formal-series-y} are the expansions at $\infty$
of rational functions of $u$. Specifically, let $t_{i,1} = \xi_{i,1} - \ds
\frac{\hbar}{2} \xi_{i,0}^2\in\Yhg^\h$. Then,
\[x_i^{\pm}(u)= \hbar u^{-1}\lp 1 \mp \frac{\ad(t_{i,1})}{2d_iu}\rp^{-1}x_{i,0}^{\pm}\]
and $\xi_i(u) = 1 + [x_i^+(u),x_{i,0}^-]$.
\item Let $\V$ be a $\qloop$--module on which the operators $K_h$ act
semisimply with \fd weight spaces. Then, for every weight $\mu$ of $\V$
and $\veps\in\{\pm\}$, the generating series
\[\Psi_i(z)^\pm\in\End(\V_\mu)[[z^{\mp 1}]]
\qquad
\X_i^\veps(z)^\pm\in\Hom(\V_\mu,\V_{\mu\pm\alpha_i}))[[z^{\mp 1}]]\]
defined in \ref{ssec: qla-fields} are the expansions of rational functions
$\Psi_i(z),\X_i^\veps(z)$ at $z=\infty$ and $z=0$. Specifically, let
$H_{i,\pm 1}=\pm\Psi_{i,0}^{\mp}\Psi_{i,\pm 1}^{\pm}/(q_i-q_i^{-1})$.
Then,
\[\X_i^\veps(z)
=\lp 1- \veps\frac{\ad(H_{i,1})}{[2]_iz}\rp^{-1}\X_{i,0}^\veps
=-z\lp 1- \veps z\frac{\ad(H_{i,-1})}{[2]_i}\rp^{-1}\X_{i,-1}^\veps\]
and
$\Psi_i(z) = \Psi_{i,0}^- + (q_i-q_i^{-1})[\X^+_i(z),\X^-_{i,0}]$.
\end{enumerate}
\end{prop}
\begin{pf} (i) The relations (Y2) and (Y3) imply that $[t_{i,1},
x_{i,r}^{\pm}] = \pm 2d_i x_{i,r+1}^{\pm}$, so that
\[[t_{i,1},x_i^\pm(u)]=\pm 2d_i(u x_i^\pm(u)-\hbar x_{i,0})\]
The relation $\xi_i(u) = 1+ [x_i^+(u),x_{i,0}^-]$ is a direct consequence of (Y5).

(ii) The relations (QL2) and (QL3) imply that 
$[H_{i,\pm 1}, \X^\veps_{i,k}] = \veps [2]_i \X^\veps_{i,k\pm 1}$, 
and therefore that
\begin{align*}
\X_i^\veps(z)^+ &= \phantom{-z}
\lp 1-\phantom{z}
\veps\frac{\ad(H_{i,1})}{[2]_iz}\phantom{z}\rp^{-1}
\X_{i,0}^\veps \\
\X_i^\veps(z)^- 	&= -z\lp 1- \veps z\frac{\ad(H_{i,-1})}{[2]_i}\rp^{-1}
\X_{i,-1}^\veps
\end{align*}
which are rational functions regular at $z=\infty$ and $z=0$ respectively.
To see that these functions are the same, it suffices to note that, as a
formal power series, and therefore as a rational function, $\X_i^\veps
(z)^-$ also satisfies
\[\lp1-\veps\frac{\ad(H_{i,1})}{[2]_iz}\rp\X_i^\veps(z)^-
= \X_{i,0}^\veps\]
Finally, note that $\Psi_i^\pm(z) = \Psi_{i,0}^- + (q_i-q_i^{-1}) [\X^+_i(z)^\pm,
\X^-_{i,0}]$.
\end{pf}

\subsection{}
%--------------

The following result will be needed elsewhere

\begin{lem}\label{le:control poles}\hfill
\begin{enumerate}
\item Let $V\in\Ryang$, and let $\mu\in\h^*$ be a weight of $V$.
Then, the poles of $\xi_i(u)_\mu$ are contained in the union of the
poles of $x_i^\pm(u)_\mu$ and those of $x_i^\pm(u)_{\mu\mp\alpha_i}$.
\item Let $\V\in\Rloop$, and let $\mu\in\h^*$ be a weight of $\V$.
Then the poles of $\Psi_i(z)_\mu$ are containes in the union of
the poles of $\X^\pm_i(z)_{\mu}$ and those of $\X^\pm_i(z)_{\mu\mp\alpha_i}$.\end{enumerate}
\end{lem}
\begin{pf}
(i) and (ii) follow from
\begin{gather*}
\xi_i(u)=1+[x^+_i(u),x^-_{i,0}]=1+[x^+_{i,0},x^-_i(u)]\\
\Psi_i(z)=
\Psi_{i,0}^-+(q_i-q_i^{-1}) [\X^+_i(z),\X^-_{i,0}]
=\Psi_{i,0}^-+(q_i-q_i^{-1}) [\X^+_{i,0},\X^-_i(z)]
\end{gather*}
respectively
\end{pf}

\subsection{Equivalent characterizations of $\Rysub$ and $\Rlsub$}\label{ss:characterise}
%-----------------------------------------------------------------------------------

In view of Proposition \ref{prop: rationality}, we can define, for a given
representation $(\rho,V)$ in $\Ryang$, a subset $\spec(V)\subset\C$
consisting of the poles of the rational functions $\{\rho(\xi_i(u))^{\pm 1},
\rho(x_i^{\pm}(u))\}$. Similarly, for a representation $(\rho,\V)$ in
$\Rloop$, we define a subset $\spec(\V)\subset\nC$ consisting of
poles of the functions $\{\rho(\Psi_i(z))^{\pm 1},\rho(\X^{\pm}_i(z))\}$.

\begin{thm}\label{thm: first-main-theorem}\hfill
\begin{itemize}
\item[(1)] Let $V\in\Ryang$. Then, the following conditions are equivalent.
\begin{enumerate}
\item $V\in \Rysub$.
\item $\spec(V)\subset \Pi$.
\item The poles of $\xi_i(u)^{\pm 1}$ are contained in $\Pi$.
\item The eigenvalues of $\xi_i(u)$ have zeroes and poles in $\Pi$.
\end{enumerate}
\item[(2)] Let $\V\in\Rloop$. Then, the following conditions are equivalent.
\begin{enumerate}
\item $\V \in \Rlsub$.
\item $\spec(\V)\subset \Omega$.
\item The poles of $\Psi_i(z)^{\pm 1}$ are contained in $\Omega$.
\item The eigenvalues of $\Psi_i(z)$ have zeroes and poles in $\Omega$.
\end{enumerate}
\end{itemize}
\end{thm}
\begin{pf}
The rest of this section is devoted to the proof of Theorem \ref{thm: first-main-theorem}
for the Yangian $\Yhg$.  The assertion for $\qloop$ is proved similarly and we
therefore omit it. The following implications are obvious: (ii) $\Rightarrow$ (iii)
$\Rightarrow$ (iv) $\Rightarrow$ (i). We prove (i) $\Rightarrow$ (ii) by induction
on the length of the composition series of $V$. More precisely, assume (i), and
fix a weight $\mu\in\h^*$. We need to prove that the operators $\{\xi_i(u)^{\pm 1}
_{\mu},x_i^{\pm}(u)_{\mu}\}_{i\in\bfI}$ have poles in $\Pi$. Fix an arbitrary $\lambda
\in\h^*$ such that $\lambda\leq \mu-\alpha_i$ for every $i\in\bfI$, and consider
a composition series given by Lemma \ref{lem: kac-lemma}. The proof of (i)
$\Rightarrow$ (ii) is by induction on the length of this composition series.
Proposition \ref{prop: irr-C} deals with the case of an irreducible $V$, while
Proposition \ref{prop: extn-C} carries out the induction step.
\end{pf}

\subsection{}\label{ssec: irr-C}
%--------------

\begin{prop}\label{prop: irr-C}
Let $V$ be an irreducible representation in $\Ryang$. Then $V\in \Rysub$
if, and only if $\spec(V)\subset \Pi$.
\end{prop}

\begin{pf}
It is clear that if $\spec(V)\subset \Pi$, the roots of $P_i(u)$ are in $\Pi$ (since roots
of $P_i(u)$ are obtained by successively shifting the poles of $\xi_i(u)$ acting on
the \hw vector by $d_i\hbar$, and $\Pi$ is closed under these shifts). Hence,
by definition, $V\in \Rysub$.

Assuming now that the roots of $P_i(u)$ are in $\Pi$, for every $i\in\bfI$,
we shall prove that all the rational functions $\{\xi_i(u),x_i^{\pm}(u)\}$ have
poles in $\Pi$. Let
\[
V = \bigoplus_{\nu\in\h^*} V_\nu
\]
be the weight space decomposition of $V$ as a $\g$--module and, for
$y\in \Yhg$, write $y_{\nu}$ for the restriction of $y$ to $V_\nu$.

Let $\mu\in\h^*$ be the \hw of $V$. Using the results of
\ref{ssec: dp} we have the following properties of $V$, which will be used
in the proof:

\begin{itemize}
\item[(P1)] Every weight space $V_\nu$, for $\nu < \mu$, is spanned by
$\{x_{i,r}^-V_{\nu+\alpha_i}\}_{i\in\bfI, r\in\N}$.
\item[(P2)] If $v\in V_\nu$ is annihilated by $x_{i,r}^+$ for every
$i\in\bfI$ and $r\in\N$, and $\nu\neq\mu$, then $v=0$.
\end{itemize}

The proof that each rational function $\xi_i(u)_{\nu},x_i^{\pm}(u)_{\nu}$
has poles in $\Pi$ is by induction on the height of $\mu-\nu$. Recall that
for an element $\alpha=\sum_i n_i\alpha_i\in Q_+$, the height of $\alpha$,
denoted by $\hit(\alpha)$, is defined as:
\[
\hit(\alpha) = \sum_i n_i \in \N
\]

More precisely, we prove the following statement by induction on $k\in\N$:\\

\noindent ${\mathbf S(k):}$ For every $i\in\bfI$, the rational functions $\xi_i(u)_{\nu}, 
x_i^+(u)_{\nu}$ have poles in $\Pi$, for every $\nu$ such that $\hit(\mu-\nu)
\leq k$; and $x_i^-(u)_{\nu}$ has poles in $\Pi$ for every $\nu$ such that
$\hit(\mu-\nu)<k$.\\

The base case ${\mathbf S(0)}$ is clear since $x_i^+(u)_{\mu}=0$ and
\[
\xi_i(u)_{\mu} = \frac{P_i(u+d_i\hbar)}{P_i(u)}
\]

Let us assume that ${\mathbf S(k')}$ holds for every $k'\leq k$, where $k\geq
0$, and prove ${\mathbf S(k+1)}$. Let $\nu$ be a weight of $V$ such that
$\hit(\mu-\nu)=k+1$ and let $i,j\in\bfI$. Using the relation ($\Y$5) we have:
\begin{equation}\label{eq: strip-cat-y5}
x_i^+(u)_{\nu}x_j^-(v)_{\nu+\alpha_j}
=x_j^-(v)_{\nu+\alpha_i+\alpha_j}x_i^+(u)_{\nu+\alpha_j}
+\frac{\delta_{ij}\hbar}{u-v}\lp \xi_i(v)_{\nu+\alpha_i}-\xi_i(u)_{\nu+\alpha_i}\rp
\end{equation}

Note that the \rhs has poles in $\Pi\times\Pi$ by the induction
hypothesis. This allows us to conclude the same for $x_i^+(u)_{\nu}$
and $x_j^-(v)_{\nu+\alpha_j}$ as follows. Assume that $x_j^-(v)_{\nu
+\alpha_j}$ has a pole at $z\in \C\setminus \Pi$ of order $n$. Multiplying
both sides of the equation \eqref{eq: strip-cat-y5} by $(v-z)^n$ and
setting $v=z$, we get:
\[
x_i^+(u)_{\nu} \lp \left. (v-z)^nx_j^-(v)_{\nu+\alpha_j}\right|_{v=z}\rp
=0
\]

Thus the image of the operator $\ds \lp \left. (v-z)^nx_j^-(v)_{\nu+\alpha_j}
\right|_{v=z}\rp$ is annihilated by all $x_{i,r}^+$, which implies that
this operator is zero, since $\nu\neq \mu$ (see property (P2) above).
This is a contradiction to the fact that $z$ is a pole of order $n$
of $x_j^-(v)_{\nu+\alpha_j}$. The proof for $x_i^+(u)_{\nu}$ is similar.

It remains to show that $\xi_i(u)_{\nu}$ has poles in $\Pi$. For this we
use the relation ($\Y$3)
\[
\xi_i(u)_{\nu}x_j^-(v)_{\nu+\alpha_j} = \frac{u-v-a}{u-v+a}
x_j^-(v)_{\nu+\alpha_j}\xi_i(u)_{\nu+\alpha_j}
+\frac{2a}{u-v+a} x_j^-(u+a)_{\nu+\alpha_j}\xi_i(u)_{\nu+\alpha_j}
\]
where $a = \hbar d_ia_{ij}/2$.

This relation implies that the poles of $\xi_i(u)$ acting on the image
of the operator $x_j^-(v)_{\nu+\alpha_j}$ in $V_\nu$ are either
the poles of $\xi_i(u)_{\nu+\alpha_j}$ or $x_j^-(u+a)_{\nu+\alpha_j}$,
which are in $\Pi$, by the induction hypothesis and the fact that $\Pi$
is stable under shift by $a$. Since $V_\nu$ is spanned by such
subspaces (property (P1)), we are done.
\end{pf}

\subsection{}\label{ssec: extn-C}
%--------------

\begin{prop}\label{prop: extn-C}
Let $0\to V_1\to V\to V_2\to 0$ be a short exact sequence of integrable 
$\Yhg$--modules. If
$\spec(V_1),\spec(V_2)\subset \Pi$ then $\spec(V)\subset\Pi$.
\end{prop}

\begin{pf}
Let us fix $\mu\in\h^*$ and $i\in\bfI$. We will prove that for every $k\in\Z$ the 
poles of $\xi_i(u)_{\mu+k\alpha_i}, x_i^{\pm}(u)_{\mu+k\alpha_i}$ 
acting on the corresponding weight space of $V$ are in $\Pi$.

We write $V = V_1\oplus V_2$ as vector spaces.
An element $y\in Y_{\hbar}(\g)$ has the following form (viewed as an element of
$\End(V)$):
\[
y = \left[ \begin{array}{cc} y^{11} & y^{12} \\ 0 & y^{22} \end{array}\right]
\]

By assumption, $\xi_i(u)^{ll}, x^{\pm}_i(u)^{ll}$ have poles in $\Pi$, for $l=1,2$. 
Let us assume that $z\not\in\Pi$ is a pole of one of the functions $\xi_i(u)_{\nu}^{12}, 
x_i^{\pm}(u)_{\nu}^{12}$ for $\nu\in\mu+\Z\alpha_i$. 
Let $N$ be the maximum of the order of the pole at $z$ of these
functions (note that by definition, there 
are only finitely many $k\in\Z$ such that 
$V_{\mu+k\alpha_i}\neq 0$). Define for every $\nu\in\mu+\Z\alpha_i$
\begin{gather*}
\mathsf{H}_{\nu} = \lim_{u\to z} (u-z)^N \xi_i(u)_{\nu}^{12}\\
\mathsf{X}^{\pm}_{\nu} = \lim_{u\to z} (u-z)^N x_i^{\pm}(u)_{\nu}^{12}
\end{gather*}

Using the relation ($\Y$5) and taking
its $(1,2)$ entry, we have:
\begin{multline*}
(u-v)(x_i^+(u)_{\nu-\alpha_i}^{11}x_i^-(v)_{\nu}^{12} + 
x_i^+(u)_{\nu-\alpha_i}^{12}x_i^-(v)_{\nu}^{22}\\
- x_i^-(v)_{\nu+\alpha_i}^{11}x_i^+(u)_{\nu}^{12} - 
x_i^-(v)_{\nu+\alpha_i}^{12}x_i^+(u)_{\nu}^{22})
 = \hbar(\xi_i(v)_\nu^{12} - \xi_i(u)_\nu^{12})
\end{multline*}
Multiplying both sides by $(u-z)^N$ and letting $u\to z$ we get
\[
(z-v)(\mathsf{X}^+_{\nu-\alpha_i} x_i^-(v)_{\nu}^{22} - 
x_i^-(v)_{\nu+\alpha_i}^{11} \mathsf{X}^+_{\nu}) = 
-\hbar\mathsf{H}_{\nu}
\]

and similarly (using $v-z$ instead of $u-z$) we have
\[
(u-z)(x_i^+(u)_{\nu-\alpha_i}^{11}\mathsf{X}_{\nu}^- 
- \mathsf{X}^-_{\nu+\alpha_i}
x_i^+(u)_{\nu}^{22}) = \hbar\mathsf{H}_{\nu}
\]
These equations clearly imply that $\mathsf{H}_{\nu} = 0$ and 
we have the following relations:
\begin{gather}
\mathsf{X}^+_{\nu-\alpha_i} x_i^-(v)_{\nu}^{22} 
= x_i^-(v)_{\nu+\alpha_i}^{11} \mathsf{X}_{\nu}^+ \label{eq: pf-extn-1} \\
x_i^+(u)_{\nu-\alpha_i}^{11}\mathsf{X}_{\nu}^- = \mathsf{X}_{\nu+\alpha_i}^-x_i^+(u)_{\nu}^{22} \label{eq: pf-extn-2}
\end{gather}

Now we consider the relation ($\Y$3) in the following form
\begin{multline*}
(u-v-d_i\hbar)\xi_i(u)x_i^+(v) - (u-v+d_i\hbar)x_i^+(v)\xi_i(u) \\
= -d_i\hbar(\xi_i(u)x_i^+(u) + x_i^+(u)\xi_i(u))
\end{multline*}
Again taking its $(1,2)$ component, multiplying with $(u-z)^N$ and letting
$u\to z$, and using the fact that $\mathsf{H}=0$ we get:
\begin{equation}\label{eq: pf-extn-3}
\xi_i(z)_{\nu+\alpha_i}^{11}\mathsf{X}_{\nu}^+ + \mathsf{X}_{\nu}^+\xi_i(z)_{\nu}^{22}= 0
\end{equation}

Similar computation with relation ($\Y$4) yields

\begin{multline}\label{eq: pf-extn-4}
(z-v-d_i\hbar)\mathsf{X}_{\nu+\alpha_i}^+x_i^+(v)_{\nu}^{22} - 
(z-v+d_i\hbar)x_i^+(v)_{\nu+\alpha_i}^{11}\mathsf{X}_{\nu}^+
\\ = -d_i\hbar(x_i^+(z)_{\nu+\alpha_i}^{11}\mathsf{X}_{\nu}^+ + 
 \mathsf{X}_{\nu+\alpha_i}^+x_i^+(z)_{\nu}^{22})
\end{multline}

We take commutator of \eqref{eq: pf-extn-4} with $x^-_{i,0}$, using the
commutativity property \eqref{eq: pf-extn-1}, the equation
\eqref{eq: pf-extn-3}  and the relation $[x_i^+(u),x_{i,0}^-]
 = \xi_i(u)-1$, to obtain:
\[
(z-v-d_i\hbar)\mathsf{X}_{\nu}^+\xi_i(v)_{\nu}^{22} - 
(z-v+d_i\hbar)\xi_i(v)_{\nu+\alpha_i}^{11}\mathsf{X}_{\nu}^+ = 0
\]
If we set $v=z-d_i\hbar$ in this equation, we get
\[
-2d_i\hbar \xi_i(z-d_i\hbar)_{\nu+\alpha_i}^{11}\mathsf{X}_{\nu}^+ = 0
\]
But since $z-d_i\hbar\not\in\Pi$ and $\xi_i(u)^{11}$ is invertible there, we obtain
that $\mathsf{X}_{\nu}^+=0$. A similar argument shows that 
$\mathsf{X}_{\nu}^-=0$,
which contradicts the fact that $N$ was the order of the pole at $z$ of
one of $\xi_i(u)_{\nu}^{12}, x_i^{\pm}(u)_{\nu}^{12}$. This contradiction proves that all the
poles are in $\Pi$ and we are done.
\end{pf}

\section{Additive difference equations}\label{sec: difference}
%=============================

In this section, we give a self--contained account of Birkhoff's 
theory of additive difference equations, including a detailed
discussion of the inverse problem for abelian ones. Our exposition,
and understanding of the subject owes much to \cite{birkhoff-difference,
borodin,krichever}. Additional information may be found in \cite
{birkhoff-riemann,van-der-Put-Singer}.

\subsection{}\label{ssec: diff-coeff}
%------------------------------------------

Let $V$ be a \fd vector space over $\C$, and let $A:\C\to\End(V)$
be a rational function. We assume that $A$ is regular at $\infty$
and such that $A(\infty)=1$. Consider the following system of
difference equations for a function $\phi:\C\to\End(V)$
\begin{equation}\label{eq: difference}
\phi(u+1) = A(u)\phi(u)
\end{equation}

\noindent Let $\A\subseteq\End(V)$ be the subalgebra generated
by $A(u)$, $u\in\C$, or equivalently by the coefficients of the Taylor
expansion
\[A(u)=1+A_0 u^{-1}+A_1 u^{-2}+\cdots \]
of $A$ near $\infty$. The system \eqref{eq: difference} is said to be
{\it non--resonant} \relto a subalgebra $\wt{\A}\subseteq\End(V)$
containing $\A$ if the eigenvalues of $\ad(A_0)$ on $\wt{\A}$ do not
lie in $\nZ$.
Note that this is the case in either of the following situations:
\begin{itemize}
\item The eigenvalues of $A_0$ do not differ by positive integers,
and $\wt{\A}=\End(V)$.
\item The system \eqref{eq: difference} is abelian, that is $[A(u),
A(v)]=0$ for any $u,v\in\C$, and $\wt{\A}=\A$.
\end{itemize}

\noindent
If \eqref{eq: difference} is non--resonant relative to $\wt{\A}$, it is
easy to see that it possesses a unique formal solution of the form
$\Upsilon(u)u^{A_0}$, where\footnote{For the formal solution
$\Upsilon(u)u^{A_0}$, the equation \eqref{eq: difference} is
understood to mean $\Upsilon(u+1)=A(u)\Upsilon(u)(1+1/u)^{-A_0}$,
where $(1+1/u)^{-A_0}=\sum_{r\geq 0}(-1)^r\frac{A_0(A_0+1)\cdots
(A_0+r-1)}{r!}u^{-r}$.}
\begin{equation}\label{eq:Upsilon}
\Upsilon(u)=
1+\Upsilon_0u^{-1}+\Upsilon_1u^{-2}+\cdots\in\wt{\A}[[u^{-1}]]
\end{equation}

\subsection{Canonical fundamental solutions}\label{ss:fundamental}
%-------------------------------------------------------

\begin{thm}\label{thm: difference0}\cite{birkhoff-difference}
If the system \eqref{eq: difference} is non--resonant \relto $\wt{\A}
\subseteq\End(V)$, there exist unique meromorphic solutions $\phi
^\pm:\C\to\wt{\A}$ such that
\begin{enumerate}
\item $\phi^\pm$ is holomorphic and invertible for $\rightleft$.
\item $\phi^\pm$ possesses an asymptotic expansion of the form
\begin{equation}\label{eq:unique asymptotic}
\phi^\pm(u)\thicksim
(1+H^\pm_0 u^{-1}+H^\pm_1 u^{-2}+\cdots)\cdot (\pm u)^{A_0}
\end{equation}
in any right (resp. left) half--plane, 
where $v^{A_0}=\exp(A_0\log(v))$ is defined by taking the
standard determination of the logarithm on $\C\setminus\R_{\leq 0}$.
\end{enumerate}
Moreover,
\begin{enumerate}
\item[(iii)] The poles of $(\phi^-)^{\pm 1}$ are contained in $\calP+\nN
$ and $\calZ+\nN$ respectively, where $\calP,\calZ$ are
the poles of $A(u),A(u)^{-1}$.
\item[(iv)] The poles of $(\phi^+)^{\pm 1}$ are contained in $\calZ-\N$
and $\calP-\N$ respectively.
\item[(v)] The asymptotics \eqref{eq:unique asymptotic} of $\phi^
\pm$ are given by $\Upsilon(u)(\pm u)^{A_0}$.
\end{enumerate}
\end{thm}
\begin{pf}
The uniqueness of $\phi^\pm$ is proved in \ref{ssec: fs-unique}.
In \ref{ssec: fs-exist I}--\ref{ssec: fs-exist II}, we give a proof of
the existence of the fundamental solutions $\phi^\pm$ under the
simplifying assumption that $A_0$ commutes with $A(u)$ for every
$u$. The general case is treated in \cite{birkhoff-difference, krichever}.

(iii)--(iv) If $\phi^-(u)$ is holomorphic and invertible for $\Re(u)<a$,
then, for $\Re(u)<a+n$, $\phi^-(u)=A(u-1)\cdots A(u-n)\phi^-(u-n)$
(resp. $(\phi^-(u))^{-1}=\phi^-(u-n)^{-1}A(u-n)^{-1}\cdots A(u-1)^{-1}$),
which has poles only if one of $u-1,\ldots,u-n$ is a pole of $A(u)$
(resp. of $A(u)^{-1}$). Similarly for $\phi^+$.

(v) Let $H^\pm(u)\in1+u^{-1}\wt{\A}[[u^{-1}]]$ be such that
$\phi^\pm(\pm u)^{-A_0}\thicksim H^\pm(u)$ for $\pm\Re(u)
>a$. Multiplying both sides of $\phi^\pm(u+1)=A(u)\phi^\pm
(u)$ by $(\pm(u+1))^{-A_0}$, using the fact that $(\pm u)^
{A_0}(\pm(u+1))^{-A_0}=(1+1/u)^{-A_0}$ for $u\notin\R_{\leq
0}$ (resp. $u\notin\R_{\geq -1}$) and taking asymptotic expansions,
shows that $H^\pm(u+1)=A(u)H^\pm(u)(1+1/u)^{-A_0}$, so
that $H^\pm(u)=\Upsilon(u)$ by uniqueness.
\end{pf}

\subsection{Uniqueness of fundamental solutions}
\label{ssec: fs-unique}
%------------------------------------------------------------

Let $\phi^+_1(u),\phi^+_2(u)$ be two solutions of \eqref{eq: difference}
satisfying the conditions (i)--(ii) of Theorem \ref{thm: difference0}, and
set $C(u)=\lp\phi_1^+(u)\rp^{-1}\phi^+_2(u)$. $C(u)$ is $1$--periodic,
and holomorphic for $\Re(u)\gg0$ and therefore on the entire complex
plane. As a holomorphic function of $z=\Exp{u}=\Exp{\Re(u)}e^{-2\pi
\Im(u)}$, $C$ has removable singularities at $z=0$ and $z=\infty$.
Indeed,
\begin{equation}\label{eq:zeta}
C(u)=
u^{-A_0}\lp\phi_1^+(u)u^{-A_0}\rp^{-1}\lp\phi_2^+(u)u^{-A_0}\rp u^{A_0}
\end{equation}
Since $C$ is 1--periodic, we may assume that $\Re u\gg0$, so that the
second and third factors tend to 1 as $\Im u\to\pm\infty$. Since the first
and last ones grow like a polynomial in $u$, it follows that $\lim_{z\to 0}
zC(z)=0=\lim_{z\to\infty}z^{-1}C(z)$. Thus, $C(u)\equiv C$ is a constant
element of $\wt{\A}$. By \eqref{eq:zeta},
\[u^{A_0}Cu^{-A_0}\thicksim 1+C^0 u^{-1}+\cdots\]
for $\Re u\gg0$. Let $A_0=S_0+N_0$ be the Jordan decomposition of
$A_0$, and $C=\sum_{\lambda\in\C}C_\lambda$ the decomposition
of $C$ into eigenvectors of $\ad(S_0)$. Since
\[\Ad(u^{A_0})C_\lambda=
u^\lambda\sum_{k\geq 0}\log(u)^k\frac{\ad(N_0)^k}{k!}C_\lambda\]
it follows that each $C_\lambda$ is an eigenvector of $\ad(A_0)$, that
$C_\lambda=0$ for any $\lambda\notin\Z_{\leq 0}$, and that $C_0=1$.
The non--resonance condition then implies that $C_n=0$ for any $n\in
\Z_{<0}$, so that $C=1$. The uniqueness of $\phi^-$ follows in the
same way.

\subsection{Existence of fundamental solutions, I}
\label{ssec: fs-exist I}
%------------------------------------------------------------

Assume first that $A_0=0$. The system \eqref{eq: difference} is then
called {\it regular}, and its fundamental solutions are given by
\begin{align*}
\phi^+(u)	&= \prod_{n\geq 0}^{\rightarrow} A(u+n)^{-1}=A(u)^{-1}A(u+1)^{-1}\cdots\\
\phi^-(u) 	&= \prod_{n\geq 1}^{\rightarrow} A(u-n)=A(u-1)A(u-2)\cdots
\end{align*}
It is easy to prove that the above products converge locally uniformly
on the complement of $\calZ-\N$, $\calP+\nN$ respectively,
and have the required asymptotics.

\subsection{Holomorphic functional calculus \cite{dunford-schwartz}}
\label{ssec: banach-calculus}
%------------------------------------------------------------------------------------

Recall that given $X\in\End(V)$, and a meromorphic function $f$
which is holomorphic on a neighborhood of the eigenvalues of $X$,
$f(X)\in\End(V)$ is defined by the Cauchy integral
\[f(X)=\cint{f(v)(v-X)^{-1}}{v}\]
where $C$ is any contour enclosing all of the eigenvalues of $X$
and none of the poles of $f(v)$. If $X$ is semisimple, with eigenvalues
$\lambda_1,\ldots,\lambda_k$ and corresponding eigenspaces
$V_1,\ldots,V_k$, the residue theorem shows that $f(X)$ is the
semisimple endomorphism acting as multiplication by $f(\lambda_i)$
on $V_i$. More generally, if $X = X_S + X_N$ is the Jordan
decomposition of $X$, expanding $(v-X)^{-1}$ as
\[(v-X)^{-1} = \sum_{k\geq 0} (v-X_S)^{-k-1} X_N^k\]
a finite sum since $X_N$ is nilpotent, shows that
\begin{equation}\label{eq:f(X)}
f(X) = \sum_{k\geq 0} \frac{f^{(k)}(X_S)}{k!} X_N^k
\end{equation}

If $f(u)=\sum_{i=0}^n a_i u^i$ is a polynomial, then $f(X)=\sum_i a_iX^i$.
Moreover, $f_1\cdot f_2(X)=f_1(X)f_2(X)$ when both sides are defined.
Finally, if $f_n$ is sequence converging locally uniformly to $f$, then $f
_n(X)$ converges to $f(X)$. It follows for example that if
\[\Gamma(v) = v^{-1}e^{-\gamma v}\prod_{n\geq 1}\lp 1 +\frac{v}{n}\rp^{-1}
e^{v/n}\]
is the Euler Gamma function, where $\gamma$ is the Euler--Mascheroni
constant, then
\begin{equation}\label{eq:Gamma X}
\Gamma(X) = X^{-1} e^{-\gamma X}\prod_{n\geq 1}
\lp 1+ \frac{X}{n}\rp^{-1}e^{X/n}
\end{equation}

\subsection{Existence of fundamental solutions, II}
\label{ssec: fs-exist II}
%------------------------------------------------------------

Assume now that $A_0$ commutes with $A(u)$. We may then regularize
\eqref{eq: difference} as follows. Set
\[\ol{A}(u)=(1-A_0u^{-1})A(u)\]
The difference equation $\ol{\phi}(u+1) = \ol{A}(u)\ol{\phi}(u)$ is regular,
and its fundamental solutions $\ol{\phi}^{\pm}$ can be constructed as in
\ref{ssec: fs-exist I}. The fundamental solutions to \eqref{eq: difference}
are $\phi^{\pm}(u)=f^\pm(u)\ol{\phi}^{\pm}(u)$, where $f^{\pm}(u)$ are
the solutions of $f(u+1)=u(u-A_0)^{-1}f(u)$ given by
\begin{align*}
f^+(u) &= \Gamma(u)\Gamma(u-A_0)^{-1}\\
f^-(u) &= \Gamma(1-u+A_0)\Gamma(1-u)^{-1}
\end{align*}
We will show in \ref{ss:f pm} that $f^{\pm}(u)$ have the required asymptotics
as $u\to\infty$. Assuming this for now, and using \eqref{eq:Gamma X}, we
can rewrite $\phi^\pm$ as
\begin{align}
\phi^+(u) &= e^{-\gamma A_0} A(u)^{-1}\prod_{n\geq 1}^{\rightarrow} A(u+n)^{-1}e^{A_0/n}
\label{eq:phi +}\\
\phi^-(u) &= e^{-\gamma A_0}\prod_{n\geq 1}^{\rightarrow}
A(u-n)\,e^{A_0/n}
\label{eq:phi -}
\end{align}

\subsection{Asymptotics of $f^\pm$}\label{ss:f pm}
%-------------------------------------------

Recall Sterling's asymptotic expansion \cite[\S 12.33]{whittaker-watson}
\[\Gamma(u)\thicksim
e^{-u}u^{u-1/2}\sqrt{2\pi}\left(1+1/12u+1/288u^2+\cdots\right)\]
valid on $\C\setminus\R_{\leq 0}$. This implies that, for any $\lambda\in
\C$
\[\frac{\Gamma(u-\lambda)}{\Gamma(u)}\thicksim u^{-\lambda}(1+h_1u^{-1}+\cdots)\]
and, differentiating \wrt $\lambda$, that, for any $k\geq 0$
\[\frac{\Gamma^{(k)}(u)}{\Gamma(u)}
\thicksim\lp\log{(u)}\rp^k(1+h_1^{(k)}u^{-1}+\cdots)\]
Let now $A_0=S_0+N_0$ be the Jordan decomposition of $A_0$.
By \eqref{eq:f(X)},
\begin{align*}
f^+(u)^{-1}
&= \Gamma(u-A_0)\Gamma(u)^{-1}\\
&= \frac{\Gamma(u-S_0)}{\Gamma(u)} \sum_{k\geq 0} \frac{\Gamma^{(k)}(u-S_0)}
{\Gamma(u-S_0)} \frac{(-N_0)^k}{k!} \\
&\thicksim u^{-A_0}(1+\cdots )
\end{align*}
Similarly, $f^-(u) \sim (1+\cdots )(-u)^{A_0}$.

\subsection{Connection Matrix}\label{ssec: diff-conn}
%-----------------------------------------------------------------

Define the {\it connection matrix} of the system \eqref{eq: difference}
to be the function $S:\C\to GL(V)$ given by
\[S(u)=\lp\phi^+(u)\rp^{-1}\phi^-(u)\]
Clearly, $S(u)$ is $1$--periodic, and therefore a function of $z=
\Exp{u}$.

\begin{prop}
S is a rational function of $z$, regular at $0,\infty$, and such
that
\[S(\infty)=e^{\pi\iota A_0}\aand S(0)=e^{-\pi\iota A_0}\]
\end{prop}
\begin{pf}
By Theorem \ref{thm: difference0}, (iii)--(iv) $S$ is a meromorphic
function with finitely many poles in the $z$--plane. To determine
its behaviour as $z\to\infty$ (resp. $z\to 0$), it suffices to
consider the limit of $S(u)$ as $\Im(u)\to\mp\infty$, with $\Re
u$ lying in an interval of the form $[a,a+1[$. Since both
asymptotic expansions $\phi^\pm\thicksim\Upsilon(u)(\pm u)
^{A_0}$ hold on the strip $\Re(u)\in[a,a+1[$, we get
\[\begin{split}
S(u)
&=
u^{-A_0}\left(\phi^+(u)u^{-A_0}\right)^{-1}\left(\phi^-(u)(-u)^{-A_0}\right)(-u)^{A_0}\\
&\thicksim
u^{-A_0}\Upsilon(u)^{-1}\Upsilon(u)(-u)^{A_0}\\
&=
u^{-A_0}(-u)^{A_0}
\end{split}\]
which is equal to $\exp(\pm\pi\iota A_0)$ depending on whether
$\Im(u)\lessgtr 0$.
\end{pf}

Note that, if $[A_0,A(u)]=0$, it follows by \eqref{eq:phi +}--\eqref
{eq:phi -} that
\[\begin{split}
S(z)
&=
\prod_{n\geq 1}^{\leftarrow} A(u+n)e^{-A_0/n}\,A(u)\,
\prod_{n\geq 1}^{\rightarrow}A(u-n)e^{A_0/n}\\[1.1ex]
&=
\cdots A(u+2)A(u+1)
\,A(u)\,
A(u-1)A(u-2)\cdots
\end{split}\]
where the last expression is understood as $\ds{\lim_
{n\to\infty}A(u+n)\cdots A(u+1)A(u)}$ $\ds{A(u-1)\cdots
A(u-n)}$.

\subsection{Example}\label{ssec: diff-ex}
%-------------------------

The following basic example illustrates the constructions given
above, and plays an important role in the computations below.
Consider the  scalar difference equation
\[\phi(u+1) = \frac{u-a}{u-b}\phi(u)\]
with coefficient matrix $A(u) = \ds \frac{u-a}{u-b}$. The canonical
fundamental solutions $\phi^\pm$ are given by
\[\phi^+(u)=\frac{\Gamma(u-a)}{\Gamma(u-b)}
\aand
\phi^-(u)=\frac{\Gamma(1-u+b)}{\Gamma(1-u+a)}\]
The connection matrix $S(u)$ is therefore equal to
\[S(u) = \frac{\Exp{u} - \Exp{a}}{\Exp{u}-\Exp{b}}\cdot e^{\pi\iota(b-a)}\]
where we used \cite[\S 12.14]{whittaker-watson}
\begin{equation}\label{eq:sine}
\Gamma(u)\Gamma(1-u)=\pi/\sin(\pi u)
\end{equation}

\subsection{Inverse problem}\label{ssec: diff-inv}
%-----------------------------------

Fix now $A_0\in\End(V)$, and let $S$ be a rational function of
$z=\Exp{u}$ such that $S(\infty)=e^{\pi\iota A_0}$ and $S(0)=
e^{-\pi\iota A_0}$. We now discuss the problem of finding a
coefficient matrix $A(u)$ of the form $1+A_0 u^{-1}+\cdots$,
and such that $S(z)$ is the connection matrix of the system
\eqref{eq: difference} determined by $A(u)$.

This may be phrased as the following factorization problem:
find two meromorphic functions $\phi^\pm(u):\C\to\End(V)$ such
that
\begin{enumerate}
\item[(a)] $\phi^{\pm}(u)$ are holomorphic and invertible for $\rightleft$.
\item[(b)] $\phi^\pm$ possesses an asymptotic expansion of
the form
\[\phi^\pm(u)\thicksim
(1+H^\pm_0 u^{-1}+H^\pm_1 u^{-2}+\cdots)\cdot (\pm u)^{A_0}\]
valid in any right (resp. left) half--plane.
\item[(c)] $S(\Exp{u}) = \phi^+(u)^{-1}\phi^-(u)$.
\end{enumerate}
Once such a factorization is found, the coefficient matrix can be
reconstructed as
\[A(u)
=\phi^-(u+1)\phi^-(u)^{-1}
=\phi^+(u+1)\phi^+(u)^{-1}
\]
which is a rational function of $u$ since its poles are contained
in a strip $|\Re u|< a$ by (a), and it tends to 1 as $u\to\infty$ by
(b). By uniqueness, the canonical fundamental solutions of the
difference equation defined by $A(u)$ are $\phi^\pm$, so that its
connection matrix is $S$ by (c).

The above factorization, and therefore the inverse problem do not
always admit a solution. We refer the reader to \cite{krichever} for
several interesting cases when a solution can be obtained. In \ref
{ss:start factorisation}--\ref{ss:end factorisation}, we shall restrict
our attention to the abelian case, that is assume that $[S(z), S(w)]
=0$. In this case a solution always exists, and can be written explicitly
in terms of Gamma function.

\subsection{Uniqueness of the factorization} 
%-----------------------------------------------------

Assume that $\phi_1^\pm(u),\phi_2^\pm(u)$ satisfy the properties
(a)--(c) of \ref{ssec: diff-inv}, and set
\[G(u)=\phi^+_2(u)\phi^+_1(u)^{-1}=\phi^-_2(u)\phi^-_1(u)^{-1}\]
$G$ is meromorphic and, by (b) tends to 1 as $u\to\infty$, so that
$G$ is a rational function. Thus, $\phi_2^\pm(u)=G(u)\phi_1^\pm(u)$,
and the corresponding coefficient matrices are related by the {\it
isomonodromic transformation}
\[A_{\phi_2}(u)=
\phi^\pm_2(u+1)\phi^\pm_2(u)^{-1}=
G(u+1)A_{\phi_1}(u)G(u)^{-1}\]

The following gives a sufficient condition for the uniqueness
of $A$. Call two subsets $B_1,B_2\subset\C$ {\it non--congruent} if
$b_1-b_2\notin\nZ$ for any $b_1\in B_1$, $b_2\in B_2$. We shall also
say that $B\subset\C$ is non--congruent if it is non--congruent to itself.

\begin{prop}\label{pr:unique factorisation}
Let $A(u),A'(u):\C\to\ GL(V)$ be rational functions of the form $1+
A_0 u^{-1}+\cdots$, and assume that the connection matrices of
the difference equations determined by $A$ and $A'$ are equal.
Let $\calP,\calZ\subset\C$ (resp. $\calP',\calZ'$) be the poles of
$A(u),A(u)^{-1}$ (resp. $A'(u),A'(u)^{-1}$). If the pairs
\[(\calZ,\calP),\quad(\calZ',\calP'),\quad(\calZ,\calZ'),\quad(\calP,\calP')\]
are non--congruent, then $A=A'$.
\end{prop}
\begin{pf}
By the foregoing, there exists a rational function $G(u):\C\to GL(V)$
with $G(\infty)=1$, and such that $A'(u)=G(u+1)A(u)G(u)^{-1}$. Thus,
$G(u+1)=A'(u)G(u)A(u)^{-1}$ so that, for any $n\in\nN$,
\[G(u)=A'(u-1)\cdots A'(u-n)\,G(u-n)\,A(u-n)^{-1}\cdots A(u-1)^{-1}\]
Since $G(u)$ is regular for $\Re(u)\ll 0$, this implies that the poles
of $G$ are contained in $(\calP'\cup\calZ)+\nN$. Similarly, $G(u)=
A'(u)^{-1}G(u+1)A(u)$, so that for any $n\in\N$
\[G(u)=A'(u)^{-1}\cdots A'(u+n-1)^{-1}\,G(u+n)\,A(u+n-1)\cdots A(u)\]
implying that the poles of $G$ are contained in $(\calZ'\cup\calP)-\N$.
The non--congruence assumption then implies that $G$ has no
poles, hence $G\equiv 1$ and $A'=A$ as claimed.
\end{pf}

\subsection{}\label{ssec: diff-lemma}
%--------------

\begin{lem}\label{lem: diff-lemma}
Let $M:\C\to GL(V)$ be a rational function such that
\[[M(u),M(v)]=0\]
for any $u,v\in\C$. Then, the semisimple and unipotent components
$M_S,M_U:\C\to GL(V)$ of the Jordan decomposition of $M$ are
rational functions.
\end{lem}
\begin{pf}
Let $\ga\subset\End(V)$ be the span of $\{M(u)|\,u\in \C\}$. By assumption,
$\ga$ is an abelian Lie subalgebra of $\End(V)$, so that $V = \bigoplus_
{\lambda\in\ga^*} V_{\lambda}$, where $V_{\lambda}$ is the generalized
eigenspace
\[V_{\lambda} = \{v\in V|\,(x-\lambda(x))^Nv=0,\ \forall x\in \ga, N\gg0\}\]
Set $M'_S(u)=\sum_{\lambda} {\mathbf 1}_{\lambda} \circ \lambda(M(u))$,
where ${\mathbf 1}_{\lambda}:V\to V_{\lambda}$ is the projection operator,
and $M_U'(u)=M'_S(u)^{-1}M(u)$. By construction, $M_S'(u),M_U'(u)$
commute, are semisimple and unipotent respectively, and satisfy $M(u)=
M'_S(u)M'_U(u)$. By uniqueness of the Jordan decomposition $M_S=M_
S'$ and $M_U=M_U'$ and the conclusion follows since $M'_S$ and $M'_U$
are rational functions.
\end{pf}

\subsection{Poles of abelian matrix functions}\label{ss:poles}
%-------------------------------------------------------

Retain the notation of \S \ref{ssec: diff-lemma}, and let $\calP,\calP_
S,\calP_U$ be the sets of poles of $M(z), M_S(z),M_U(z)$, and $\calZ,
\calZ_S,\calZ_U$ those of $M(z)^{-1},M_S(z)^{-1},M_U(z)^{-1}$ respectively.

\begin{lem}
The following holds
\begin{enumerate}
\item $\calP=\calP_S\cup\calP_U$ and $\calZ=\calZ_S\cup\calZ_U$.
\item $\calP_U=\calZ_U$.
\item The eigenvalues of $M(z)$ are rational functions of $z$, and
$\calP_S$ (resp. $\calZ_S$) consist of the poles (resp. zeros) of
those eigenvalues.
\end{enumerate}
\end{lem}
\begin{pf}
(ii) follows from $M_U(z)^{\pm 1}=\exp(\pm m_U(z))$, with
$m_U=\log(M_U(z))$. (i) and (iii) follow by choosing a basis of
$V$ in which $M_S(z)$ are all diagonal and $M_U(z)$ strictly
upper triangular.
\end{pf}

\subsection{}\label{ss:start factorisation}
%--------------

Fix $A_0\in\End(V)$, and let $S:\C\to GL(V)$ be a rational function
such that
\begin{itemize}
\item $\ds{S(\infty)=\exp(\pi\iota A_0)=S(0)^{-1}}$\\[-2ex]
\item $\ds{[S(z),S(w)]=0}$\\[-2ex]
\item $\ds{[A_0,S(z)]=0}$
\end{itemize}

\noindent
To factorise $S(z)$, we shall need to choose logarithms of its poles
and zeros. Specifically, apply \S \ref{ss:poles} to $M(z)=S(z)$ and
set
\[\calZ_S=\{\alpha_i\} _{i\in I}\qquad\qquad
\calP_S=\{\beta_j\}_{j\in J}\qquad\qquad
\calP_U=\{\delta_k\}_{k\in K}=\calZ_U\]
Let $(\lambda(z),\mu)$ be a joint eigenvalue of $(S(z),A_0)$. Then,
\[\lambda(z)=
e^{\pi\iota\mu}\prod_{i\in I',j\in J'}(z-\alpha_i)^{n_i}/(z-\beta_j)^{m_j}\]
for some $I'\subset I,J'\subset J$, and $n_i,m_j\in\nN$.
Since $S(\infty)=S(0)^{-1}$, we get $\prod\beta_j^{m_j}\alpha_i^{-n_i}
=\Exp{\mu}$. Choose complex numbers $\{a^{(\lambda,\mu)}_i\}_
{i\in I'}$ and $\{b^{(\lambda,\mu)}_j\}_{j\in J'}$ such that 
$\Exp{a_i^{(\lambda,\mu)}}=\alpha_i$, $\Exp{b_j^{(\lambda,\mu)}}=\beta_j$,
and
\begin{equation}\label{eq:consistent choice}
\sum_j m_jb_j^{(\lambda,\mu)}-\sum_i n_ia^{(\lambda,\mu)}_i=\mu
\end{equation}
Let $\{\delta_k\}_{k\in K'}$ be the poles of the restriction of $S_U(z)$ 
to the generalised eigenspace of $S(z)$ corresponding to $(\lambda,\mu)$,
and choose $\{d_k^{(\lambda,\mu)}\}_{k\in\wt{K}}$ such that $\Exp{d_k
^{(\lambda,\mu)}}=\delta_k$. Finally, set
\begin{gather*}
\{a_i\}_{i\in\wt{I}}=\bigcup_{(\lambda,\mu)}\{a^{(\lambda,\mu)}_i\}_{i\in I'}\qquad\qquad
\{b_j\}_{j\in\wt{J}}=\bigcup_{(\lambda,\mu)}\{b^{(\lambda,\mu)}_j\}_{j\in J'}\\
\{d_k\}_{k\in\wt{K}}=\bigcup_{(\lambda,\mu)}\{d^{(\lambda,\mu)}_k\}_{k\in K'}
\end{gather*}

\begin{thm}\label{thm:unique coeff}
For any consistent choice as above, there is a rational function $A:\C\to
GL(V)$ such that
\begin{itemize}
\item $A=1+A_0 u^{-1}+\cdots$\\[-2ex]
\item $[A(u),A(v)]=0$\\[-2ex]
\item the connection matrix of the system defined by $A(u)$ is $S(z)$\\[-2ex]
\item the poles of $A(u)$ (resp. $A(u)^{-1}$) are $\{b_j\}\cup\{d_k\}$
(resp. $\{a_i\}\cup\{d_k\}$).
\end{itemize}
Moreover, if none of the sets $\{a_i\},\{b_j\},\{d_k\}$ are congruent to each
other, or to themselves, $A$ is unique.
\end{thm}
\begin{pf}
The uniqueness of $A(u)$ follows from Proposition \ref{pr:unique factorisation}.
To prove the existence of $A(u)$, let $A_0=A_0^S+A_0^N$ be the Jordan
decomposition of $A_0$. Since
\[S_S(\infty)=e^{\pi\iota A_0^S}=S_S(0)^{-1}
\aand
S_U(\infty)=e^{\pi\iota A_0^N}=S_U(0)^{-1}\]
it suffices to treat the cases when $S(z)$ is semisimple and unipotent separately,
which we do in \ref{ss:semisimple} and \ref{ss:end factorisation} respectively.
\end{pf}

\subsection{Semisimple case}\label{ss:semisimple}
%-----------------------------------

Assume that $A_0$ and $S(z)$ are semisimple, let $\ga\subset\End(V)$
be the abelian Lie algebra spanned by $A_0$ and $S(w)$, $w\in\IP^1$,
and decompose $V=\bigoplus_{\lambda\in\ga^*}V_\lambda$ as the direct
sum of weight spaces under $\ga$. Since it suffices to factorise $S(z)$ on
each $V_\lambda$, the problem reduces to the scalar case. We may
therefore assume that
\[A_0\in\C\aand S(z)=e^{\pi\iota A_0}\prod_i\frac{z-\alpha_i}{z-\beta_i}\]
where $\alpha_i,\beta_i$ are such
that $\prod_i\alpha_i\beta_i^{-1}=e^{-2\pi\iota A_0}$, so that $S(0)=e
^{-\pi\iota A_0}$. Let $a_i,b_i\in\C$ be such that
\[\Exp{a_i}=\alpha_i\qquad\quad
\Exp{b_i}=\beta_i\qquad\quad
\sum_i b_i-a_i=A_0\]
Then, by Example \ref{ssec: diff-ex}, the factorization problem
is solved by
\[\phi^+=\prod_i\frac{\Gamma(u-a_i)}{\Gamma(u-b_i)}
\aand
\phi^-=\prod_i\frac{\Gamma(1-u+b_i)}{\Gamma(1-u+a_i)}\]
The corresponding coefficient matrix is
\[A(u)=\phi^+(u+1)\phi^+(u)^{-1}=
\prod_i\frac{u-a_i}{u-b_i}\]

\subsection{Unipotent case}\label{ss:end factorisation}
%--------------------------------- 

If $A_0$ is nilpotent, and $S(z)$ unipotent, taking logarithms reduces
the factorization problem to finding meromorphic functions $\varphi^
\pm:\C\to\End(V)$ such that
\begin{enumerate}
\item[(a)] $\varphi^\pm$ is holomorphic for $\rightleft$
\item[(b)] $\varphi^\pm(u)$ has an asymptotic expansion of the
form
\[\varphi^\pm(u)\thicksim
A_0\log(\pm u)+h^\pm_0 u^{-1}+h^\pm_1 u^{-2}+\cdots\]
valid in any right (resp. left) half--plane.
\item[(c)] $\varphi^-(u) - \varphi^+(u) = \log(S(\Exp{u}))$
\end{enumerate}

\noindent
Since $S(\infty)=e^{\pi\iota A_0}=S(0)^{-1}$, it follows that
\[\log (S(z))=
\pi\iota A_0 + \sum_k \sum_{r\geq 1}\frac{C_{r,k}}{(z-\delta_k)^r}\]
where $\{\delta_k\}$ are the poles of $S(z)$, and $C_{r,k}\in\End(V)$
are such that
\begin{equation}\label{eq:relation}
\sum C_{r,k}(-\delta_k)^{-r} = -2\pi\iota A_0
\end{equation}
To obtain an additive factorization of $\log(S(z))$, we first determine
meromorphic functions $\varphi^\pm_r(u,d)$ for any $r\geq 1$ and
$d\in\C$, which are holomorphic for $\rightleft$, and such that
\begin{equation}\label{eq:recursion}
\varphi_r^-(u,d)-\varphi_r^+(u,d)=\frac{1}{(z-\delta)^r}
\end{equation}
where $z=\Exp{u}$ and $\delta=\Exp{d}$.

Set $\Psi^+(u)=\Gamma'(u)/\Gamma(u)$ and $\Psi^-(u)=\Gamma'
(1-u)/\Gamma(1-u)$. Taking logarithmic derivatives in $\Gamma
(1-u)\Gamma(u)=\pi/\sin(\pi u)$, we get
\[\Psi^-(u)-\Psi^+(u)=
\pi\cot(\pi u)=\pi\iota\left(\frac{2}{z-1}+1\right)\]
so that the functions
\[\varphi^\pm_1(u,d)=
\frac{1}{2\delta}\left(\frac{\Psi^\pm(u-d)}{\pi\iota}\pm\frac{1}{2}\right)\]
solve \eqref{eq:recursion} for $r=1$. Next, if $\varphi^\pm_r(u,d)$ solve
\eqref{eq:recursion} for a given $r$, then
\[\begin{split}
\partial_u\varphi^-_r(u,d)-\partial_u\varphi^+_r(u,d)
&=
-\frac{2\pi\iota r z}{(z-\delta)^{r+1}}\\
&=-2\pi\iota r\left(\frac{\delta}{(z-\delta)^{r+1}}+\frac{1}{(z-\delta)^r}\right)
\end{split}\]
so that $\ds{\varphi_{r+1}^\pm(u,d)=-\frac{1}{\delta}\left(\frac{\partial_u}
{2\pi\iota r}+1\right)\varphi^\pm_r(u,d)}$ solve \eqref{eq:recursion} for
$r+1$. We therefore recursively get
\begin{equation}\label{eq:varphi_r}
\varphi^\pm_r(u,d)=
-\frac{(-\delta)^{-r}}{2}\prod_{k=1}^{r-1}\left(\frac{\partial_u}{2\pi\iota k}+1\right)
\left(\frac{\Psi^\pm(u-d)}{\pi\iota}\pm\frac{1}{2}\right)
\end{equation}
Since $\Psi^\pm(u)=\log(\pm u)+O(u^{-1})$ in any right (resp. left) half--plane
\cite[\S 12.33]{whittaker-watson}, where $O(u^{-1})$ denotes an element of
$u^{-1}\C[[u^{-1}]]$, it follows that
\[\varphi_r^{\pm}(u,d)
\thicksim
-(-\delta)^{-r}\left(
\frac{\log(\pm u)}{2\pi\iota}\pm \frac{1}{4}\right)+O(u^{-1})\] 

Choose now $d_k\in\C$ such that $\Exp{d_k}=\delta_k$, and set
\[\varphi^\pm(u)=\mp\frac{\pi\iota A_0}{2}+\sum_{k,r}\varphi_r^\pm(u,d_k)C_{k,r}\]
then $\varphi^-(u)-\varphi^+(u)=\log S(z)$ and, in any right (resp. left) half--plane,
\[\begin{split}
\varphi^\pm
&\thicksim
\mp\frac{\pi\iota A_0}{2}-
\left(\frac{\log(\pm u)}{2\pi\iota}\pm\frac{1}{4}\right)\sum C_{k,r}(-\delta_k)^{-r}+O(u^{-1})\\
&=
\mp\frac{\pi\iota A_0}{2}+
2\pi\iota A_0\left(\frac{\log(\pm u)}{2\pi\iota}\pm\frac{1}{4}\right)+O(u^{-1})\\
&=
A_0\log(\pm u)+O(u^{-1})
\end{split}\]
where we used \eqref{eq:relation}, as required.

To determine the coefficient matrix of the difference equation, let $\Delta$
be the operator defined by $\Delta f(u)=f(u+1)-f(u)$. Since $\Delta\Psi^{\pm}(u)=
1/u$, \eqref{eq:varphi_r} yields
\[\Delta\varphi^\pm_r(u,d)=
-\frac{(-\delta)^{-r}}{2\pi\iota}\prod_{k=1}^{r-1}\left(\frac{\partial_u}{2\pi\iota k}+1\right)
\frac{1}{u-d}\]
so that $A(u)=\exp(\Delta\varphi^\pm)$ is given by
\[A(u)
=\exp\left(
-\frac{1}{2\pi\iota}\sum_{k,r}(-\delta_k)^{-r}\prod_{k=1}^{r-1}\left(\frac{\partial_u}{2\pi\iota k}+1\right)
\frac{C_{k,r}}{u-d_k}\right)\]

\section{The functor $\Fh{}$}\label{sec: functor}
%=====================

In this section, we construct a functor $\Fh{}$ from a dense
subcategory of $\Ryang$ to $\Rloop$.

\subsection{Non--congruent representations}\label{ssec: category-ch}
%------------------------------------------------------

We shall say that $V\in\Ryang$ is {\em non--congruent} if,
for any $i\in\bfI$, the poles of $x^+_i(u)$ (resp. $x^-_i(u)$)
are not congruent modulo $\Z$. Let $\Rync$ be the full
subcategory of $\Ryang$ consisting of non--congruent
representations.

\subsection{Contours}\label{ssec: contours}
%-------------------------

By a Jordan curve $\calC$, we shall mean a disjoint union
of simple, closed curves in $\C$, the inner domains of which
are pairwise disjoint. By definition, the inner (resp. outer)
domain of $\calC$ is the union (resp. intersection) of those
of its connected components. If $\calC$ is a Jordan curve,
and $f$ a continuous function on $\calC$, we shall set
\[\ointC f(u)du=\frac{1}{2\pi\iota}\intC f(u)du\]

The definition of $\Fh{}$ relies upon the following contours
of integration. Given $V\in\Ryang$, a weight $\mu\in\h^*$
of $V$, and $i\in\bfI$, we shall denote by $\calC^\pm_{i,\mu}
\subset\C$ a Jordan curve such that
\begin{enumerate}
\item $\calC^\pm_{i,\mu}$ encloses the poles of $x^\pm_i(u)_\mu$.
\item $\calC^\pm_{i,\mu}$ encloses no $\nZ$--translate of the
poles of $x^\pm_i(u)_\mu$ and $x^\pm_i(u)_{\mu\pm\alpha_i}$.
\end{enumerate}
Note that such a curve exists for any $i\in\bfI$ and weight $\mu
\in\h^*$ if $V$ is non--congruent.

\subsection{Quantum loop operators}\label{ssec: ql-operators}
%---------------------------------------------

Given $V\in\Rync$, and a weight $\mu$ of $V$, we define the
action of the generators of $\qloop$ on the $\mu$--weight space
$V_\mu$ of $\Fh{}(V)=V$ as follows.
\begin{enumerate}
\item For any $h\in \h$, set $(K_h)_{\mu}=q^{\mu(h)}$.
\\[.1ex]%

\item For any $i\in\bfI$, consider the abelian, additive difference equation
\begin{equation}\label{eq: diff-eq-functor}
\phi_i(u+1)_{\mu} = \xi_i(u)_{\mu} \phi_i(u)_{\mu}
\end{equation}
with coefficient matrix $\xi_i(u)_\mu=1+\hbar\xi_{i,0}u^{-1}+\cdots$.
Let $S_i(u)_{\mu}$ be its connection matrix, as defined in \ref
{ssec: diff-conn}. Set
\[\Psi_i(z)_{\mu}=\left.S_i(u)_{\mu}\right|_{e^{2\pi\iota u}=z}\]
and let $\Psi_i(z)^{\pm}_{\mu}\in\End(V_\mu)[[z^{\mp 1}]]$
be its Taylor expansions at $\infty$ and $0$ respectively.
\end{enumerate}

To define the action of the remaining generators of $\qloop$,
consider the fundamental solutions $\phi_i^\pm(u)_\mu$ of
\eqref{eq: diff-eq-functor}, as defined in \ref{ss:fundamental}.
Define $g_i^\pm(u)_\mu:\C\to GL(V_\mu)$ by
\begin{equation}\label{eq:gi+-}
g_i^-(u)_{\mu} = \phi_i^-(u)_{\mu}
\aand
g_i^+(u)_{\mu} = \phi_i^+(u+1)_{\mu}^{-1}
\end{equation}
Note that, by \eqref{eq:phi +}--\eqref{eq:phi -}, 
\begin{align}
g_i^-(u)_{\mu}&=
e^{-\gamma\hbar\xi_{i,0}}\prod_{n\geq 1}^{\rightarrow}
\xi_i(u-n)_\mu\,e^{\hbar\xi_{i,0}/n}
\label{eq:gi-}\\
g_i^+(u)_{\mu}&=
\lp\prod_{n\geq 1}^{\leftarrow}
e^{-\hbar\xi_{i,0}/n}\xi_i(u+n)_\mu\rp\,e^{\gamma\hbar\xi_{i,0}}
\label{eq:gi+}
\end{align}

\begin{enumerate}
\item[(iii)] Let $c_i^{\pm}\in\nC$ be scalars such that
\begin{equation}\label{eq: constants}
c_i^-c_i^+=
d_i\Gamma(\hbar d_i)^2
\end{equation}
For any $i\in\bfI$ and $k\in \Z$, the action of $\X_{i,k}
^\pm$ as an operator $V_\mu\to V_{\mu\pm\alpha_i}$ is given
by
\begin{equation}
(\X^{\pm}_{i,k})_\mu
= c_i^{\pm}
\oint_{\calC^\pm_{i,\mu}}
e^{\pp ku}g_i^{\pm}(u)_{\mu\pm\alpha_i}x_i^{\pm}(u)_{\mu}\, du 
\label{eq: fun-Ek}
\end{equation}
where the contours $\calC^\pm_{i,\mu}$ are chosen as in
\ref{ssec: contours}.
\end{enumerate}

\begin{rem}
The above definition implies that the generating functions
\[\X_i^{\pm}(z)^+_{\mu}=\sum_{k\geq 0}\X^{\pm}_{i,k}z^{-k}
\aand
\X_i^{\pm}(z)^-_{\mu}=-\sum_{k<0}\X^{\pm}_{i,k}z^{-k}\]
are the expansions at $z=\infty$ and $z=0$ of the rational
function 
\begin{equation}\label{eq:fun X(z)}
\X_i^{\pm}(z)_\mu 
= c_i^{\pm}
\oint_{\calC^\pm_{i,\mu}}
\lp \frac{z}{z-e^{\pp u}}\rp g_i^{\pm}(u)_{\mu\pm\alpha_i}
x_i^{\pm}(u)_{\mu}\,du
\end{equation}
where $z$ lies outside the contour $\exp(2\pi\iota\calC^\pm_{i,\mu})$. 
\end{rem}

\begin{rem}\label{rem: independence}
The operator \eqref{eq: fun-Ek} is independent of the choice of
the contour $\calC^\pm_{i,\mu}$ satisfying the conditions of \S
\ref{ssec: contours}. Indeed, by \eqref{eq:gi-}--\eqref{eq:gi+} and
Lemma \ref{le:control poles}, the poles of $g_i^\pm(u)_{\mu\pm
\alpha_i}$ are contained in the $\Z_{\lessgtr 0}$--translates of
the poles of $x^\pm_i(u)_\mu$ and $x^\pm_i(u)_{\mu\pm\alpha_i}$.
In particular, $g_i^\pm(u)$ is holomorphic inside $\calC^\pm_{i,\mu}$.

It also follows that, for any finite collection of weights $\Xi\subset
\h^*$ of $V$, a contour $\calC^\pm_i$ can be chosen which
satisfies (i)--(ii) of \S \ref{ssec: contours} for all $\mu\in\Xi$.
\end{rem}

\subsection{}\label{ss:1st main thm}
%--------------

\begin{thm}\label{thm: second-main-theorem}\hfill\break
\begin{enumerate}
\item The above operators give rise to an action of $\qloop$ on $V$.
They therefore define an exact, faithful functor
\[\Fh{}:\Rync\longrightarrow\Rloop\]
\item The functor $\Fh{}$ is compatible with shift automorphisms. That
is, for any $V\in\Rync$ and $a\in\C$, 
\[\Fh{}(V(a)) = \Fh{}(V)(\Exp{a})\]
\item Let $\Pi\subset\C$ be a non--congruent subset such that $\Pi\pm
\frac{1}{2}\hbar\subset\Pi$. Then, $\Rysub$ is a subcategory
of $\Rync$, and $\Fh{}$ restricts to a functor 
\[\Fh{\Pi} : \Rysub\longrightarrow\Rlsub\]
where $\Omega=\exp(2\pi\iota\Pi)$.
\end{enumerate}
\end{thm}

\begin{pf}
(i) The relations (QL0), (QL1) and (QL2) clearly hold. We prove (QL3)
in \ref{ssec: pf-QL23}, (QL4) in \ref{ssec: pf-QL4} and (QL5) in \ref
{ssec: pf-QL5}. Relation (QL6) is a consequence of Proposition 
\ref{prop: QL6-reduction}.\\

{\em A word on notation:} the relations of $\qloop$ need to be verified
on a given weight space. For notational convenience, we omit the subscript
indicating the weight space in the proofs below. Moreover, the contours
$\calC^\pm_i$ are chosen to satisfy the conditions of \S \ref{ssec: contours}
for all the weight spaces under consideration, consistently with Remark
\ref{rem: independence}.

(ii) is proved in \ref{ssec: functor-shift}.

(iii) follows from Theorem \ref{thm: first-main-theorem} and the monodromy
computation in \ref{ssec: diff-ex}.
\end{pf}

\subsection{}
%--------------
Set $\Gamma^\pm(u)=\Gamma(1\pm u)$, and note that the identity \eqref
{eq:sine} reads
\begin{equation}\label{eq:sine'}
\Gamma^-(u)\Gamma^+(u)=
\frac{2\pi\iota u}{e^{\pi\iota u}-e^{-\pi\iota u}}
\end{equation}
Let $V\in\Ryang$, $i,j\in\bfI$, and $a=\hbar d_ia_{ij}/2$. We shall need the following

\begin{prop}\label{pr:comm int}
Let $\calC$ be a Jordan curve with interior domain $D$, $\Omega_1,
\Omega_2\subset\C$ two open subsets with $\ol{D}\subset\Omega_2$, and $f:
\Omega_1\times\Omega_2\to\End(V)$ a holomorphic function such that $[\xi_i
(u),f(u,v)]=0$ for any $u,v$. Then, the following holds for any $\eps,\eta\in\{\pm\}$.
\\[-2ex]
\begin{enumerate}
\item If $u\notin\ol{D}\pm\eps a$,
\[\Ad(\xi_i(u))^{\pm 1}\ointC f(u,v)\,x^\eps_j(v)dv=
\ointC\left(\frac{u-v+\eps a}{u-v-\eps a}\right)^{\pm 1}f(u,v)\,x^\eps_j(v)dv\]\\[-2.2ex]
\item If $u\notin\ol{D}\pm\eps a-\eta\nN$
\[\Ad(g_i^\eta(u))^{\pm 1}\ointC f(u,v)\,x^\eps_j(v)dv=
\ointC\left(\frac{\Gamma^\eta(u-v-\eps a)}{\Gamma^\eta(u-v+\eps a)}\right)^{\pm 1}
f(u,v)\,x^\eps_j(v)dv\]\\[-2.2ex]
\item If $u\notin\ol{D}\pm\eps a+\Z$
\[\Ad(S_i(u))^{\pm 1}\ointC f(u,v)\,x^\eps_j(v)dv=
\ointC\left(\frac{q_i^{\eps a_{ij}}z-w}{z-q_i^{\eps a_{ij}}w}\right)^{\pm 1}f(u,v)\,x^\eps_j(v)dv\]
where $z=e^{2\pi\iota u}$ and $w=e^{2\pi\iota v}$.\\[-2.2ex]
\end{enumerate}
\end{prop}
\begin{pf}
(i) follows by left multiplying ($\Y$3) by $f(u,v)$ and integrating along $\calC$.

(ii) A repeated application of (i) shows that, for any $n\geq 1$
\begin{multline*}
\Ad(\xi_i(u+\eta n)\cdots\xi_i(u+\eta 1))^{\pm 1}\ointC f(u,v)\,x^\eps_j(v)dv\\
=\ointC\prod_{m=1}^n\left(\frac{u+\eta m-v+\eps a}{u+\eta m-v-\eps a}\right)^{\pm 1}
f(u,v)\,x^\eps_j(v)dv
\end{multline*}
so long as $u\notin\ol{D}\pm\eps a-\{\eta 1,\ldots,\eta n\}$. The result then follows
from \eqref{eq:gi-}--\eqref{eq:gi+}, $[\hbar\xi_{i,0},x_j^\eps(v)]=\eps2ax_j^\eps(v)$
and the fact that, for any $\alpha,\beta\in\C$
\[e^{\gamma(\alpha-\beta)}\prod_{m=1}^n\frac{m+\alpha}{m+\beta}e^{-\frac{\alpha-\beta}{m}}
\longrightarrow
\frac{\beta}{\alpha}\frac{\Gamma(\beta)}{\Gamma(\alpha)}
=\frac{\Gamma^+(\beta)}{\Gamma^+(\alpha)}
=\frac{\Gamma^-(-\beta)}{\Gamma^-(-\alpha)} \qquad
\]
as $n\to\infty$.
(iii) follows from (i)--(ii), the fact that
\begin{equation}\label{eq:Si}
S_i(u)=g_i^+(u)\xi_i(u)g_i^-(u)
\end{equation}
and \eqref{eq:sine'}.
\end{pf}

\subsection{Proof of (QL3)}\label{ssec: pf-QL23}
%--------------------------------

We shall check that
\begin{equation}\label{eq:QL3 to prove}
z\left(\Ad(\Psi_i(z))\X_{j,\ell}^\pm-q_i^{\pm a_{ij}}\X_{j,\ell}^\pm\right)
=q_i^{\pm a_{ij}}\Ad(\Psi_i(z))\X_{j,\ell+1}^\pm-\X_{j,\ell+1}^\pm
\end{equation}
The relation (QL3) then follows by expanding around $z=\infty$
and $z=0$.

Since \eqref{eq:QL3 to prove} is an equality of rational functions,
it suffices to prove it when $z$ lies outside of $\exp(2\pi\iota\calC
^\pm_j)$. Since $g_j^\pm$ is holomorphic inside $\calC^\pm_j$ by
Remark \ref{rem: independence}, Proposition \ref{pr:comm int}
implies that, for any $k\in\Z$,
\[\Ad(\Psi_i(z))\X_{j,k}^\pm= c_j^{\pm}
\ointCi{j}{\pm}w^k\frac{q_i^{\pm a_{ij}}z-w}{z-q_i^{\pm a_{ij}}w}g^\pm_j(v)x^\pm_j(v)dv\]
where $w=\exp(2\pi\iota v)$. The \lhs of \eqref{eq:QL3 to prove} is therefore
equal to
\[c_j^{\pm}\ointCi{j}{\pm}zw^\ell\left(\frac{q_i^{\pm a_{ij}}z-w}{z-q_i^{\pm a_{ij}}w}-q_i^{\pm a_{ij}}\right)g^\pm_j(v)x^\pm_j(v)dv\]
and the \rhs to
\[c_j^{\pm}\ointCi{j}{\pm}
w^{\ell+1}
\left(q_i^{\pm a_{ij}}\frac{q_i^{\pm a_{ij}}z-w}{z-q_i^{\pm a_{ij}}w}-1\right)g^\pm_j(v)x^\pm_j(v)dv\]

\subsection{Proof of (QL4)}\label{ssec: pf-QL4}
%--------------------------------

We shall prove (QL4) in the form
\begin{equation}\label{eq:XX}
\X^\pm_{i,k+1}\X^\pm_{j,\ell}-q_i^{\pm a_{ij}}\X^\pm_{i,k}\X^\pm_{j,\ell+1} = 
q_i^{\pm a_{ij}}\X^\pm_{j,\ell}\X^\pm_{i,k+1} - \X^\pm_{j,\ell+1}\X^\pm_{i,k}
\end{equation}
Set $a=\hbar d_ia_{ij}/2$, so that $\exp(2\pi\iota a)=q_i^{a_{ij}}$. By Proposition
\ref{pr:comm int},
\begin{multline}\label{eq:pre cancellation}
\ointCi{i}{\pm}
\left(e^{2\pi\iota(k+1)u}e^{2\pi\iota\ell v}-q_i^{\pm a_{ij}}e^{2\pi\iota k u}e^{2\pi\iota(\ell+1)v}\right)
g_i^\pm(u)x_i^\pm(u)du\cdot g_j^\pm(v)\\
=
g_j^\pm(v)\cdot\ointCi{i}{\pm}
e^{2\pi\iota(ku+\ell v)}\left(e^{2\pi\iota u}-e^{2\pi\iota(v\pm a)}\right)
\frac{\Gamma^\pm(v-u\pm a)}{\Gamma^\pm(v-u\mp a)}g_i^\pm(u)x_i^\pm(u)du
\end{multline}
so long as $v\notin\ol{D}^\pm_i\mp a\mp\nN$. As a function of $v$, the integrand
on the \rhs of \eqref{eq:pre cancellation} is regular on $\calC^\pm_i\mp a\mp\nN$
since the (simple) poles of $\Gamma^\pm(v-u\pm a)$ are cancelled by the zeros
of $\left(e^{2\pi\iota u}-e^{2\pi\iota(v\pm a)}\right)$. It follows that \eqref
{eq:pre cancellation} holds for any $v$.

Right multiplying the above by $x_j^\pm(v)$ and integrating along $v\in\calC^\pm_j$
then shows that the \lhs of \eqref{eq:XX} is equal to
\begin{multline*}
c_i^{\pm}c_j^{\pm}\ointCi{j}{\pm}\ointCi{i}{\pm} e^{2\pi\iota(ku+\ell v)}
\left(e^{2\pi\iota u}-e^{2\pi\iota(v\pm a)}\right)\\
\cdot\frac{\Gamma^\pm(v-u\pm a)}{\Gamma^\pm(v-u\mp a)}g_i^\pm(u)g_j^\pm(v) x_i^\pm(u)x^\pm_j(v)dudv
\end{multline*}
Similarly, the \rhs is equal to
\begin{multline*}
c_i^{\pm}c_j^{\pm}\ointCi{j}{\pm}\ointCi{i}{\pm}
e^{2\pi\iota(ku+\ell v)}\left(e^{2\pi\iota(u\pm a)}-e^{2\pi\iota v}\right)\\
\cdot\frac{\Gamma^\pm(u-v\pm a)}{\Gamma^\pm(u-v\mp a)}g_i^\pm(u)g_j^\pm(v) x^\pm_j(v)x_i^\pm(u)dudv
\end{multline*}

By ($\Y$4), the latter expression may be rewritten as
\begin{multline*}
c_i^{\pm}c_j^{\pm}\ointCi{j}{\pm}\ointCi{i}{\pm}
e^{2\pi\iota(ku+\ell v)}\left(e^{2\pi\iota(u\pm a)}-e^{2\pi\iota v}\right)\\
\cdot\frac{v-u\pm a}{v-u\mp a}\frac{\Gamma^\pm(u-v\pm a)}{\Gamma^\pm(u-v\mp a)}
g_i^\pm(u)g_j^\pm(v) x_i^\pm(u)x^\pm_j(v)dudv
\end{multline*}
since the contribution of the terms $[x_{i,0}^\pm,x_j^\pm(v)]$ and $[x_i^\pm(u),
x_{j,0}^\pm(v)]$ to the double integral is zero given that $(\Exp{(u\pm a)}-\Exp
{v})\Gamma^\pm(u-v\pm a)$ is holomorphic in $D_i\times D_j$.

(QL4) now follows from the fact that $\Gamma^+(x)=\Gamma^{-}(-x)$ and
\eqref{eq:sine'}.

\subsection{}
%--------------

We shall need the following result to show that relation (QL5) holds.

\begin{lem}\label{le:double}
Let $\calC\subset\C$ be a Jordan curve with interior domain $D$, $\xi(u)$
a meromorphic function with no singularities on $\calC$, and $f(u,v)$ a
meromorphic function with no singularities in $\ol{D}\times\ol{D}$. Then,
\[\ointC\ointC f(u,v)\frac{\xi(u)-\xi(v)}{u-v}du\,dv=
-\ointC f(u,u)\xi(u)du\]
\end{lem}
\begin{pf}
Let $\calC_<$ be a small deformation of $\calC$ contained in $\calC$. Then,
\[\begin{split}
\ointC\ointC f(u,v)\frac{\xi(u)-\xi(v)}{u-v}du\,dv
\phantom{-}
&=\phantom{-}
\ointC\left(\oint_{\calC_<} f(u,v)\frac{\xi(u)-\xi(v)}{u-v}du\right)\,dv\\
&=\phantom{-}
\ointC\left(\oint_{\calC_<} f(u,v)\frac{\xi(u)}{u-v}du\right)\,dv\\
&=\phantom{-}
\oint_{\calC_<}\xi(u)\left(\ointC \frac{f(u,v)}{u-v}dv\right)du\\
&=
-\oint_{\calC_<}\xi(u)f(u,u)du
\end{split}\]
\end{pf}

\subsection{Proof of (QL5)}\label{ssec: pf-QL5}
%-------------------------------

We have
\[\begin{split}
[\X^+_{i,k},\X^-_{j,\ell}]
&= c_i^+c_j^-
\ointCi{i}{+}\ointCi{j}{-}\Exp{k u}\Exp{\ell v}g_i^+(u)x_i^+(u)g_j^-(v)x_j^-(v)dudv\\
&-c_i^+c_j^-
\ointCi{j}{-}\ointCi{i}{+}\Exp{k u}\Exp{\ell v}g_j^-(v)x_j^-(v)g_i^+(u)x_i^+(u)dudv
\end{split}\]
We wish to apply Proposition \ref{pr:comm int} to permute $x_i^+(u)$ and
$g_j^-(v)$ in the first integral, and $x_j^-(v)$ and $g_i^+(u)$ in the second.
This cannot be done directly however, since $\calC^-_j\ni v$ may lie inside
$\ol{D}^+_i-a+\nN$, where $a=\hbar d_ia_{ij}/2$, and $\calC^+_i\ni u$
may lie inside $\ol{D}^-_j+a-\nN$.

To circumvent this issue, let $\delta\in\C$ and consider the shifted integral
\begin{multline*}
I(\delta)=c_i^+c_j^-
\ointCi{i}{+}\ointCi{j}{-}\Exp{k u}\Exp{\ell v}g_i^+(u+\delta)x_i^+(u)g_j^-(v-\delta)x_j^-(v)dudv\\
-c_i^+c_j^-
\ointCi{j}{-}\ointCi{i}{+}\Exp{k u}\Exp{\ell v}g_j^-(v-\delta)x_j^-(v)g_i^+(u+\delta)x_i^+(u)dvdu
\end{multline*}
$I(\delta)$ is holomorphic in a disk $|\delta|<R$, where $R$ is such that $g_i
^+$ and $g_j^-$ are holomorphic on $\ol{D}^+_i+\delta'$ and $\ol{D}^-_j-\delta'$
for any $\delta'$ such that $|\delta'|<R$ respectively. Moreover, if $\calC^+_i,
\calC^-_j$ are small enough, there is an $r<R$ such that if $|\delta|>r$, $\ol{D}
^-_j-\delta$ is disjoint from $\ol{D}^+_i-a+\nN$, and $\ol{D}^+_i+\delta$ is
disjoint from $\ol{D}^-_j+a-\nN$.

Assuming that $r<|\delta|<R$, we may now apply Proposition \ref{pr:comm int}
to find that $I(\delta)$ is equal to
\begin{multline*}
c_i^+c_j^-\ointCi{i}{+}\ointCi{j}{-}\Exp{k u}\Exp{\ell v}
g_i^+(u+\delta)g_j^-(v-\delta)
\frac{\Gamma^-(v-\delta-u+a)}{\Gamma^-(v-\delta-u-a)}x_i^+(u)x_j^-(v)dudv\\
-c_i^+c_j^- 
\ointCi{j}{-}\ointCi{i}{+}\Exp{k u}
\Exp{\ell v}g_j^-(v-\delta)g_i^+(u+\delta)
\frac{\Gamma^+(u+\delta-v-a)}{\Gamma^+(u+\delta-v+a)}x_j^-(v)x_i^+(u)dvdu
\end{multline*}
Since $\Gamma^+(x)=\Gamma^{-}(x)$, this yields
\[\begin{split}
\frac{I(\delta)}{c_i^+c_j^-}
&=
\ointCi{i}{+}\ointCi{j}{-}\Exp{k u}\Exp{\ell v}
g_i^+(u+\delta)g_j^-(v-\delta)
\frac{\Gamma^+(u+\delta-v-a)}{\Gamma^+(u+\delta-v+a)}[x_i^+(u),x_j^-(v)]dudv\\
&=
-\hbar\delta_{ij}\ointCi{i}{+}\ointCi{i}{-}\Exp{k u}\Exp{\ell v}
g_i^+(u+\delta)g_i^-(v-\delta)
\frac{\Gamma^+(u+\delta-v-\hbar d_i)}{\Gamma^+(u+\delta-v+\hbar d_i)}
\frac{\xi_i(u)-\xi_i(v)}{u-v}dudv
\end{split}\]
where we used ($\Y$5).

Since the function $\Gamma^+(u+\delta-v-\hbar d_i)$ is holomorphic on
$\ol{D}^+_i\times\ol{D}^-_i$ for $|\delta|>r$, the only singularities of the
integrand are the poles of $\xi_i(u),\xi_i(v)$. Choosing the contours $\calC
_i^\pm$ to be equal to a common contour $\calC_i$ which contains all
poles of  $\xi_i$, Lemma \ref{le:double} then yields
\begin{equation}\label{eq:almost}
I(\delta)
=c_i^+c_j^-
\hbar\delta_{ij}\ointCi{i}{}\Exp{(k+\ell)u}
g_i^+(u+\delta)g_i^-(u-\delta)
\frac{\Gamma^+(\delta-\hbar d_i)}{\Gamma^+(\delta+\hbar d_i)}\xi_i(u)du
\end{equation}
Since both $I(\delta)$ and the \rhs $J(\delta)$ of \eqref{eq:almost} are holomorphic
functions on the disk $|\delta|<R$, $I(0)$ is equal to
\[\begin{split}
J(0)
&=\delta_{ij}\hbar c_i^+c_i^-\frac{\Gamma(1-d_i\hbar)}{\Gamma(1+d_i\hbar)}
\ointCi{i}{}\Exp{(k+\ell)u}g_i^+(u)g_i^-(u)\xi_i(u)du\\
&=\delta_{ij}\frac{2\pi\iota}{q_i-q_i^{-1}}\ointCi{i}{} \Exp{(k+\ell)u}S_i(u)du\\
&=\frac{\delta_{ij}}{q_i-q_i^{-1}}\oint_{\wt{\calC}_i} z^{k+\ell-1}\Psi_i(z)dz\\
&=\frac{\delta_{ij}}{q_i-q_i^{-1}}\left(\Psi^+_{i,k+\ell}-\Psi^-_{i,k+\ell}\right)
\end{split}\]
where the second equality follows from \eqref{eq:Si}, \eqref{eq: constants}
and \eqref{eq:sine'}, the third from the change of variables $z=\exp(2\pi\iota u)$
\footnote{we are assuming $\calC_i$ to be small enough that the exponential
map is injective on each of the connected components of its interior domain
$D_i$}, and the fourth from the fact that $\wt{\calC}_i=\exp(2\pi\iota\calC_i)$
encloses all the poles of $\Psi_i$ by \ref{ssec: contours} (i),
and neither $z=0$ nor $z=\infty$.

\subsection{Compatibility with shift automorphisms}
\label{ssec: functor-shift}
%--------------------------------------------------------------

Let $\rho:\Yhg\to\End(V)$ be the action of $\Yhg$ on $V$ and
$\rho_a : \Yhg \to \End(V)$ the pull--back of $\rho$ via the shift
automorphism $\tau_a$ of $\Yhg$. Thus, for $y = \xi_i, x_i^{\pm}$
we have
\[\rho_a(y(u)) = \rho(y(u-a))\]

We need to show that for $Y=\Psi_i, \X_i^{\pm}$, the following
holds on $V$
\[\Fh{}(\rho_a)(Y(z)) = \Fh{}(\rho)(Y(e^{-2\pi\iota a}z))\]
For notational convenience, we denote the \lhs by $\wt{Y}(z)$
and the \rhs by $Y(e^{-2\pi\iota a}z)$. By definition, $\wt{\Psi_i}(z)$
is obtained from the monodromy of the difference equation
\[\wt{\phi_i}(u+1) 
= \xi_i(u-a)\wt{\phi_i}(u)\]
Clearly if $\phi_i^{\pm}(u)$ and $S_i(u)$ are the fundamental solutions
and the connection matrix of $\phi_i(u+1) = \xi_i(u)\phi_i(u)$, then
$\phi_i^{\pm}(u-a), S_i(u-a)$ are those of the difference equation
$\wt{\phi_i}(u+1) = \xi_i(u-a)\wt{\phi_i}(u)$. Therefore, we get
\[
\wt{\Psi_i}(z) = \left. S_i(u-a)\right|_{z=\Exp{u}}
= \left. S_i(u)\right|_{z=\Exp{(u+a)}} = 
\Psi_i(e^{-2\pi\iota a}z)
\]

Next, if $\calC$ is the contour used to define $\X_i^{\pm}(z)$, then
$\calC+a$ can be taken as the contour for $\wt{\X_i^{\pm}}(z)$. It
follows that
\begin{align*}
\wt{\X_i^{\pm}}(z) &= c_i^{\pm}\oint_{\calC+a} 
\frac{z}{z-\Exp{u}} g_i^{\pm}(u-a)x_i^{\pm}(u-a)du \\
&= c_i^{\pm}\oint_{\calC} \frac{z}{z-\Exp{(u+a)}} 
g_i^{\pm}(u)x_i^{\pm}(u)du\\
&= \X_i^{\pm}(e^{-2\pi\iota a}z)
\end{align*}

\section{The inverse functor $\Gh{}$}\label{sec: inverse-functor}
%==========================

Let $\Omega\subset\C^\times$ be a subset stable under multiplication
by $q$, and $\Pi\subset\C$ a non--congruent subset stable under
shifts by $\frac{1}{2}\Z\hbar$ such that $\Omega=\exp(2\pi\iota\Pi)$. In
this section, we construct a functor 
\[\Gh{\Pi}:\Rlsub\to\Rysub\]
which is an inverse to the functor $\Fh{\Pi}$ given by Theorem
\ref{thm: second-main-theorem}.

\subsection{Inverse monodromy problem}\label{ss:inv-mon-qloop}
%---------------------------------------------------

Let $\V\in\Rlsub$, and let $\mu\in\h^*$ be a weight of $\V$. For any
$i\in\bfI$, Proposition \ref{prop: rationality} yields a $GL(\V_\mu)$--valued
rational function $\Psi_i(z)_{\mu}$ satisfying
\[\Psi_i(\infty)_\mu=q_i^{\mu(\alpha_i^\vee)}=\Psi_i(0)^{-1}_\mu\]
and such that the poles of $\Psi_i(z)^{\pm 1}_{\mu}$ are contained
in $\Omega$. Set  $A_0=\hbar d_i\mu(\alpha_i^\vee)$ so that $\Psi_i(\infty)
_\mu=\exp(\pi\iota A_0)$, and $S_i(u)=\Psi\lp\Exp{u}\rp_{\mu}$.

\begin{lem}\label{lem:inv-mon-qloop}
There exists a unique rational function $A(u):\C\to GL(\V_\mu)$ such
that
\begin{itemize}
\item $A(u)=1+A_0 u^{-1}+\cdots$.
\item $[A(u),A(v)]=0$ for any $u,v$.
\item the poles of $A(u)^{\pm 1}$ are contained in $\Pi$.
\item $S_i(u)$ is the connection matrix of $\phi(u+1) = A(u)\phi(u)$.
\end{itemize}
\end{lem}

\begin{pf} By Theorem \ref{ss:start factorisation}, the inverse monodromy
problem may be solved by making a choice of logarithms of the zeros and
poles of the eigenvalues of $\Psi_i(z)_\mu$,  and of the poles of the unipotent
part of $\Psi_i(z)_\mu$. Choosing the logarithms to lie in $\Pi$ leads to
a unique coefficient matrix $A(u)$ since $\Pi$ is non--congruent. We do
however need to check that the choices are consistent in the sense of
equation \eqref{eq:consistent choice}.

Let $\gamma(z)$ be an eigenvalue of $\Psi_i(z)$ on the weight space
$\V_{\mu}$. By \cite{frenkel-reshetikhin-qchar},\cite[Prop. 16]{hernandez-affinizations},
there are two monic polynomials $Q(z),R(z)$, not vanishing at $z=0$,
and such that
\[\gamma(z)=
q_i^{\deg(R)-\deg(Q)}
\frac{Q(q_i^2z)}{Q(z)}\frac{R(z)}{R(q_i^2z)}\]
Let $Q(z)=\prod_j(z-\beta_j)$
and $R(z)=\prod_k(z-\beta_k')$. Let $b_j,b_k'\in\Pi$ be such that $\Exp
{b_j}=\beta_j$ and $\Exp{b_k'} = \beta_k'$. Since $\Pi$ is stable under
shifts by $\hbar/2$, $b_j-d_i\hbar,b'_k-
d_i\hbar$ are logarithms of $q_i^{-2}\beta_j$ and $q_i^{-2}\beta'_k$
respectively, which lie in $\Pi$. Since $A_0=\hbar d_i\mu(\alpha_i^\vee)$,
the consistency requirement \eqref{eq:consistent choice} reads
\[\sum_j b_j + \sum_k (b'_k-d_i\hbar) - \sum_j (b_j-d_i\hbar) - \sum_k
b'_k=\hbar d_i\mu(\alpha_i^\vee)\]
The latter holds since $\Psi_i(\infty)$ acts on $\V_\mu$ as multiplication
by $q_i^{\mu(\alpha_i^\vee)}$, $\gamma(\infty)=q_i^{\deg Q-\deg R}$,
and $q$ is not a root of unity. 
\end{pf}

\subsection{Yangian operators}\label{ssec: y-operators}
%-------------------------------------

Given $\V\in\Rlsub$, and a weight $\mu\in\h^*$ of $\V$, define
the action of the generators of $\Yhg$ on $\V_\mu$ as follows.

\begin{enumerate}
\item The action of $\h$ on $\V_{\mu}$ is by weight $\mu$. Note
that $K_h = q^{\mu(h)}$ on $\V_{\mu}$ for every $h\in\h$ determines
$\mu\in\h^*$ uniquely.
\item $\xi_i(u)_{\mu}$ acts as the coefficient matrix $A(u)$ given
by Lemma \ref{lem:inv-mon-qloop}.
\end{enumerate}
Let $\phi_i^{\pm}(u)_{\mu}$ be the fundamental solutions of the
difference equation determined by $\xi_i(u)_\mu$, and define $
g_i^{\pm}(u)_{\mu}$ by \eqref{eq:gi+-}. Let $\calC^{\pm}_{i,\mu}\subset
\C$ a Jordan curve such that
\begin{enumerate}
\item[(a)] $\calC^\pm_{i,\mu}$ encloses all the poles of $\X_i^\pm(\Exp{u})$
contained in $\Pi$, and none of their $\nZ$--translates.
\item[(b)] $\calC^\pm_{i,\mu}$ encloses none of the $\nZ$--translates
of the poles of $\xi_i(u)_\mu^{-1}$.
\end{enumerate}
Note that such a curve exists since the poles of $\xi_i(u)_\mu^
{-1}$ are contained in $\Pi$, and the latter is non--congruent.

\begin{enumerate}
\item[(iii)] The generators $x_{i,r}^\pm$ act on $\V_\mu$ by
\begin{equation}\label{eq: i-fun-x+r}
\lp x_{i,r}^{\pm}\rp_{\mu}=
\frac{1}{c_i^{\pm}\hbar}\int_{\calC_{i,\mu}^\pm}
v^r g_i^{\pm}(v)_{\mu\pm\alpha_i}^{-1} 
\X^{\pm}_i\lp\Exp{v}\rp_{\mu}\, dv 
\end{equation}
Note that these operators are independent of the choice of the contour
$\calC_{i,\mu}^{\pm}$ since, by assumption (b) and \eqref{eq:gi-}--\eqref{eq:gi+},
$g_i^{\pm}(v)_{\mu\pm\alpha_i}^{-1}$ are holomorphic inside $\calC
_{i,\mu}^{\pm}$. In terms of the fields $x_i^\pm
(u)$, \eqref{eq: i-fun-x+r} may be rewritten as
\begin{gather}
x_i^{\pm}(u)_{\mu} = \frac{1}{c_i^{\pm}}\int_{\calC_{i,\mu}^\pm}\frac{1}{u-v}
g_i^{\pm}(v)_{\mu\pm\alpha_i}^{-1}
\X^{\pm}_i\lp\Exp{v}\rp_{\mu}\, dv \label{eq: i-fun-x+} 
\end{gather}
for $u$ not enclosed by $\calC_{i,\mu}^{\pm}$.
\end{enumerate}

\vspace*{-0.2cm}

\subsection{}
%--------------

\begin{thm}\label{thm: third-main-theorem}\hfill\break
\begin{enumerate}
\item The above operators define an action of $\Yhg$ on $\V$.
\item As a representation of $\Yhg$, $\V$ lies in $\Rysub$.
\item The corresponding functor $\Gh{\Pi}:\Rlsub\to\Rysub$
satisfies
\[\Fh{\Pi}\circ\Gh{\Pi}=\id_{\Rlsub}\aand
\Gh{\Pi}\circ\Fh{\Pi}=\id_{\Rysub}\]
In particular, if $\Pi$ is a fundamental domain for $u\mapsto u+1$,
$\Fh{\Pi}$ gives rise to an isomorphism of categories $\Rysub\isom
\Rloop$.
\item The functor $\Gh{\Pi}$ is compatible with shift automorphisms.
That is, for any $\V\in\Rlsub$, and $a\in\C$
\[\Gh{\Pi+a}(\V(\Exp{a}))=\Gh{\Pi}(\V)(a)\]
\end{enumerate}
\end{thm}
\begin{pf}
(i) The relations (Y0), (Y1) and (Y2) hold trivially. After some preparatory
work in \ref{ss:lin indep}--\ref{ss:converse comm int}, we prove (Y3)
in \S \ref{ssec: pf-Y23}, (Y4) in \S \ref{ssec: pf-Y4} and (Y5) in \S \ref
{ssec: pf-Y5}. Relation (Y6) follows from Proposition \ref{prop: Y6-reduction}.

(ii) By construction, $\{\xi_i(u)\}_{i\in\bfI}$ have zeros and poles in
$\Pi$. It follows by Theorem \ref{thm: first-main-theorem} part (ii)
that $\Gh{\Pi}(\V)$ lies in $\Rysub$.

(iii) The fact that $\Fh{\Pi}\circ\Gh{\Pi}(\V)=\V$ as modules over the
subalgebra of $\qloop$ generated by $\{\Psi^\pm_{i,\pm r}\}_{i\in\bfI,r\in\N}$
follows from the definition of $\Fh{}$, and the fact that Lemma \ref
{lem:inv-mon-qloop} solves an inverse monodromy problem. The
fact that $\Gh{\Pi}\circ\Fh{\Pi}(V)=V$ as modules over the
commutative subalgebra of $\Yhg$ generated by $\{\xi_{i,r}\}_{i\in
\bfI,r\in\N}$ follows from the uniqueness statement of Lemma \ref
{lem:inv-mon-qloop}. The corresponding assertions for the generators
$x_{i,r}^\pm,\X_{i,k}^\pm$ of $\Yhg,\qloop$ are checked in \ref
{ssec: pf-inv} and \ref{ssec: pf'-inv} respectively.

(iv) Follows by a direct check similar to that in \ref{ssec: functor-shift},
or by applying $\Fh{\Pi+a}$ and using Theorem \ref{thm: second-main-theorem}
(ii).
\end{pf}

\subsection{}
%--------------

\begin{rem}\label{rem: zorn}
We observe that one can always find $\Pi\subset\C$ which is stable under
shifts by $\frac{\hbar}{2}$ and is a fundamental domain for $u\mapsto u+1$,
assuming that $\hbar$ is not rational. If $\Im(\hbar)\neq 0$, $\Pi$ may be
taken to be the strip bounded between $\R\hbar$ and $\R\hbar + 1$ including
$\R\hbar$.

If $\hbar\in\R\setminus\Q$, we claim that there exists $U\subset\R$ such that
$U\pm\frac{\hbar}{2}\subset U$ and $u\mapsto \Exp{u}$ is a bijection between
$U$ and $S^1$. Assuming this, take $\Pi = \{u+\iota t: u\in U, t\in \R\}$. To prove
the claim, we shall use Zorn's lemma as follows. Let $\mathcal{U}$ to be the
set of all subsets $V \subset\R$ such that $V$ is stable under shifts
by $\hbar/2$ and $u\mapsto \Exp{u}$ is an injective map from $V$ to $S^1$. The
set $\mathcal{U}$ is non--empty, since $\Z\frac{\hbar}{2}$ is an element of
$\mathcal{U}$. Moreover every chain in $\mathcal{U}$ has the largest element
in $\mathcal{U}$ and hence by Zorn's lemma $\mathcal{U}$ has a maximal
element $U$. We claim that $u\mapsto \Exp{u}$ is a bijection between
$U$ and $S^1$. If not, then there exists $x\in \R$ such that $x\not\equiv
y \text{ (mod } \Z\text{)} $ for every $y\in U$. The same is true for every $x+n\hbar$
($n\in \frac{1}{2}\Z$) and hence $U \cup \{x+ n\hbar\}_{n\in (1/2)\Z}$ is also 
in $\mathcal{U}$, contradicting the maximality of $U$.
\end{rem}

\subsection{}\label{ss:lin indep}
%--------------

\begin{lem}\label{le:lin indep}
Let $\calC$ be a Jordan curve the interior domain $D$ of which
is non--congruent, and $f$ a meromorphic function on a neighborhood
of $\ol{D}$ with no singularities on $\calC$.
Let $a\in\C\setminus\Z$, and assume that
\[\intC\frac{\Exp{(u+a)}-\Exp{v}}{\Exp{u}-\Exp{(v+a)}}f(v)dv=0\]
for any $u\notin\ol{D}+a+\Z$. Then, $f$ has no poles in $D$.
\end{lem}
\begin{pf}
Expanding the above integral near $\Exp{u}=\infty$ shows that
\[\intC\Exp{k v}f(v)dv=0\]
for any $k\in\N$, and therefore that $\intC p(\Exp{v})f(v)dv=0$
for any polynomial $p$. Let us assume on the contrary that
 $\{b_j\}_{j\in J}$ are the poles of $f$
in $D$, and $\{n_j\}_{j\in J}$ their orders. For any $j\in J$, set
\[p_j(w)=(w-\beta_j)^{n_j-1}\prod_{j'\neq j}(w-\beta_{j'})^{n_{j'}}\]
where $\beta_k=\Exp{b_k}$. Then $p_j(\Exp{v})f(v)$ is holomorphic
on $D\setminus\{b_j\}$ and has a pole of order 1 at $v=b_j$ since
$D$ is non--congruent. It follows that $\intC p_j(\Exp{v})f(v)dv\neq 0$,
a contradiction.
\end{pf}

\subsection{}\label{ss:converse comm int}
%--------------

Let $i,j\in\bfI$, and set $a=\hbar d_ia_{ij}/2$. The following is a 
converse to Proposition \ref{pr:comm int}.

\begin{prop}\label{pr:converse comm int}
Let $\calC$ be a Jordan curve with interior domain $D$, and $f:\C\to
\End(\V)$ a meromorphic function with no singularities in $\ol{D}$,
and such that $[\Psi_i(\Exp{u}),f(v)]=0$ for any $u,v$. Then, the
following holds for any $\eps\in\{\pm\}$.\\[-2ex]
\begin{enumerate}
\item If $u\notin\ol{D}\pm\eps a+\Z$
\begin{multline*}
\Ad(\Psi_i(\Exp{u}))^{\pm 1}\intC f(v)\,\X^\eps_j(\Exp{v})dv\\
=
\intC\left(\frac{\Exp{(u+\eps a)}-\Exp{v}}{\Exp{u}-\Exp{(v+\eps a)}}\right)^{\pm 1}
f(v)\,\X^\eps_j(\Exp{v})dv
\end{multline*}
\end{enumerate}
Assume moreover that $D$ is such that
\begin{enumerate}
\item[(1)] $D$ is non--congruent.
\item[(2)] $\ol{D}+a$ is non--congruent to $\ol{D}-a$.
\item[(3)] $\ol{D}\pm a$ are non--congruent to the set of poles of $\xi_i(u)^{\pm 1}$.
\end{enumerate}
Then, the following holds for any $\eps,\eta\in\{\pm\}$.
\begin{enumerate}
\item[(ii)] If $u\notin\ol{D}\pm\eps a-\eta\nN$
\[\Ad(g_i^\eta(u))^{\pm 1}\intC f(v)\,\X^\eps_j(\Exp{v})dv=
\intC\left(\frac{\Gamma^\eta(u-v-\eps a)}{\Gamma^\eta(u-v+\eps a)}\right)^{\pm 1}
f(v)\,\X^\eps_j(\Exp{v})dv\]\\[-2.2ex]
\item[(iii)] If $u\notin\ol{D}\pm\eps a$,
\[\Ad(\xi_i(u))^{\pm 1}\intC f(v)\,\X^\eps_j(\Exp{v})dv=
\intC\left(\frac{u-v+\eps a}{u-v-\eps a}\right)^{\pm 1}f(v)\,\X^\eps_j(\Exp{v})dv\]
\end{enumerate}
\end{prop}
\begin{pf}
(i) follows by left multiplying (\QL3) by $f$ and integrating along $\calC$.

(ii)--(iii) Fix a weight $\mu$ of $\V$, and let $E^\eps\subseteq\Hom(\V_
\mu,\V_{\mu+\eps\alpha_j})$ be the subspace spanned by elements
of the form $\intC f(v)\X^\eps_j(\Exp{v})dv$, where $f(v)$ is meromorphic
with no singularities in $\ol{D}$, and commutes with $\Psi_i(\Exp{u})$ for
all $u,v$. By (i), $E^\eps$ is stable under $\Ad(\Psi_i
(\Exp{u}))$ for $u\notin\ol{D}+\eps a+\Z$, and therefore for any $u$ since
$\Psi_i(\Exp{u})$ is meromorphic in $u$. It follows that $E^\eps$ is stable
under $\Ad(\xi_i(u))$ and $\Ad(g_i^\eta(u))$. We now
wish to show that the action of $\Ad(g_i^\eta(u))$ and $\Ad(\xi_i(u))$ on
$E^\eps$ is given by the formulae (ii)--(iii). We shall do so by showing
that the formulae (ii) define meromorphic automorphisms of $E^\eps$
which factorise $\Ad(\Psi_i(\Exp{u}))$, and then relying on the uniqueness
statement of Proposition \ref{pr:unique factorisation}.

For $u\notin\ol{D}+\eps a-\eta\nN$, define $\calG^\eta(u)\in
GL(E^\eps)$ by
\[\intC f(v)\X^\eps_j(\Exp{v})dv\longrightarrow
\intC\frac{\Gamma^\eta(u-v-\eps a)}{\Gamma^\eta(u-v+\eps a)}
f(v)\X^\eps_j(\Exp{v})dv\]
Note that this is well--defined by (i), Lemma \ref{le:lin indep} and the
fact that $D$ is non--congruent. Since the \rhs is a meromorphic function
of $u\in\C$, and lies in
$E^\eps$ for $u\notin\ol{D}+\eps a-\eta\nN$, it lies in $E^\eps$ for
all $u$.

It is clear that $\calG^\eta(u)$ are invertible and holomorphic for 
$\Re(\pm u)\gg 0$, and that $\calG^+(u)^{-1}\calG^-(u)=\left.\Ad
(\Psi_i(\Exp{u}))\right|_{E^\eps}$. Set
\[\Xi(u)=\calG^+(u+1)\calG^+(u)^{-1}=\calG^-(u+1)\calG^-(u)^{-1}\]
By inspection, $\Xi(u)^{\pm 1}$ is given by
\[\intC f(v)\X^\eps_j(\Exp{v})dv\longrightarrow
\intC \left(\frac{u-v+\eps a}{u-v-\eps a}\right)^{\pm 1}f(v)\X^\eps_j(\Exp{v})dv\]
for any $u\notin\ol{D}\pm\eps a$. Since the poles of
$\Xi(u)^{\pm 1}$ are contained in $\ol{D}\pm\eps a$,
and those of $\left.\Ad(\xi_i(u))\right|_{E^\eps}$ are contained in
$\Pi$ which is non--congruent, it follows from our assumptions
on $D$ and Proposition \ref{pr:unique factorisation} that $\Xi(u)
=\left.(\Ad(\xi_i(u))\right|_{E^\eps}$, and therefore that $\calG^
\eta(u)=\left.(\Ad(g^\eta_i(u))\right|_{E^\eps}$.
\end{pf}

\subsection{Proof of (Y3)}\label{ssec: pf-Y23}
%------------------------------

It suffices to prove $(\Y3)$ for $u,v$ large enough. Let $i,j\in\bfI$,
and set $a=\hbar d_ia_{ij}/2$. Since the contour $\calC^\pm_j$ defining
$x^\pm_j(v)$ in \eqref{eq: i-fun-x+r} is only required to contain the
poles of $\xi_j(v)$, and the latter are contained in $\Pi$ which is
non--congruent and stable under shift by $\hbar/2$, we may choose
$\calC^\pm_j$ so that it also satisfies the conditions (1)--(3) of Proposition
\ref{pr:converse comm int}. It then follows that
\[\xi_i(u)x^\pm_j(v)\xi_i(u)^{-1}
= \frac{1}{c_j^{\pm}}
\intCi{j}{\pm}\frac{u-v'\pm a}{(v-v')(u-v'\mp a)}
g_j^\pm(v')^{-1}
\X^\pm_j\lp\Exp{v'}\rp\, dv'\]
so long as $u\mp a$ and $v$ do not lie within $\calC_j$. $(\Y3)$
now follows from the fact that
\[\frac{u-v'\pm a}{(v-v')(u-v'\mp a)}
=
\frac{u-v\pm a}{u-v\mp a}\,\frac{1}{v-v'}
\mp
\frac{2a}{u-v\mp a}\,\frac{1}{u\mp a-v'}\]

\subsection{Proof of (Y4)}\label{ssec: pf-Y4}
%------------------------------

We shall prove that
\begin{equation}
x_{i,r+1}^\pm x_{j,s}^\pm - x^\pm_{i,r}x^\pm_{j,s+1}\mp ax_{i,r}^\pm x_{j,s}^\pm = 
x_{j,s}^\pm x^\pm_{i,r+1} - x_{j,s+1}^\pm x^\pm_{i,r}\pm ax^\pm_{j,s}x^\pm_{i,r}
\label{eq:Y4 to prove}
\end{equation}

The left--hand side, multiplied by $c_i^{\pm}c_j^{\pm}$ is equal to
\[\frac{1}{\hbar^2}\intCi{i}{\pm}\intCi{j}{\pm}
u^rv^s(u-v\mp a) g_i^\pm(u)^{-1}
\X^\pm_i\lp\Exp{u}\rp g_j^\pm(v)^{-1}\X^\pm_j\lp\Exp{v}\rp\,du\,dv\]
Choose $\calC_i^\pm,\calC_j^\pm$ so that they satisfy conditions (1)--(3) of
Proposition \ref{pr:converse comm int} \wrt $\xi_j(v)$ and $\xi_i(u)$ respectively,
as well as
\[\calC^\pm_j\cap\lp\ol{D}_i\pm a\mp\nN\rp=
\emptyset=
\calC^\pm_i\cap\lp\ol{D}_j\pm a\mp\nN\rp\]
where $D_i,D_j$ are the interior domains of $\calC^\pm_i,\calC^\pm_j$. Then,
the above is equal to
\[\frac{1}{\hbar^2}\intCi{i}{\pm}\intCi{j}{\pm}
u^rv^s(u-v\mp a)\frac{\Gamma^\pm(v-u\mp a)}{\Gamma^\pm(v-u\pm a)}
g_i^\pm(u)^{-1}g_j^\pm(v)^{-1}\X^\pm_i\lp\Exp{u}\rp\X^\pm_j\lp\Exp{v}\rp\,du\,dv\]

Similarly, the \rhs of \eqref{eq:Y4 to prove} multiplied by $c_i^{\pm}c_j^{\pm}$ 
is equal to
\[\frac{1}{\hbar^2}\intCi{i}{\pm}\intCi{j}{\pm}
u^rv^s (u-v\pm a)
\frac{\Gamma^\pm(u-v\mp a)}{\Gamma^\pm(u-v\pm a)}
g_i^\pm(u)^{-1} g_j^\pm(v)^{-1} \X^\pm_j\lp\Exp{v}\rp \X^\pm_i\lp\Exp{u}\rp\,du\,dv\]
By (\QL4), this may be rewritten as
\begin{multline*}
\frac{1}{\hbar^2}\intCi{i}{\pm}\intCi{j}{\pm}
u^rv^s (u-v\pm a)
\frac{\Gamma^\pm(u-v\mp a)}{\Gamma^\pm(u-v\pm a)}
\frac{\Exp{(v\pm a)}-\Exp{u}}{\Exp{v}-\Exp{(u\pm a)}}
\\
\cdot g_i^\pm(u)^{-1} g_j^\pm(v)^{-1}\X^\pm_i\lp\Exp{u}\rp\X^\pm_j\lp\Exp{v}\rp\,du\,dv
\end{multline*}
Indeed, neither of the boundary terms
$\Exp{v}\lp\X_{j,0}^\pm\X_i^\pm(\Exp{u})-q_j^{\pm a_{ji}}\X_i^\pm(\Exp{u})\X_{j,0}^\pm\rp$
and
$\Exp{u}\lp\X_{i,0}^\pm\X_j^\pm(\Exp{v})-q_j^{\pm a_{ji}}\X_j^\pm(\Exp{v})\X_{i,0}^\pm\rp$
from (\QL4) contribute to the double integral since $(u-v\pm a)/(\Exp{v}-\Exp{(u\pm a)})$
does not have any poles inside $\calC^\pm_i\times\calC^\pm_j$.

The result now follows from the fact that $\Gamma^\mp(x)=\Gamma^\pm(-x)$,
and \eqref{eq:sine'}.

\subsection{Proof of (Y5)}\label{ssec: pf-Y5}
%------------------------------

We need to prove that for any $i,j\in \bfI$ and $r,s\in \N$ we have
\[[x_{i,r}^+, x_{j,s}^-] = \delta_{ij} \xi_{i,r+s}\]
Choosing the contours $\calC^+_i,\calC^-_j$ as in \ref{ssec: pf-Y4}, 
Proposition \ref{pr:converse comm int} yields
\begin{align*}
x_{i,r}^+x_{j,s}^- &=
\frac{1}{c_i^+c_j^-\hbar^2} 
\intCi{i}{+}\intCi{j}{-} u^rv^s g_i^+(u)^{-1}\X^+_i\lp\Exp{u}\rp 
g_j^-(v)^{-1}\X^-_j\lp\Exp{v}\rp\, du\, dv \\
&= \frac{1}{c_i^+c_j^-\hbar^2} \intCi{i}{+}\intCi{j}{-} u^rv^s g_i^+(u)^{-1}g_j^-(v)^{-1}
\frac{\Gamma^-(v-u-a)}{\Gamma^-(v-u+a)}
\X^+_i\lp\Exp{u}\rp \X^-_j\lp\Exp{v}\rp\, du\, dv
\end{align*}
Similarly,
\[x_{j,s}^-x_{i,r}^+ =\frac{1}{c_i^+c_j^-\hbar^2}\intCi{i}{+}\intCi{j}{-}
u^rv^s g_i^+(u)^{-1}g_j^-(v)^{-1}
\frac{\Gamma^+(u-v+a)}{\Gamma^+(u-v-a)}
\X^-_j\lp\Exp{v}\rp \X^+_i\lp\Exp{u}\rp\, du\, dv\]
Since $\Gamma^+(x) = \Gamma^-(-x)$ it follows that
\begin{multline*}
[x_{i,r}^+, x_{j,s}^-]=
\frac{1}{c_i^+c_j^-\hbar^2} 
\intCi{i}{+}\intCi{j}{-} u^rv^s g_i^+(u)^{-1}g_j^-(v)^{-1}\\
\cdot\frac{\Gamma^+(u-v+a)}{\Gamma^+(u-v-a)}
[\X^+_i\lp\Exp{u}\rp,\X^-_j\lp\Exp{v}\rp]\, du\, dv 
\end{multline*}

By (\QL 5), we get
\begin{multline*}
[x_{i,r}^+,x_{j,s}^-] =
\frac{\delta_{ij}}{q_i-q_i^{-1}} 
\frac{1}{c_i^+c_j^-\hbar^2}\intCi{i}{+}\intCi{j}{-}u^rv^s g_i^+(u)^{-1}g_j^-(v)^{-1}\\
\cdot
\frac{\Gamma^+(u-v+a)}{\Gamma^+(u-v-a)}
\frac{\Psi_i\lp\Exp{v}\rp - e^{\pp (v-u)}\Psi_i\lp\Exp{u}\rp}
{1-e^{\pp (v-u)}}\, du\, dv
\end{multline*}
since the contribution of the boundary term $\Psi_{i,0}^-$ to the double
integral is zero given that the function $\Gamma^+(u-v+a)$ is holomorphic
in $D_i\times D_j$. We may therefore assume that $i=j$.

Note that
\[\frac{\Psi_i\lp\Exp{v}\rp - e^{\pp (v-u)}\Psi_i\lp\Exp{u}\rp}{1-e^{\pp (v-u)}}=
\frac{\Psi_i\lp\Exp{v}\rp-\Psi_i\lp\Exp{u}\rp}{1-e^{\pp (v-u)}}+\Psi_i\lp\Exp{u}\rp\]
and that the second summand on the \rhs does not contribute to the double
integral. Choosing the contours $\calC_i^\pm$ to be equal to a common
contour $\calC_i$ which contains all poles of $\xi_i(u)$, we may apply
Lemma \ref{le:double} to find
\[[x_{i,r}^+,x_{i,s}^-]
=
\frac{1}{q_i-q_i^{-1}}
\frac{1}{c_i^+c_i^-\hbar^2} 
\frac{\Gamma^+(\hbar d_i)}{\Gamma^-(\hbar d_i)}
\intCi{i}{}
u^{r+s}g_i^+(u)^{-1}g_i^-(u)^{-1}
\Psi_i\lp\Exp{u}\rp\, du
\]
where we used the fact that
\[\lim_{u\to v} \frac{u-v}{1-e^{\pp (v-u)}} = \frac{1}{2\pi\iota}\]
Using \eqref{eq:Si}, \eqref{eq: constants} and \eqref{eq:sine}, we then get
\[[x_{i,r}^+,x_{i,s}^-]
=\frac{1}{2\pi\iota\hbar} \intCi{i}{} u^{r+s} \xi_i(u)\, du\\
=\xi_{i,r+s}\]

\subsection{Proof that $\Gh{\Pi}\circ\Fh{\Pi}=\id$}\label{ssec: pf-inv}
%------------------------------------------------------------

Let $V\in \Rysub$. By \eqref{eq:fun X(z)}, $x_i^\pm(u)$ acts on $\Gh
{\Pi}\circ\Fh{\Pi}(V)$ as
\begin{multline*}
\frac{1}{c_i^{\pm}}\intCi{i}{\pm}\frac{1}{u-v} g_i^\pm(v)^{-1} \X^\pm_i(v)\,dv
\\
= \intCi{i}{\pm} \frac{1}{u-v} g_i^\pm(v)^{-1} \lp \oint_{{\calC'}^\pm_i}
\frac{1}{1-e^{-\pp (v-v')}} g_i^\pm(v')x_i^\pm(v')\, dv' \rp \, dv
\end{multline*}
where $u$ lies outside of $\calC^\pm_i$, and $\calC^\pm_i$ lies outside of
${\calC'}^\pm_i+\Z$. Assuming that $\calC^\pm_i$ contains ${\calC'}^\pm_i$ and
none of its $\nZ$--translates, and integrating in $v$ first
yields
\[\frac{1}{2\pi\iota}\int_{{\calC'}^\pm_i}\frac{1}{u-v'} x_i^\pm(v')\, dv'=x_i^\pm(u)\]
where the equality follows from the fact that ${\calC'}^\pm_i$ contains all
the poles of $x_i^\pm(u)$.

\subsection{Proof that $\Fh{\Pi}\circ\Gh{\Pi}=\id$}\label{ssec: pf'-inv}
%-----------------------------------------------------------

Let $\V\in\Rlsub$. Then, $\X_i^\pm(z)$ acts on $\Fh{\Pi}\circ\Gh
{\Pi}(\V)$ as
\begin{multline*}
c_i^{\pm}\ointCi{i}{\pm}\frac{z}{z-\Exp{v}}g_i^\pm(v)x_i^\pm(v)\,dv\\
=
\ointCi{i}{\pm}\frac{z}{z-\Exp{v}}g_i^\pm(v)
\int_{{\calC'}^\pm_i}\frac{1}{v-v'}g_i^\pm(v')^{-1}\X_i^\pm(\Exp{v'})\,dv'\,dv
\end{multline*}
where $z$ lies outside of $\exp(2\pi\iota\calC^\pm_i)$, and $\calC^\pm_i$ lies outside
of ${\calC'}^\pm_i$. Assuming that $\calC^\pm_i$ contains
${\calC'}^\pm_i$, and integrating in $v$ first yields
\[\int_{{\calC'}^\pm_i}\frac{z}{z-\Exp{v'}}\X_i^\pm(\Exp{v'})\,dv'=
\oint_{\wt{\calC'}^\pm_i}\frac{z}{z-w}\X_i^\pm(w)\,\frac{dw}{w}\]
where $\wt{\calC'}^\pm_i=\exp(2\pi\iota{\calC'}^\pm_i)$ contains all the poles of $\X_i^\pm$.
Since the one--form $dw/(z-w)w$ is regular on $\IP^1\setminus\{z,0\}$, and
$\X_i^\pm(0)=0$, the above is equal to $\X_i^\pm(z)$.

\section{$q$--Characters}
\label{sec: q-characters}
%==================

\subsection{$q$--characters for $\qloop$}\label{ssec: q-char-qla}
%--------------------------------------------------------------

The $q$--characters of \fd representations of the quantum loop
algebra $\qloop$ of a semisimple Lie algebra were introduced by
Frenkel--Reshetikhin in \cite{frenkel-reshetikhin-qchar}, and later
extended to the case of a symmetrisable \KM algebra by Hernandez
in \cite{hernandez-algebraic, hernandez-affinizations}. We follow
the presentation in \cite{hernandez-algebraic,hernandez-affinizations}.

Let $\V\in\Rloop$. Given a collection $\gamma=(\gamma_{i,\pm m}
^{\pm})_{i\in \bfI, m\in\N}$ of complex numbers, and $\mu\in\h^*$
such that $q_i^{\pm\mu(\alpha_i^{\vee})} = \gamma_{i,0}^{\pm}$, 
define $\V[\mu,\gamma]$ as follows
\begin{multline*}
\V[\mu,\gamma]=
\{v\in \V|\,K_hv = q^{\mu(h)}v, \text{ and }\\
\exists p>0 \text{ such that }(\Psi_{i,\pm m}^{\pm} - \gamma_{i,\pm m}^{\pm})^pv=0\,
\forall i\in\bfI, m\in\N\}
\end{multline*}

\begin{prop}\label{pr:q-char-qla1}
For a given pair $(\mu,\gamma)$, set
\[\gamma_i(z)^{\pm}=\sum_{m\geq 0} \gamma_{i,\pm m}^{\pm} z^{\mp m}\]
Then, there exists $\V\in\Rloop$ such that $\V[\mu,\gamma]\neq 0$ if, and only
if 
\begin{enumerate}
\item There exist monic polynomials $\{\QQ_i(w), \RR_i(w)\}$
such that $\QQ_i(0)\neq 0$,  $\RR_i(0)\neq 0$ and for every $i\in \bfI$
we have
\begin{equation}\label{eq: q-char-qla}
\gamma_i(z)^{\pm} = q_i^{-\deg(\QQ_i) + \deg(\RR_i)} 
\frac{\QQ_i(q_i^2z)}{\QQ_i(z)} \frac{\RR_i(z)}{\RR_i(q_i^{2}z)}
\end{equation}
\item $\mu\leq \lambda$ for some $\lambda\in P_+$.
\end{enumerate}
\end{prop}

\subsection{}\label{ssec: the-ring-qla}
%----------------

Consider the collection of variables $\{Y_{i,a}\}_{i\in\bfI, a\in\nC}$, and let 
$\mathcal{A}$ be the abelian group consisting of elements $e(\mu)\cdot M$, 
where $\mu\in\h^*$ and $M$ is a Laurent monomial in the variables $\{Y_{i,a}
\}$ such that
\[\mu(\alpha_i^{\vee}) = \sum_{a\in\nC} \text{degree of } Y_{i,a} \text{ in } M\]
The group operation on $\mathcal{A}$ is $e(\mu)\cdot M e(\mu')\cdot M' = e
(\mu+\mu')\cdot MM'$, with $e(0)\cdot 1$ being the unit element.

Let $\mathcal{Y}$ be the group algebra of $\mathcal{A}$ over $\Z$, completed 
in the following sense. An element of $\mathcal{Y}$ is a formal (possibly infinite)
linear combination $\chi$ of elements of $\mathcal{A}$ with $\Z$ coefficients, satisfying
the following condition: there exist $\lambda_1,\cdots, \lambda_r\in P_+$ such
that if $e(\mu)\cdot M$ appears in $\chi$ with non--zero coefficient, then $\mu
\leq \lambda_j$ for some $j\in\{1,\ldots,r\}$.

\subsection{}\label{ssec: q-char-qla2}
%--------------

To each $(\mu,\gamma)$ of the form given by Proposition \ref{pr:q-char-qla1},
we associate a monomial $M(\mu,\gamma)\in\mathcal{Y}$ as follows. Let $\{\QQ
_{i,\gamma}, \RR_{i,\gamma}\}$ be a collection of polynomials associated to
$\gamma$ so that \eqref{eq: q-char-qla} holds. If

\[\QQ_{i,\gamma} = \prod_j (w-a_j^{(i)}) \aand 
\RR_{i,\gamma} = \prod_k (w-b_k^{(i)})\]
set
\[
M(\mu,\gamma)=e(\mu)\cdot
\prod_{i\in\bfI} \lp\prod_j Y_{i,a_j^{(i)}}\rp\lp\prod_k Y_{i,b_k^{(i)}}\rp^{-1}
\]

The $q$--character of $\V\in\Rloop$ is defined by
\begin{equation}\label{eq: q-char-qla-2}
\qchar(\V) = \sum_{(\mu,\gamma)} \dim(\V[\mu,\gamma]) M(\mu,\gamma) \in \mathcal{Y}
\end{equation}
It gives rise to a homomorphism of abelian groups
\[\qchar : K\lp\Rloop\rp \to \mathcal{Y}\]

\begin{rem}
When $\g$ is a simple Lie algebra, the tensor structure on $\Rloop$
induces a ring structure on $K\lp\Rloop\rp$. In this case, it is proved
in \cite{frenkel-reshetikhin-qchar} that $\qchar$ is an injective ring
homomorphism.
\end{rem}

\subsection{$q$--characters for $\Yhg$}\label{ssec: q-char-y}
%------------------------------------------------

$q$--characters for the Yangian $\Yhg$ of a semisimple Lie algebra
were introduced by Knight in \cite{knight}. The treatment in \cite
{hernandez-algebraic,hernandez-affinizations} readily carries over
to the Yangian, and allows to extend their definition to the case of
a symmetrisable \KM algebra.

Given a collection $\beta = (\beta_{i,r})_{i\in\bfI,r\in \N}$ of complex
numbers, $\mu\in\h^*$ such that $d_i\mu(\alpha_i^{\vee}) = \beta_
{i,0}$, and $V\in\Ryang$, define
\begin{multline*}
V[\mu,\beta]=\{v\in V : h.v = \mu(h)v,\,\text{ and }\\
\exists p>0 \text{ such that } (\xi_{i,r}-\beta_{i,r})^pv=0\,
\forall\,i\in\bfI, r\in\N\}
\end{multline*}

\begin{thm}\label{thm: q-char-y}\cite{knight}
Given a pair $(\mu,\beta)$ as above, set
\[\beta_i(u) = 1 + \hbar\sum_{r\geq 0} \beta_{i,r}u^{-r-1}\in1+u^{-1}\C[[u^{-1}]]\]
Then, there exists $V\in\Ryang$ such that 
$V[\mu,\beta]\neq 0$ if, and only if 
\begin{enumerate}
\item There exist monic polynomials $\{Q_i, R_i\}_{i\in\bfI}$ such that for every
$i\in \bfI$ we have

\begin{equation}\label{eq: q-char-y}
\beta_i(u) = \frac{Q_i(u+\hbar d_i)}{Q_i(u)} \frac{R_i(u)}{R_i(u+\hbar d_i)}
\end{equation}

\item $\mu\leq \lambda$ for some $\lambda\in P_+$.
\end{enumerate}
\end{thm}

\subsection{}\label{ssec: the-ring-y}
%--------------

Define a commutative ring $\mathcal{X}$ analogous to $\mathcal{Y}$,
over the collection of variables $\{X_{i,a}\}_{i\in\bfI, a\in\C}$. Namely,
consider the abelian group of monomials $e(\mu)\cdot m$ as in \S
\ref{ssec: the-ring-qla}, except that $m$ is now a Laurent monomial
in $\{X_{i,a}\}$. Then, $\mathcal{X}$ is the group algebra of this abelian
group over $\Z$, completed in the same way as in \S \ref{ssec: the-ring-qla}.

\subsection{}\label{ssec: q-char-y2}
%--------------

To a pair $(\mu,\beta)$ of the form given by Theorem \ref{thm: q-char-y}, 
we associate a monomial $m(\mu,\beta)\in \mathcal{X}$ as follows. Let
$\{Q_i,R_i\}$ be a set of monic polynomials associated to $\beta$ so
that \eqref{eq: q-char-y} holds. If
\[Q_i = \prod_j (u-a_j^{(i)}) \aand R_i = \prod_k (u-b_k^{(i)})\]
set
\[m(\mu,\beta) = e(\mu)\cdot
\prod_{i\in\bfI} \lp\prod_j X_{i,a_j^{(i)}}\rp\lp\prod_k X_{i,b_k^{(i)}}\rp^{-1}\]

The $q$--character of $V\in\Ryang$ is now defined as
\begin{equation}\label{eq: q-char-y2}
\ychar(V) := \sum_{(\mu,\beta)} \dim(V[\mu,\beta]) m(\mu,\beta) \in \mathcal{X}
\end{equation}
and gives rise to a homomorphism of abelian groups
\[\ychar: K\lp\Ryang\rp \to \mathcal{X}\]
For a semisimple Lie algebra $\g$, $\ychar$ is an injective ring homomorphism
\cite{knight}.

\subsection{}\label{ssec: app-thm3}
%-------------

Let $\mathcal{X}_{\Pi} \subset \mathcal{X}$ be the subring generated by
$\h^*$ and the variables $\{X_{i,a}\}_{i\in\bfI, a\in\Pi}$. By Theorem \ref
{thm: first-main-theorem}, $\ychar$ restricts to a homomorphism
\[\ychar : K(\Rysub) \to \mathcal{X}_{\Pi}\]

Similarly, let $\mathcal{Y}_{\Omega}\subset \mathcal{Y}$ be the subring
generated by $\h^*$ and $\{Y_{i,\alpha}\}_{i\in \bfI, \alpha\in\Omega}$.
Again by Theorem \ref{thm: first-main-theorem}, $\qchar$ restricts to a
homomorphism
\[\qchar : K(\Rlsub) \to \mathcal{Y}_{\Omega}\]

Consider the isomorphism $e_{\Pi}:\mathcal{X}_{\Pi} \to \mathcal{Y}_
{\Omega}$ which is identity on $\h^*$ and maps $X_{i,a}$ to $Y_{i,e^
{\pp a}}$. The following result follows from the definition of the functor
$\Fh{}$ and the computation in \ref{ssec: diff-ex}.

\begin{prop}\label{prop: app-thm3}
The following diagram is commutative
\[\xymatrix{
K(\Rysub) \ar[rr]^{(\Fh{\Pi})_*} \ar[d]_{\ychar} && K\lp\Rlsub\rp \ar[d]^{\qchar} \\
\mathcal{X}_{\Pi} \ar[rr]_{e_{\Pi}} && \mathcal{Y}_{\Omega}
}\]
In particular, for a simple Lie algebra $\g$, $(\Fh{\Pi})_*$ is an isomorphism
of commutative rings.
\end{prop}

\end{document}